\documentclass[12pt]{amsart}

\usepackage{amssymb,amscd,amsthm, verbatim,amsmath, color,fancyhdr, mathrsfs}
\usepackage{graphicx}
\usepackage{turnstile}

\usepackage{hyperref}
\usepackage{cleveref}
\usepackage{rotating}
\usepackage{subcaption}
\usepackage{algorithm2e}
\usepackage{tikz-cd}

\usepackage[letterpaper, left=2.5cm, right=2.5cm, top=2.5cm,
bottom=2.5cm,dvips]{geometry}

\newtheorem{theorem}{Theorem}[section]

\newtheorem{cor}[theorem]{Corollary}

\newcommand{\bR}{\mathbb{R}}
\newcommand{\im}{\operatorname{im}}
\newcommand{\eps}{\epsilon}

\begin{document}

\author{Taejin Paik}
\address{Department of Mathematical Sciences and Research Institute of Mathematics, Seoul National University}
\email{paiktj@snu.ac.kr}

\author{Jaemin Park}
\address{Department of Mathematical Sciences and Research Institute of Mathematics, Seoul National University}
\email{woalsee@snu.ac.kr}

\title{Circular Coordinates for Density-Robust Analysis}

\date{}

\begin{abstract}
Dimensionality reduction is a crucial technique in data analysis, as it allows for the efficient visualization and understanding of high-dimensional datasets.
The circular coordinate is one of the topological data analysis techniques associated with dimensionality reduction but can be sensitive to variations in density.
To address this issue, we propose new circular coordinates to extract robust and density-independent features.
Our new methods generate a new coordinate system that depends on a shape of an underlying manifold preserving topological structures.
We demonstrate the effectiveness of our methods through extensive experiments on synthetic and real-world datasets.

\end{abstract}

\maketitle


\section{Introduction}



Dimensionality reduction allows us to understand high-dimensional data and gives us intuitive information about a dataset.
One of the key challenges in this area is preserving the intrinsic topological structure.
Different dimensionality reduction strategies try to handle this problem in different ways.

Principal component analysis (PCA) \cite{pearson1901liii,hotelling1933analysis} is one of the most basic techniques for linear dimensionality reduction.
Given a dataset, PCA aims to find a projection to a low-dimensional vector space maximizing variance.
As a non-linear dimensionality reduction, t-Distributed Stochastic Neighbor Embedding (t-SNE) \cite{van2008visualizing} is commonly used to visualize high-dimensional data.
It is a stochastic approach that aims to capture close points as close points in the low-dimensional embedding. 
Dimensionality reduction methods including the above dimensional reduction techniques, however, often fail to maintain the original topological structure when data points are sampled from an underlying manifold $M$ with complex topology.

The circular coordinate, which was introduced in \cite{de2009persistent}, deals in part with this problem by capturing $1$-dimensional holes; if there are $1$-dimensional holes in the underlying manifold $M$, the coordinates give circle-valued maps $\{\theta: X\rightarrow \mathbb{R}/\mathbb{Z}\}$ to identify the holes.
This approach is motivated by the bijection
$$\langle \mathcal{K}, K(\mathbb{Z}, 1)\rangle \cong H^1(\mathcal{K} ; \mathbb{Z})$$
for every CW-complex $\mathcal{K}$ where $\langle \mathcal{K}, K(\mathbb{Z}, 1)\rangle$ is basepoint-preserving homotopy classes of maps from $\mathcal{K}$ to the Eilenberg-MacLane space $K(\mathbb{Z}, 1)$, which is the circle $S^1$.
That is, for each cocycle in $H^1(X; \mathbb{Z})$, we can get a map $f: X\rightarrow S^1$.
Practically, for each cocycle $\alpha$, the circular coordinate is obtained by finding $L^2$-norm minimizer among cocycles that are cohomologous to $\alpha$.
We show how circular coordinates help to reveal topological structures hidden in low-dimensional embeddings in \Cref{subsec: exp-coil100}.

\begin{figure}[t]
\begin{subfigure}{.30\textwidth}
    \centering
    \includegraphics[width=\linewidth]{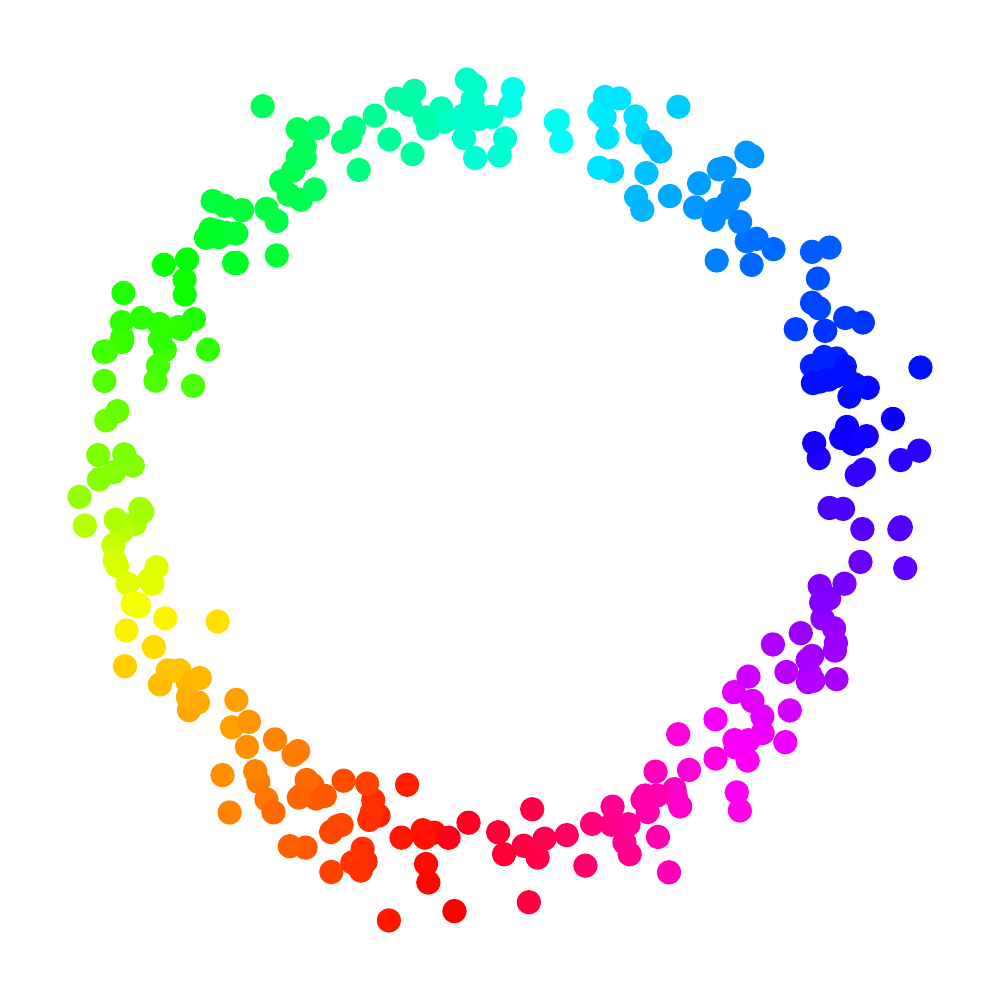}
\end{subfigure}
\hspace{2cm}
\begin{subfigure}{.30\textwidth}
    \centering
    \includegraphics[width=\linewidth]{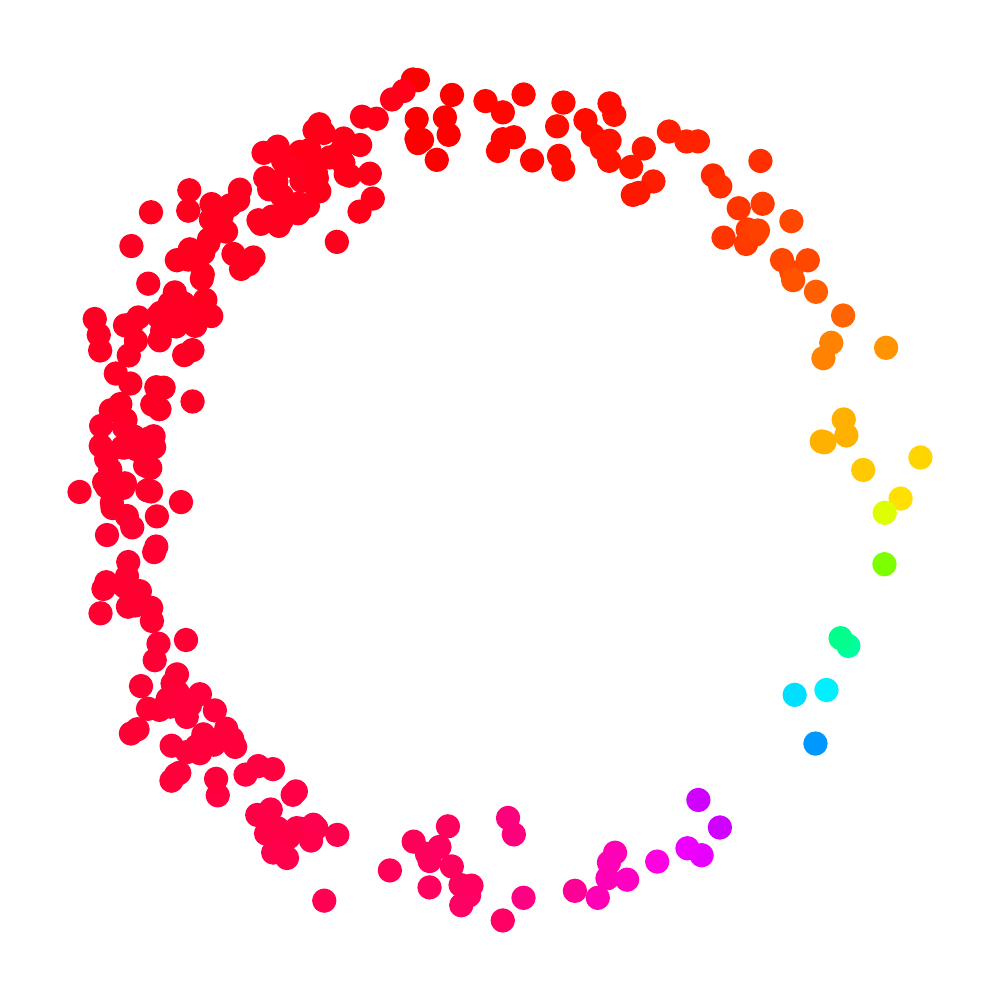}
\end{subfigure}
\caption{Changes in circular coordinates with probability density on $S^1$.}
\label{fig: intro circle example}
\end{figure}

It should be noted that the circular coordinates depend not only on the manifold's shape but also on the probability density function on the manifold.
For instance, as shown in \Cref{fig: intro circle example}, though we sample data points on the same circle, the circular coordinates are different depending on the probability density functions; in the low-density region, the circular coordinate changes quickly, whereas in the high-density region, it changes extremely slowly.
The results may be difficult to analyze if they vary with density if we want to explore the shape of manifolds embedded in Euclidean space.

In this research, we propose new circular coordinates that are dependent on the shape of the underlying Riemannian submanifold, rather than the probability density function on the manifold.
We propose two distinct approaches for achieving this goal:

\begin{enumerate}
    \item Obtaining a circular coordinate as a solution to a Dirichlet problem using a Laplacian matrix, which is obtained by approximating the Laplace-Beltrami operator of the underlying manifold.
    \item Utilizing the $L^p$-norm with $p>2$ instead of the $L^2$-norm in the optimization process to obtain a new circular coordinate.
\end{enumerate}
In the \Cref{sec: method}, we provide some justifications for why each approach generates a new circular coordinate that is robust to changes in the probability density function.
In the \Cref{sec: experiments}, we demonstrate the robustness of these methods through evaluations on synthetic datasets and a real-world dataset.
Our source is available under \url{https://github.com/TJPaik/CircularCoordinates}.

\section{Preliminaries}

In this section, we briefly review three well-known theories:
Hodge theory \cite{lim2020hodge}, persistent cohomology \cite{de2009persistent}, and the circular coordinate.

\subsection{Cohomology and Hodge theory} \label{subsection: Cohomology and Hodge theory}

Let $X$ be a finite simplicial complex.
The space of $i$-cochains $\mathcal{C}^i(X;\bR)$ is defined as the vector space dual to the vector space of $i$-chains $\mathcal{C}_i(X; \bR)$, and the coboundary map $d_i$ is dual to the boundary map between $\mathcal{C}_i(X;\bR)$ and $\mathcal{C}_{i+1}(X;\bR)$.
If the complex $X$ is clear from the context, we abbreviate $\mathcal{C}^i(X;\bR)$ by $\mathcal{C}^i$.
Specifically, the coboundary maps $d_i : \mathcal{C}^i \to \mathcal{C}^{i+1}$ for $i=0,1$ are
\begin{align*}
(d_0 f)(xy) &= f(y) - f(x), \text{ and}\\
(d_1 \alpha)(xyz) &= \alpha(xy) - \alpha(xz) + \alpha(yz).
\end{align*}
For $\alpha \in \mathcal{C}^1$, we call $\alpha$ a \textit{cocycle} if $d_1 \alpha = 0$, i.e. $\alpha \in \ker d_1$.
We call $\alpha$ a \textit{coboundary} if $\alpha = d_0 f$ for $f \in \mathcal{C}^0$, i.e. $\alpha \in \im d_0$.
Since $d_1 d_0 f = 0$ for all $f \in \mathcal{C}^0$, we have $\im d_0 \subset \ker d_1$.
We now define \textit{1-cohomology} of $X$ by
\[ H^1(X;\bR) = \ker d_1 / \im d_0. \]
We say $\alpha, \beta \in \mathcal{C}^1$ are \textit{cohomologous} if $[\alpha] = [\beta] \in H^1(X;\bR)$.

In this paper, we use the standard basis on $\mathcal{C}^i$ and the dual basis of the standard basis on $(\mathcal{C}^i)^*$.
We assign an inner product on $\mathcal{C}^i$ for each $i$.
Whenever an inner product is not specified, it will be assumed to be the standard inner product.
It follows from linear algebra that $\mathcal{C}^1$ can be decomposed (called the \textit{Fredholm alternative}) as
\begin{equation} \label{Fredholm alternative}
\mathcal{C}^1 \cong \ker d_1 \oplus \im d_1^*,
\end{equation}
where $d_i^*$ is the adjoint operator of $d_i$ for $i=1,2$.
Note that the adjoint operator depends on inner products.
Moreover, there is the \textit{Hodge decomposition} 
\begin{equation} \label{Hodge decomposition}
    \mathcal{C}^1 \cong \im d_1^* \oplus \ker ( d_1^* d_1 + d_0 d_0^* ) \oplus \im d_0.
\end{equation}
Here, we call $d_1^* d_1 + d_0 d_0^*$ the \textit{(1-dimensional) Hodge Laplacian} and denote it by $\Delta_1$.
Combining two decompositions \eqref{Fredholm alternative} and \eqref{Hodge decomposition}, we have
\[ \ker ( d_1^* d_1 + d_0 d_0^* ) \oplus \im d_0 = \ker d_1. \]
Thus, we have
\[ H^1(X;\bR) = \ker d_1 / \im d_0 \cong \ker ( d_1^* d_1 + d_0 d_0^* ). \]
We note that each $\alpha \in \ker \Delta_1$ is the representative of the corresponding equivalent class in $\ker d_1 / \im d_0$ with the minimal $2$-norm.
That is, for every $[\alpha] \in \ker d_1 / \im d_0$, the representative 
\[ \alpha_H = \underset{\overline{\alpha}}{\operatorname{argmin}}\left\{\|\overline{\alpha}\|_2 \mid \exists f \in \mathcal{C}^0, \overline{\alpha}=\alpha+d_0 f\right\}  \]
is contained in $\ker \Delta_1$ (see \cite{lim2020hodge}) since $\alpha_H \perp \im d_0$.
We call $\alpha_H$ the \textit{harmonic cocycle}.

\subsection{Persistent cohomology}
Now we consider a dataset $X$ in the Euclidean space and a $1$-parameter family of Vietoris-Rips complexes $\{ X^\eps \}$ where $\eps$ is a scale parameter with the dataset.
Let $\eps_1,\dots, \eps_m$ be critical values where the homotopy type of $X^\eps$ changes.
We can write this situation as follows:
\begin{center}
    \begin{tikzcd}
        X^{\epsilon_1} \arrow[r, "i_1"] & X^{\epsilon_2} \arrow[r, "i_2"] & \cdots \arrow[r, "i_{m-1}"] & X^{\epsilon_m}
    \end{tikzcd}
\end{center}
where $\to$ denotes the inclusion maps.
The inclusion maps between $X^\eps$'s induce the homomorphisms between $H^1(X^\eps;\bR)$:
\begin{center}
    \begin{tikzcd}
        H^1(X^{\epsilon_1};\mathbb{R}) & H^1(X^{\epsilon_2};\mathbb{R}) \arrow[l, "i_1^*"'] & \cdots \arrow[l, "i_2^*"'] & H^1(X^{\epsilon_m};\mathbb{R}) \arrow[l, "i_{m-1}^*"']
    \end{tikzcd}
\end{center}
For a nonzero cocycle class $[\alpha] \in H^1(X^{\eps_k};\bR)$, let 
\begin{align*}
    b_\alpha &= \inf \left\{ \eps \in \{\eps_1, \dots, \eps_k\} : \exists [\beta] \neq 0 \in H^1(X^\eps;\bR) \text{ such that } i^*[\beta] = [\alpha] \right\} \\
    d_\alpha &= \sup \left\{ \eps \in \{\eps_k,\dots,\eps_m\} : i^*[\alpha] \neq 0 \in H^1(X^\eps;\bR) \right\},
\end{align*}
where $i^*$ is an induced homomorphism by an inclusion map.
We call $b_\alpha$ ($d_\alpha$, respectively) a \textit{birth} (\textit{death}, respectively) and the value $d_\alpha - b_\alpha$ \textit{life time}.
If we collect all $b_\alpha$ and $d_\alpha$ for all nonzero cocycle and draw points $\{(b_\alpha, d_\alpha)\}_\alpha$ on a coordinate plane, that is called \textit{persistence diagram}.


\subsection{Circular coordinate}\label{subsec: prelim Circular coord}
The circular coordinate is introduced in \cite{de2009persistent}.
Given a dataset $X$, its circular coordinate $\theta : X \to \mathbb{R}/\mathbb{Z}\cong S^1$ is defined using a nonzero cocycle class $[\alpha] \in H^1(X^\epsilon;\mathbb{Z}_\mathbf{p})$ with a fixed $\epsilon$ and a prime $\mathbf{p}$ if $[\alpha]$ lies in the image of the coefficient homomorphism $H^1(X^\epsilon;\mathbb{Z}) \to H^1(X^\epsilon;\mathbb{Z}_\mathbf{p})$.
Practically, we take a large prime $\mathbf{p}$.
The brief algorithm is as follows.
First, we find the harmonic cocycle $\alpha_H$ which is cohomologous to $\alpha$, i.e.
\begin{equation} \label{minimal norm cocycle}
    \alpha_H = \underset{\overline{\alpha}}{\operatorname{argmin}}\left\{\|\overline{\alpha}\|_2 \mid \exists f \in \mathcal{C}^0(X^\eps; \mathbb{R}), \overline{\alpha}=\alpha+d_0 f\right\}.
\end{equation}
Then we fix a vertex $x$ and assign $\theta (x) = 0$.
For a vertex $y$, we assign $\theta(y) = \sum_{i=1}^n \alpha_H(e_i)$, where $e_1 \cdots e_n$ is an edge path starting from $x$ to $y$.
The circular coordinate $\theta$ is well-defined since $\alpha_H$ is cohomologous to a cocycle with integer coefficients.
We note that we can use $f$ in \eqref{minimal norm cocycle} as the circular coordinate.

Now, we interpret the circular coordinate in terms of Laplacian.
Since $\alpha_H \in \im d_0^\perp = \ker d_0^*$, we have
\[ 0 = d_0^* \alpha_H = d_0^* \alpha + d_0^* d_0 f. \]
Here, we call $d_0^* d_0$ the \textit{(0-dimensional) graph Laplacian} and denote it by $\Delta_0$.
Thus the circular coordinate $f$ is a solution of the Dirichlet problem
\begin{equation} \label{Dirichlet problem}
    \Delta_0 f = - d_0^* \alpha.
\end{equation}
The solution of \eqref{Dirichlet problem} exists and is unique.
More formally,

\begin{theorem}\label{thm: unique and exist of the solution}
    Let $\alpha\in \mathcal{C}^1$ be a cocycle.
    Then, a solution to the Dirichlet problem on the graph(1-skeleton of the Vietoris-Rips complex)
    $$
    \Delta_0 f = d_0^*\alpha
    $$
    exists, and if the graph is connected, then the solution is unique up to scalar addition.
\end{theorem}
See \Cref{sec: proofOfTheorems} for the proof.

\section{Methods}\label{sec: method}

\subsection{Weighted circular coordinate}\label{subsec: weighted circular}

A \textit{weighted circular coordinate} is a generalization of the circular coordinate that can take into account a density of a dataset.
By assigning weights to each edge in a Vietoris-Rips complex, we can get a new weighted circular coordinate that considers a density of a point cloud in the following way.

Let us assume that we have a simplicial complex $X$ and positive weights $w(e)$ on each edge $e$ of $X$.
We define a weight matrix $W$ by a diagonal matrix whose diagonal entries are the weights.
Formally, we write
\[ W = \mathrm{diag } \, \{ w(e) \; : \; e \text{ is an edge of } X \}. \]

We consider a \textit{weighted cochain complex}
\begin{equation}\label{eq: chain_complex_weighted}
\begin{tikzcd}
0 \arrow[r, "0"] & \mathcal{C}^0 \arrow[r, "Wd_0"] & \mathcal{C}^1 \arrow[r, "d_1W^{-1}"] & \mathcal{C}^2
\end{tikzcd}
\end{equation}
where $Wd_0$ and $d_1W^{-1}$ are defined as follows:
\begin{align*}
    (Wd_0) f (x y) &= w(x y) ( f(y) - f(x) ), && \forall f \in \mathcal{C}^0; \\
    (d_1 W^{-1}) \alpha (xyz) &= \frac{1}{w(xy)} \alpha(xy) - \frac{1}{w(xz)} \alpha(xz) + \frac{1}{w(yz)} \alpha(yz), && \forall \alpha \in \mathcal{C}^1.
\end{align*}
The \textit{weighted Hodge Laplacian} $\Delta_{W,1}$ is defined by
\[ \Delta_{W,1} = (d_1 W^{-1})^* (d_1 W^{-1}) + (W d_0) (W d_0)^*. \]
Note that the weighted harmonic cocycle $\alpha_{W, H} \in \ker \Delta_{W,1}$ is the minimal 2-norm cocycle among its cohomology class.
With this observation, we introduce a weighted circular coordinate.

Since $W\alpha \in \ker(d_1W^{-1})$ for every $\alpha \in \ker d_1$, we can think of the harmonic element for $W\alpha$ on this chain complex and $f' \in \mathcal{C}^0$ satisfying 
$$\|W\alpha + Wd_0f'\|_2 = \|W(\alpha + d_0f')\|_2 \leq  \|W(\alpha + d_0g)\|_2 $$
for every $g\in \mathcal{C}^0$ as we do above.

We can understand this in terms of an inner product in the following way.
Assume that we have a different inner product $\langle\cdot,\cdot \rangle_q$, instead of the standard inner product on $\mathcal{C}^1$.
Then, we can represent the inner product as a matrix from
$$
\langle v, w\rangle_q = v^t Q w
$$
where $Q$ is a real symmetric positive definite matrix.
From the Cholesky decomposition, we have $Q = W_0^tW_0$ for a real upper triangular matrix $W_0$.
Therefore, we have
\begin{align*}
\langle\alpha + d_0f'',\alpha + d_0f'' \rangle_q &= (\alpha + d_0f'')^tW_0^tW_0(\alpha + d_0f'')\\
&= (W_0\alpha + W_0d_0f'')^t(W_0\alpha + W_0d_0f'')\\
&= \|W_0\alpha + W_0d_0f''\|_2^2
\end{align*}
for $f''\in \mathcal{C}^0$, and the corresponding harmonic element using the inner product $\langle\cdot,\cdot\rangle_q$ is the same as the harmonic element on the chain complex \eqref{eq: chain_complex_weighted} if $W=W_0$.

For a dataset, we seek to obtain a weighted circular coordinate that is robust to the density of the dataset.
To achieve our goal of obtaining a circular coordinate that depends only on the shape of the dataset, we adjust the weights of the edges in the graph constructed from a Vietoris-Rips complex.
Now, we introduce how to set weights for each edge by understanding the graph Laplacian as an approximation of a Laplace-Beltrami operator $\Delta_M$ of a manifold $M$.

Data points in Euclidean space can be thought of as samples from a probability density function on a manifold.
In other words, the data points are drawn independently and identically from a probability distribution defined on a submanifold of Euclidean space.
Vietoris-Rips complex is a topological space constructed by connecting nearby data points with edges and higher simplices.
It can be used to approximate the topology of the underlying manifold.

The key idea is that a weighted graph Laplacian that is defined on the $1$-skeleton of the Vietoris-Rips complex approximates the Laplace-Beltrami operator of the underlying manifold.
If we approximate the Laplace-Beltrami operator and find a harmonic solution on the underlying manifold, it will be a solution that depends only on the shape of the manifold regardless of the density of the data.
Therefore, we need to approximate the Laplace-Beltrami operator with the given data.
Before proceeding, we introduce some notations and give some explanations.

Assume that we have $n$ points $x_1, \dots, x_n \in M$ where $M$ is a $k$-dimensional Riemannian submanifold of $\mathbb{R}^m$ and a probability density function $P:M\rightarrow\mathbb{R}$ with $\inf_{x\in M} P(x) > 0$.
We take a Vietoris-Rips complex constructed from the $n$ points with a scale parameter $\epsilon> 0$ and the $1$-skeleton of the Vietoris-Rips complex which is a graph.
In our setting, we consider each edge of the graph to have a direction so that we can give a different weight for each direction.

Given the $n$ points and $t>0$, we construct a weighted directed graph taking the weight of the edge connecting $x_i$ and $x_j$ to be $w_{ij} = \frac{1}{P(x_j)}g_{ij}$ where $g_{ij} = \frac{1}{(4\pi t)^{k/2}}e^{-\frac{\|x_i - x_j \|^2}{4t}}$.
We set $w_{ij} = 0$ if $x_i$ and $x_j$ are not connected.
Then the corresponding \textit{weighted directed graph Laplacian matrix} is defined as
$$
\left(L_n^t\right)_{ij}= \begin{cases}
-w_{ij} & \text { if } i\neq j \\
\sum_{k\neq i} w_{i k} & \text { if } i=j.
\end{cases}
$$
Note that the matrix $L_n^t$ can be understood as an operator on functions defined on the vertices:
$$
L_n^t f\left(x_i\right)=\sum_j w_{ij} (f(x_i) - f(x_j)),
$$
and naturally, on functions defined on the manifold $M$:
$$
\mathbf{L}_n^t f\left(x\right)=\sum_j w_t(x, x_j) (f(x) - f(x_j))
$$
where $w_t(x, y) = \frac{1}{P(y)}g_t(x, y)$ and $g_t(x, y) = \frac{1}{(4 \pi t)^{k / 2}} e^{-\frac{\|x-y\|^2}{4t}}$ if $\| x - y \| < \epsilon$ and $g_t(x, y) = 0$ for $\|x -y\| \geq \epsilon$ defined on $M\times M$.
We can easily see that $L_n^tf(x_i) = \mathbf{L}_n^tf(x_i)$ for every $i=1, \dots, n$.

Another note is that the matrix $L_n^t$ can be decomposed as $L_n^t = P(P^{-1}D - P^{-1}GP^{-1})$ where $P$ is a diagonal matrix with $P_{ii} = P(x_i)$, $D$ is a diagonal matrix with $(D)_{ii} = \sum_j g_{ij}/P(x_j)$, and $(G)_{ij} = g_{ij}$.
We denote the matrix $P^{-1}D - P^{-1}GP^{-1}$ by $L_{n, W}^t$.
Now, we can approximate the Laplace-Beltrami operator of the underlying manifold:
\begin{theorem}\label{thm: LB approxi}
Let $M$ be a $k$-dimensional compact Riemannian submanifold in $\mathbb{R}^n$ and $P:M\rightarrow \mathbb{R}$ be a smooth probability distribution function on $M$ with $\inf_{x\in M} P(x) > 0$.
If data points $x_1, \dots, x_n$ are independent and identically distributed samples drawn from the distribution $P$ and $f:M\rightarrow \mathbb{R}$ is a smooth function, then for $x\in M$ and $t_n=n^{-\frac{1}{k+2+\alpha}}$, where $\alpha>0$, we have
$$
\frac{1}{nt_n}\mathbf{L}_n^{t_n} f(x) \overset{p}{\rightarrow} \Delta_M f(x)
$$
as $n$ goes to infinity.
\end{theorem}
For proof, see \Cref{sec: proofOfTheorems}. 
From the above theorem, we can think of $\frac{1}{nt_n}L_n^t$ as a discrete approximation of $\Delta_M$ for appropriate $t$.
Since $\Delta_M$ is the Laplace-Beltrami operator of the underlying Riemannian manifold regardless of the density of the data, if we use $L_n^t$ to solve a Dirichlet problem for a cocycle, we find a circular coordinate that relies more on the shape of the manifold and robust to the density distribution.

We provide an interpretation of the matrix $L_n^t$ in the theorem that follows.
Before that, we define a diagonal matrix $Q_1$ where $(Q_1)_{ii} = (L_{n, W}^t)_{jk}$ where $i$ is the index of the edge connecting $x_j$ and $x_k$.
\begin{theorem}\label{thm: PL_nwt is dtd}
Let us assign an inner product on $\mathcal{C}^0$ by $\langle v, w\rangle = v^tP^{-1}w$, and on $\mathcal{C}^1$ by $\langle v, w\rangle = v^tQ_1w$.
Then the graph Laplacian $d_0^*d_0$ is $L_n^t$.
\end{theorem}
See \Cref{sec: proofOfTheorems} for the proof.
Since $L_n^t$ can be understood as a graph Laplacian, we can make and solve a Dirichlet problem from \Cref{thm: unique and exist of the solution} using a cocycle as follows:
\begin{cor}
Let $\alpha$ be a cocycle in $\mathcal{C}^1$ and take the inner products defined on \Cref{thm: PL_nwt is dtd}.
A solution to the Dirichlet problem on the graph
$$
L_n^tf = d_0^*\alpha
$$
exists, and if the graph is connected, then the solution is unique up to scalar addition.
\end{cor}

Practically, we approximate $P(x_i)$ by $\frac{1}{n}\sum_j g_{ij}$.
Indeed, $\frac{1}{n}\sum_j g_{ij}$ converges to 
$$
d_t(x_i) := \int_M g_t(x_i, x)P(x) \; dV(x)
$$
in probability as $n\rightarrow \infty$ where $dV$ is the volume form of $M$ induced from the ambient Euclidean space satisfying $\int_M P(x) \;dV(x)=1$, and $d_t(x_i)$ converges to $P(x_i)$ as $t$ goes to $0$.
Specifically, we denote the weighted circular coordinate using the weighted directed graph Laplacian matrix by \textit{WDGL-circular coordinate}.
For each cocycle, we can obtain the WDGL-circular coordinate equivalently by finding the harmonic element corresponding to the cocycle using the inner product $\langle v,w\rangle = v^tQ_1w$ on $\mathcal{C}^1$ from the relation between the Dirichlet problem and circular coordinates.

We also experiment with several other weight candidates.
The weights on edges can be chosen based on various criteria, such as the distance between the points connected by the edge or the number of points in the region surrounding the edge.
In our experiments, we show that weighted circular coordinates give a density-robust result when we use the following weights:
\begin{enumerate}
    \item The weight of an edge is the $\frac{1}{D_0 + D_1}$ where $D_0$ and $D_1$ are the degrees of vertices attached to the edge.
    \item The weight of an edge is the $\frac{1}{\sqrt{D_0D_1}}$ where $D_0$ and $D_1$ are the degrees of vertices attached to the edge.
\end{enumerate}
We present results using the above two weights in figures in \Cref{sec: app_exp_results}, and use notations ``1/(D0 + D1)'' and ``1/sqrt(D0 D1)'' to denote the above two weights respectively in the figures.

\subsection{$L^p$-circular coordinate}
A circular coordinate can be considered in terms of a harmonic cocycle representing a cohomology class \cite{de2009persistent}.
Let us denote $L^p$-norm by $\|\cdot\|_p$.
Given a simplicial complex $X$ and a cochain $\alpha \in \mathcal{C}^1$, the corresponding harmonic cocycle is 
$$
\underset{\overline{\alpha}}{\operatorname{argmin}}\left\{\|\overline{\alpha}\|_2 \mid \exists f \in \mathcal{C}^0, \overline{\alpha}=\alpha+d_0 f\right\},
$$
the minimal cocycle among all cohomologous cocycle with $\alpha$ concerning the $L^2$-norm.
The function $f\in\mathcal{C}^0$ that attains the minimum is the circular coordinate.
One of the variations we can consider is a minimizer with $L^p$-norm with different $p$, $L^p$-\textit{circular coordinate}.
That is, we can look at 
\begin{equation}
\label{equ: p_norm_argmin}
\underset{\overline{\alpha}}{\operatorname{argmin}}\left\{\|\overline{\alpha}\|_p \mid \exists f \in \mathcal{C}^0, \overline{\alpha}=\alpha+d_0 f\right\}.
\end{equation}
Previous research \cite{luo2021generalized} shows that a circular coordinate with $L^1$-norm or mixed norm becomes more ``locally constant''.
In this work, we study the behavior of a circular coordinate with different $p$.
Our several experiments show that the higher the $p$ is, the more robust the circular coordinate to the density.

To get the corresponding element $f\in\mathcal{C}^0$ for $p=1$, we can consider the problem 
as a linear programming problem and solve it in linear time \cite{jiang2011statistical,luo2021generalized}.
For $p=2$, the problem \eqref{equ: p_norm_argmin} is the well-known least-squares problem, and LSQR \cite{paige1982lsqr} can be used as has been observed in \cite{de2009persistent}.
In the case of a general $p>2$, we can use the gradient descent algorithm.
We specify our computation in Algorithm \ref{alg: p_circular}.

Now we introduce an iterative way to get a $L^\infty$-circular coordinate.
Note that $\|\overline{\alpha}\|_\infty= \underset{e}{\max}\{|\overline{\alpha}(e)|\}$.
We have the following theorem:
\begin{theorem}\label{thm: lp_converge}
    Suppose $f$ is a complex measurable function on $X$ and $\mu$ is a positive measure on $X$.
    Assume that $\|f\|_r <\infty$ for some $r<\infty$. Then 
    $$\|f\|_p\rightarrow \|f\|_\infty\quad \text{as }\;p\rightarrow \infty.$$
\end{theorem}
For proof, see \Cref{sec: proofOfTheorems}.
From the above theorem, a $L^p$-circular coordinate with large $p$ can be an approximation of a $L^\infty$-circular coordinate.
From this observation, we propose an iterative way to get a $L^\infty$-circular coordinate.
That is, we can obtain a $L^\infty$-circular coordinate by sequentially optimizing $L^p$-circular coordinate while increasing $p$ in some range, and finally, we optimize the function $f$ to get the $L^\infty$-circular coordinate.
We write down our algorithm in Algorithm \ref{alg: infty_circular_p}.
In practice, we increase $p$ up to about $50$ to obtain an estimated loss because of underflow or overflow issues.
We confirm that for many hyperparameter settings, this algorithm gives faster convergence to the $L^\infty$-circular coordinate.

As a similar methodology, we can use the softmax function with temperatures:
\begin{theorem}\label{thm: softmax}
    Let $h:\bR^n\rightarrow \bR^n$ be a function satisfying $(h(x))_i = |x_i|$, and $s$ be the softmax function.
    Then for $v\in \bR^n$, $(s\circ h)(tv) \cdot h(v)$ converges to $\|v\|_\infty$ as $t$ goes to infinity.
\end{theorem}
See \Cref{sec: proofOfTheorems} for the proof.
Similarly, after initializing a vector $f\in\bR^n$, we can optimize $(s\circ h)(t f) \cdot h(f)$ while increasing the parameter $t$.
We specify our method in Algorithm \ref{alg: infty_circular_softmax}.
In our experiment, we do not benefit much when using the softmax function; the convergence is usually not faster than in other methods.

In our algorithms, we use a hyperparameter $\tau$ to determine whether the function $f$ is converged or not.
If the difference between the previous loss and the current loss is smaller than $\tau$, we determine that the function is converged.

\subsubsection{$p$-harmonic on manifolds}
In order to gain a deeper understanding of the concept of the $L^p$-circular coordinate, let us examine an example of the analog on manifolds.
We consider a $n$-dimensional compact Riemannian manifold $(M,g)$, so we can define the usual smooth Hodge theory.
Fixing a cohomology class $[\alpha] \in H_{dR}^k(M)$ and a representative $\alpha$, we want to choose a representative $\alpha_H$ such that its $p$-norm with $p>1$ rather than its $2$-norm is minimal:
$$
\alpha_H := \underset{\overline{\alpha}}{\operatorname{argmin}}\left\{
    \int_M \| \overline{\alpha}\|^p \; dV \mid \exists f\in \Omega^{k-1}(M),  \overline{\alpha}=\alpha+df
    \right\}
$$
where $dV$ is a volume form.
Therefore, we have 
$$
\frac{d}{dt} \int_M \Vert \alpha_H+t df \Vert^p \; dV \;\bigg|_{t=0} =0.
$$
for every $f\in \Omega^{k-1}(M)$, and from the Euler-Lagrange equation, we have
\begin{equation}\label{eq: k-form harmonic}
    d\alpha_H = 0 \quad \text{and} \quad \delta(\|\alpha_H\|^{p-2} \alpha_H) = 0
\end{equation}
since
$$
\int_M \langle\alpha_H, df\rangle \|\alpha_H\|^{p-2} \; dV = 
\int_M \langle\delta(\|\alpha_H\|^{p-2}\alpha_H), f\rangle \; dV
$$
where $\delta$ is the codifferential operator and $\|\alpha\|^2$ is defined as $\sum_{i, j} g^{ij}\alpha_i \alpha_j$ for a $1$-form $\alpha =\sum_{i=1}^n\alpha_i\; dx^i$.
An $k$-form $\alpha_H$ satisfying \eqref{eq: k-form harmonic} is called \textit{$p$-harmonic $k$-form}.
We note that this definition is reduced to the $p$-laplacian in the case of Euclidean space.
We will not investigate the regularity of the above equation and only look at a simple toy example to see what the effect of $p$ on a minimizer is.
As a manifold, we take $M=S^1 \times D^{n-1}$, and we take $x^1$ to be the coordinate on the circle.
For $\rho > 0$ on $M$, define the metric 
$$
g= \rho^{2/n}\sum_{j=1}^n dx^j\otimes dx^j,
$$
which is simply a conformal rescaling of the standard Euclidean metric.
We understand the change in metric as the change in probability density function on the manifold.
The resulting density is $\sqrt{\det g}=\rho$.
For simplicity, we assume that the density only depends on $x^1$.

Now let's see what $p$-harmonic $1$-form $\alpha$ generating $H^1_{dR}(M)$ looks like.
We have $\alpha_H=\sum_i \alpha_i dx^i$.
To write out the $p$-harmonic equations, we need the Hodge star $\star$ to write out $\delta=(-1)^{n(k-1)+1} \star d \star$ for $k$-forms. 
Since we have
$$
\star dx^i=\sum_{j=1}^n (-1)^{j+1} g^{ij} \sqrt{\det g} \;dx^1\wedge \ldots \wedge\widehat{dx^j}\wedge \ldots \wedge dx^n,
$$
the equation $\delta(\|\alpha_H\|^{p-2}\alpha_H) = 0$ is equivalent to
$$
d\left(\|\alpha_H\|^{p-2}\sum_{i, j = 1}^n (-1)^{j+1}\alpha_i g^{ij}\sqrt{\det g}\;dx^1\wedge \ldots \wedge\widehat{dx^j}\wedge \ldots \wedge dx^n \right) = 0.
$$
If we use the form of our metric, this simplifies to
$$
d\left(
    \rho^{\frac{n-p}{n}}\sum_{i=1}^n |\alpha_i|^{p-2}\sum_{j=1}^n(-1)^{j+1} \alpha_j \; dx^1\wedge \ldots \wedge \widehat{dx^j}\wedge \ldots \wedge dx^n
\right) = 0.
$$
We find that if we set
$$
\alpha_H=\rho^{\frac{p-n}{n(p-1)}}(x^1) dx^1,
$$
then it satisfies \eqref{eq: k-form harmonic} since
$$
\rho^{\frac{n-p}{n}}\left(\sum_{i=1}^n |\alpha_i|^{p-2}\right)\left(\sum_{j=1}^n(-1)^{j+1} \alpha_j\right) = \rho^{\frac{n-p}{n}}(\rho^{\frac{p-n}{n(p-1)}})^{p-2} \rho^{\frac{p-n}{n(p-1)}} = 1.
$$

For $n=1$, the $p$-harmonic representative is just $\rho(x^1)dx^1$, so independent of $p$, and has linear dependence on the density induced by the Riemannian metric. 
For fixed large $n$, the $p$-harmonic representative has density dependence of the form $\rho^{1/n}$, which means that it is almost independent of the density for large $n$.
On the other hand, the usual harmonic representative has density dependence of the form $\rho^{2/n - 1}$, which is very different. 
We don't know how well this toy example generalizes, but the toy example suggests that the $p$-harmonic representative of a cohomology class has less density dependence than the usual harmonic representative for large $p$, provided the dimension $n$ is also large.
In the discrete setting, we may think of the Vietoris-Rips complex as a high-dimensional manifold if the scale parameter is sufficiently large.

\section{Experiments}\label{sec: experiments}

In this section, we test our methods on several synthetic data including a noisy circle, a noisy trefoil knot, two conjoined circles, and a torus dataset.
For real data analysis, we present experimental results on the COIL-100 dataset \cite{nene1996columbia} with low-dimensional embeddings.

For the synthetic datasets, we use correlation scatter plots to show the results.
A correlation scatter plot is a graphical representation of a relationship between two circular coordinates.
It consists of a set of points plotted on a coordinate plane, with one circular coordinate plotted on the $x$-axis and the other circular coordinate plotted on the $y$-axis.
The scatter plot is used to compare two circular coordinates.
In our experimental results, most scatter plots are made to show the comparison between the actual circular coordinate and the circular coordinate inferred from our algorithms.
Each circular coordinate is placed on the original dataset using a color map to represent the results.
In our figures, we use the cyclic HSV color wheel.

\subsection{Experimental details}
We have options to use various software to implement parts of our algorithms, and Ripser \cite{Bauer_2021} and SciPy \cite{2020SciPy-NMeth} are specifically used to get a cocycle and for the LSQR algorithm.
For optimizing $L^p$-norm, we use PyTorch \cite{paszke2019pytorch}.
To visualize our experimental results, we use the t-SNE of Scikit-learn \cite{scikit-learn} and matplotlib \cite{Hunter:2007} library.

For a non-zero cocycle $\alpha$, we choose $\epsilon = \frac{b_\alpha + d_\alpha}{2}$ where $b_\alpha$ and $d_\alpha$ are the birth and death of $\alpha$ respectively.
To get a persistence diagram, we use a prime of $47$.
If we use the WDGL method, we heuristically set $t$ to 0.2 times the average Euclidean distance between connected vertices.

\subsection{Noisy circle}
\begin{figure}[t]
\begin{subfigure}{.328\textwidth}
    \centering
    \includegraphics[width=\linewidth]{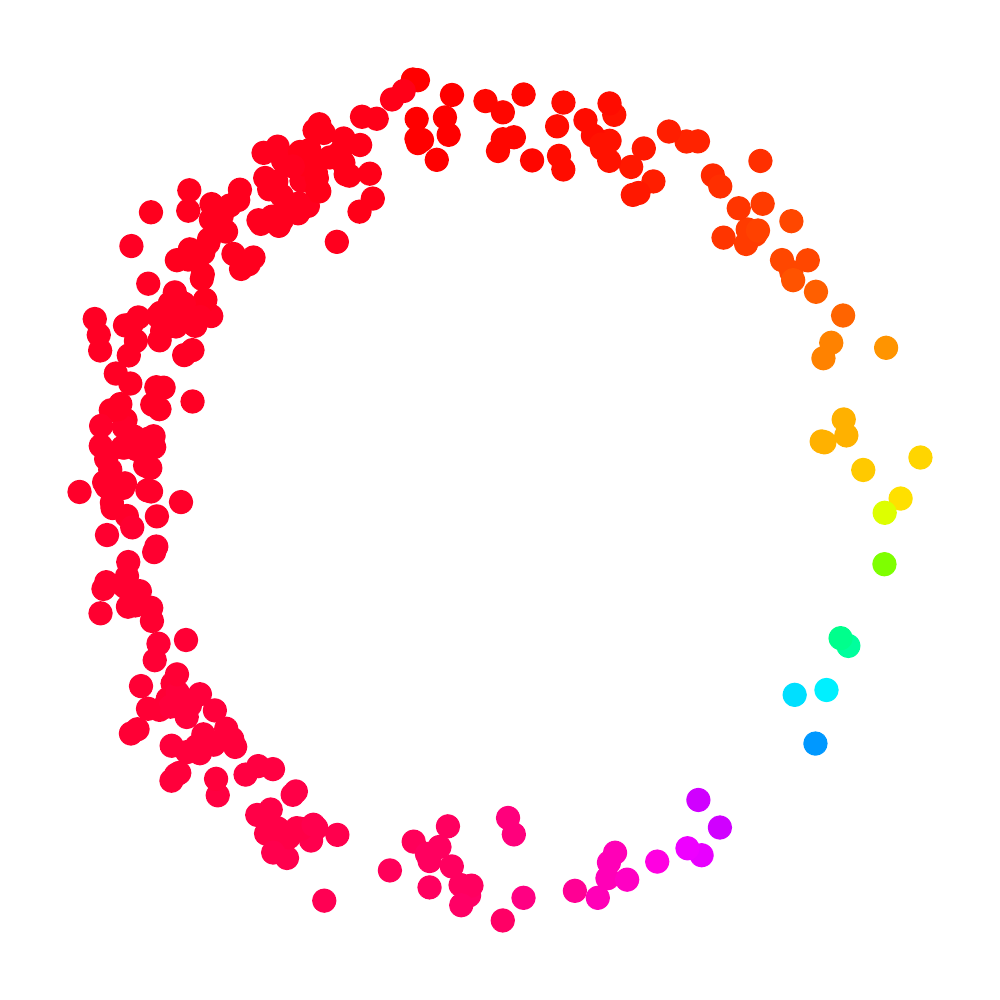}
    \caption{Circular coordinate}
\end{subfigure}
\begin{subfigure}{.328\textwidth}
    \centering
    \includegraphics[width=\linewidth]{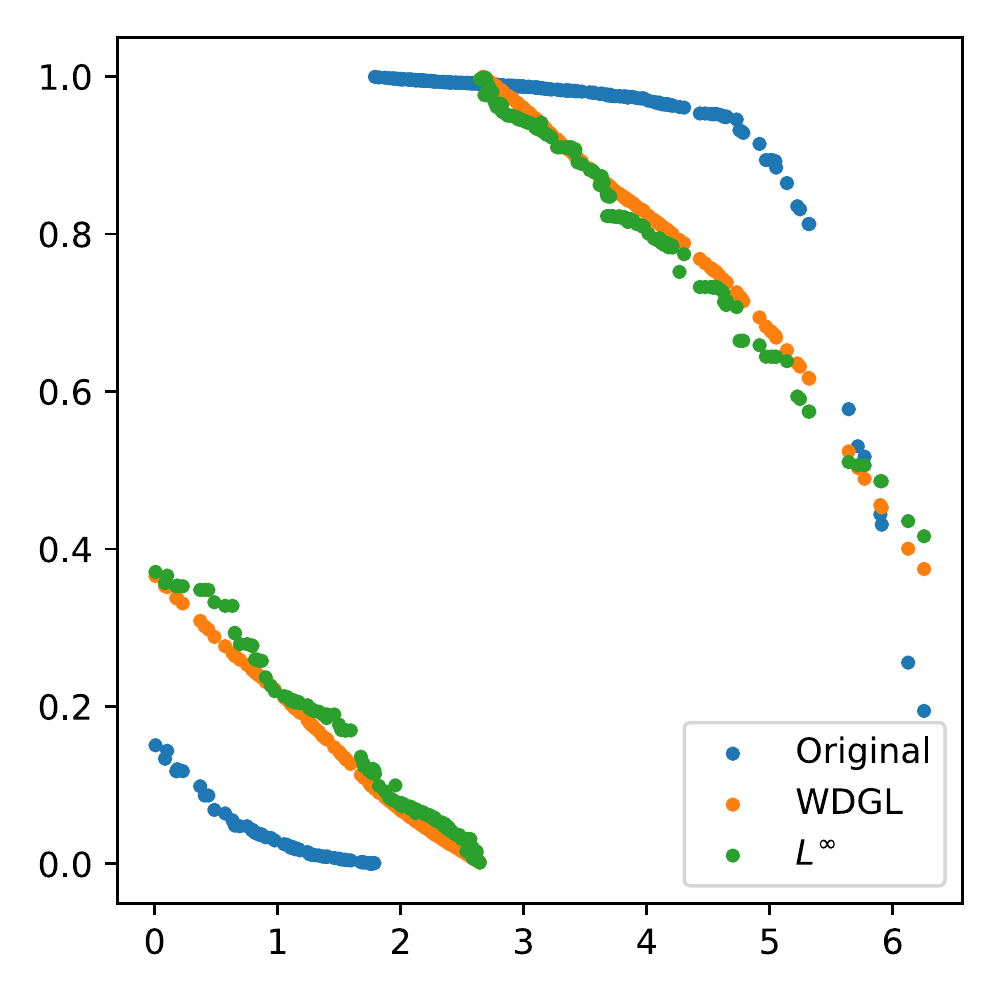}
    \caption{Correlation scatter plot}
\end{subfigure}
\begin{subfigure}{.328\textwidth}
    \centering
    \includegraphics[width=\linewidth]{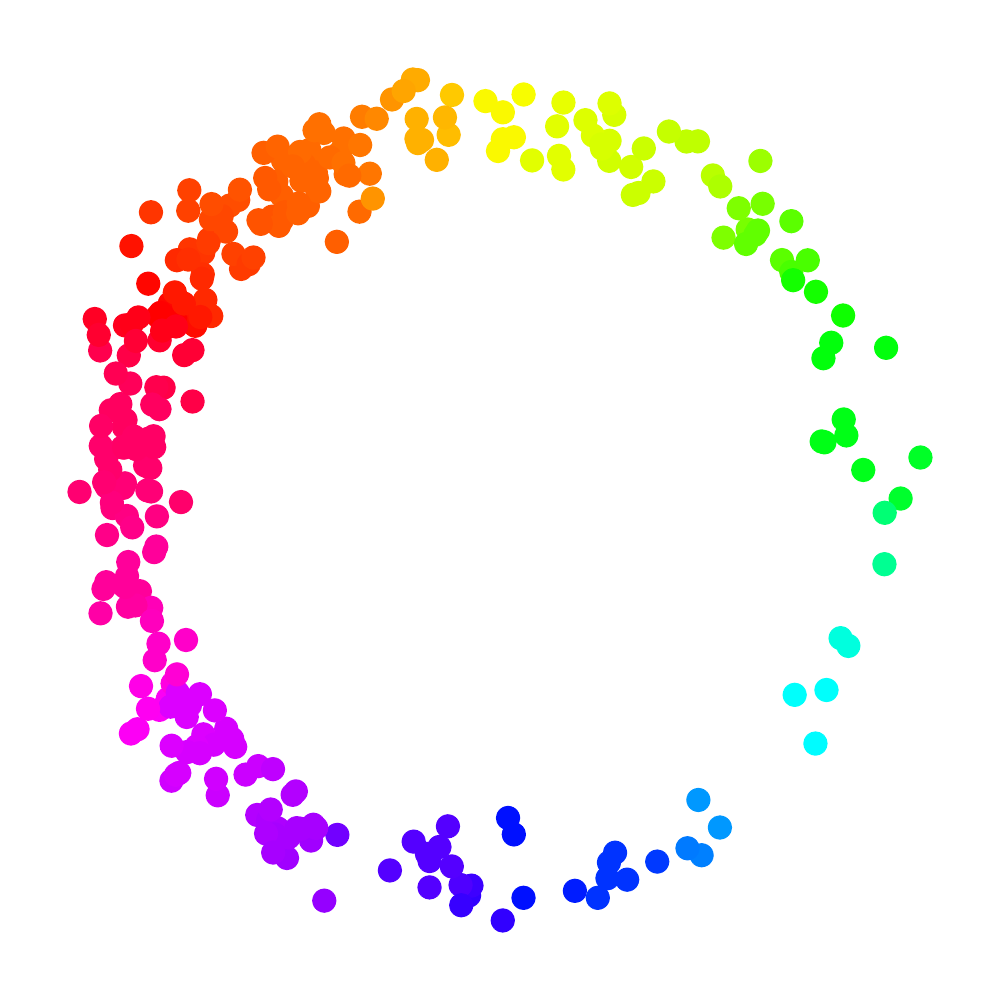}
    \caption{Weighted circular coordinate}
\end{subfigure}
\caption{Results for noisy circle dataset; original circular coordinate (left), correlation scatter plots when we use original method, WDGL, and $L^\infty$-norm optimization (middle), WDGL-circular coordinate (right).}
\label{fig: circle_result_exp}
\end{figure}
We begin with a noisy circle with uneven density.
With the circle, we test our density-robust circular coordinate algorithms.
The circle is parametrized by 
$$
\left\{\begin{array}{l}
x=\sin(t)\\
y=\cos(t),\\
\end{array}\right.
$$
and $t$ is sampled from $\mathcal{N}\left(\pi, \left(0.4\pi\right)^2\right)$, and we add Gaussian noise with a mean of $0$ and a standard deviation of $0.07$.
In the experiment, we sample $300$ points.

We present the original circular coordinate on the dataset, correlation scatter plots when using various methods, and the WDGL-circular coordinate in \Cref{fig: circle_result_exp}.
In this experiment, the results of $L^\infty$-norm optimization and WDGL method are similar as shown in the correlation scatter plot, and both look linear.

We also experiment using different weights to get different weighted circular coordinates.
As we explain in the last part of \Cref{subsec: weighted circular}, we try two weights to get weighted circular coordinates, and those results are similar to the result when we use the WDGL method.

For $L^p$ norm variation, we try $p =2, 4, 6, 10, 20$, and $\infty$.
Note that $L^2$ norm optimizing produces the same result as the original circular coordinate.
From the experiment on this dataset, we confirm that the higher the $p$ value, the more linear the correlation scatter plot appears.

To provide more experimental results, we collect the correlation scatter plots for weighted circular coordinates and $L^p$ norm variation, and present the plots in \Cref{fig: circle_result_scatter_app}.

Now, we study the convergence speed for a $L^\infty$-circular coordinate.
A simple approach to get a $L^\infty$-circular coordinate is just using Algorithm \ref{alg: p_circular}.
In the algorithm, the $L^\infty$ loss is $\|\alpha+d_0 f \|_\infty = \underset{e}{\max}\{|\alpha(e) + (d_0f)(e)|\}$, therefore you can optimize the function $f$.

Instead, we can optimize the function $f$ using Algorithm \ref{alg: infty_circular_p}.
We experiment with several hyperparameters of $\Delta$ and $\eta$, and we observe that of all the results, the algorithm that increases $p$ from $2$ to $50$ appears to be the fastest.
For the experimental details, see \Cref{fig: lp_loss_comp1}

We can also get a $L^\infty$-circular coordinate using Algorithm \ref{alg: infty_circular_softmax}.
See \Cref{fig: lp_loss_comp2} for details.
From the experimental results, we find that the algorithm with softmax with temperature sometimes converges faster than Algorithm \ref{alg: p_circular}. Still, in most cases, they exhibit similar or slower convergence speed than Algorithm \ref{alg: p_circular}.

From the experimental results, we conclude that we can speed up the convergence if we initialize the function $f$ with the original circular coordinate, which is $L^2$-circular coordinate, and then we use Algorithm \ref{alg: infty_circular_p}.

\subsection{Noisy trefoil knot}
\begin{figure}[t]
\begin{subfigure}{.328\textwidth}
    \centering
    \includegraphics[width=\linewidth]{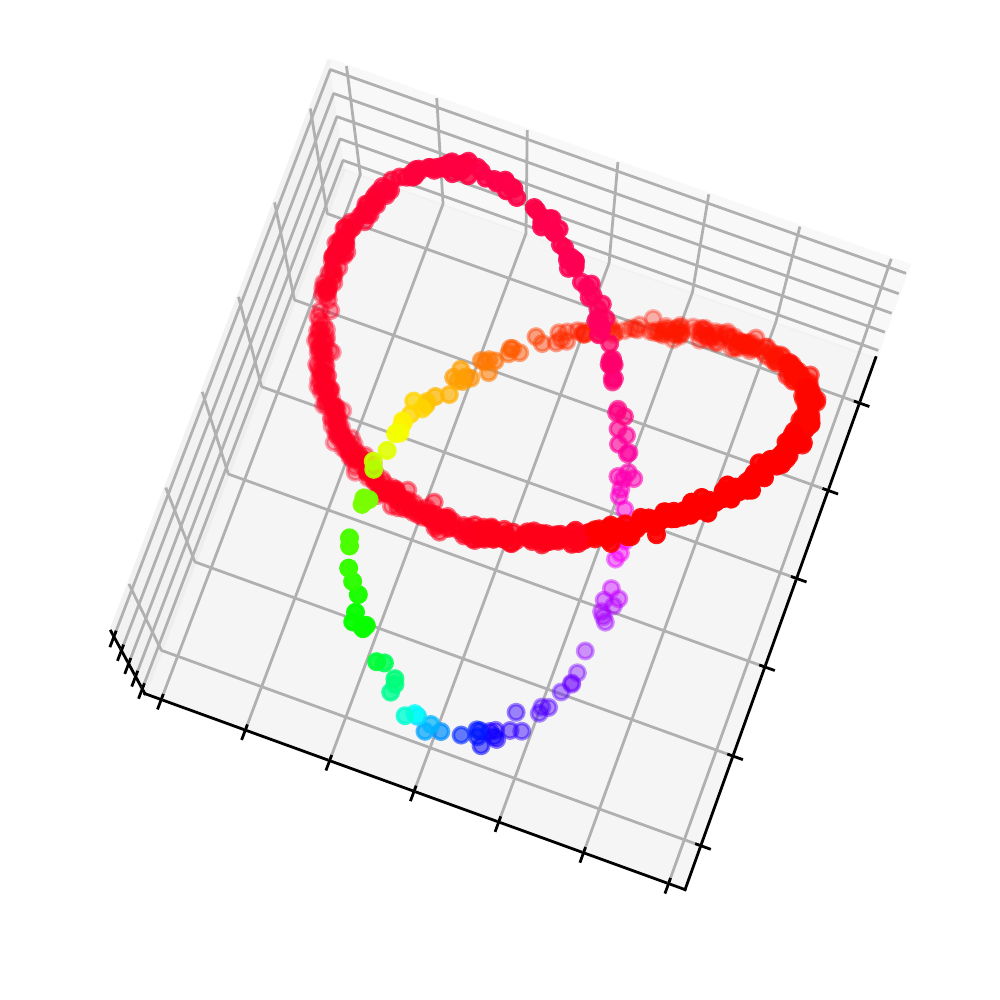}
    \caption{Circular coordinate}
\end{subfigure}
\begin{subfigure}{.328\textwidth}
    \centering
    \includegraphics[width=\linewidth]{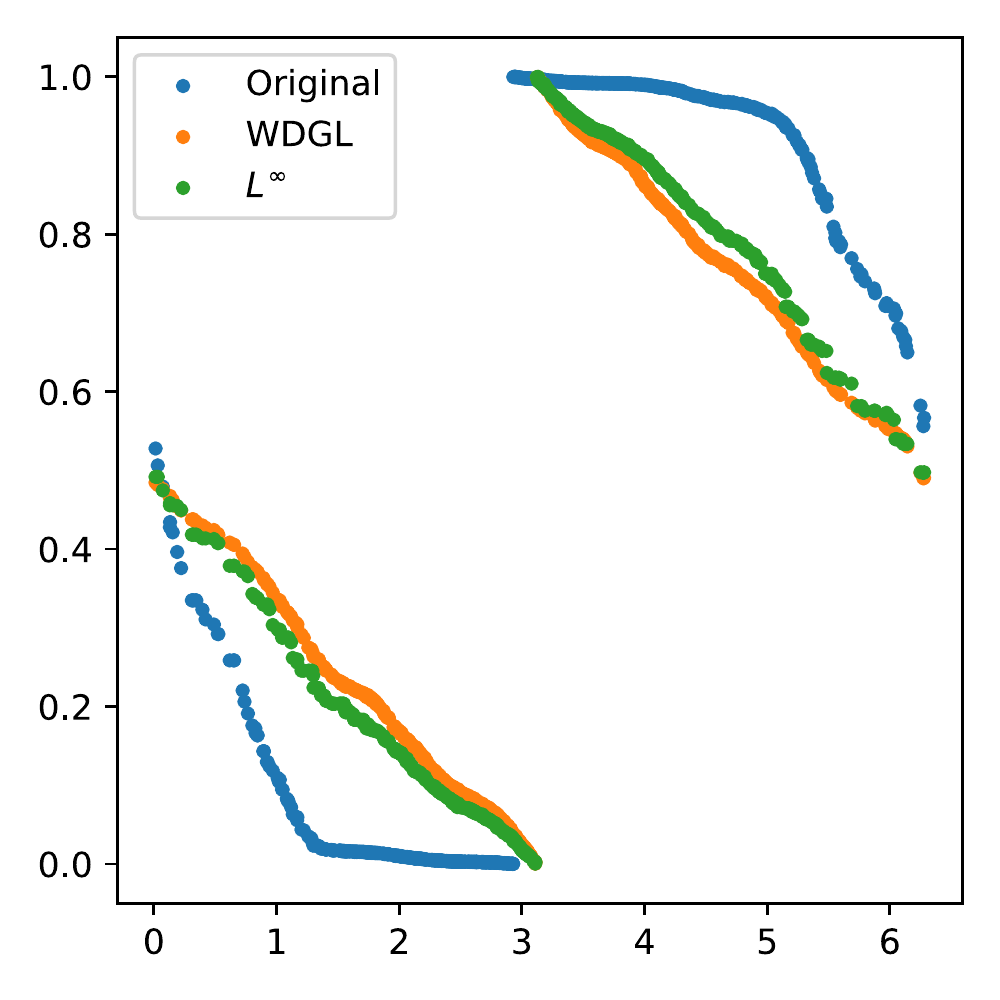}
    \caption{Correlation scatter plot}
\end{subfigure}
\begin{subfigure}{.328\textwidth}
    \centering
    \includegraphics[width=\linewidth]{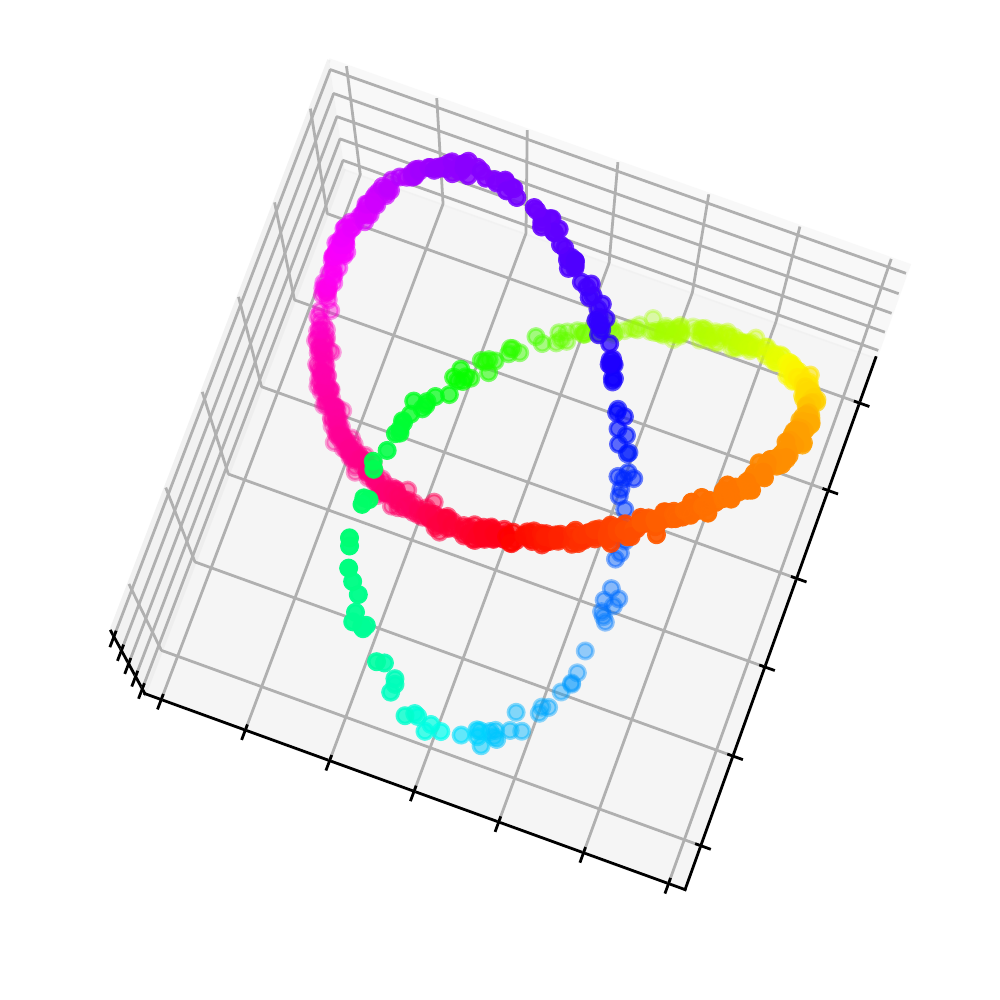}
    \caption{Weighted circular coordinate}
\end{subfigure}
\caption{Results for noisy trefoil knot dataset; original circular coordinate (left), correlation scatter plots when we use original method, WDGL, and $L^\infty$-norm optimization (middle), WDGL-circular coordinate (right).}
\label{fig: knot_result_exp}
\end{figure}

Next, we test our algorithms on a noisy trefoil knot.
We use parametrization 
$$
\left\{\begin{array}{l}
x=\cos (t)+2 \cos (2 t) \\
y=\sin (t)-2 \sin (2 t) \\
z=2 \sin (3 t),
\end{array}\right.
$$
and $t$ is sampled from $\mathcal{N}\left(\pi,\left(0.4\pi\right)^2\right)$, and we add Gaussian noise with a mean of $0$ and a standard deviation of $0.04$.
In the experiment, we sample $900$ points.

We show the original circular coordinate on the dataset, correlation scatter plots with different methods, and the WDGL-circular coordinate in \Cref{fig: knot_result_exp}.
As with the previous experimental results on the noisy circle dataset, the results of $L^\infty$-norm optimization and WDGL-circular coordinate are similar, and those look linear.

We carry out two further experiments on the weighted circular coordinate and experiment on the $L^p$-norm variation technique using the identical $p$, both of which are similar to the experimental setup of the prior noisy circle dataset.
We present the results in \Cref{fig: knot_result_scatter_app}, and the conclusion is also similar to the noisy circle dataset; the three weighted circular coordinates look similar to each other, and a higher $p$ value leads to a more linear correlation scatter plot.

\subsection{Two conjoined circles}
Next, we test our algorithms on two conjoined circles dataset.
In the dataset, there are clear two circles and we can observe those in a persistence diagram induced from the dataset.
Therefore, we conduct the algorithms two times for each cocycle.
To generate the dataset, we make $2$ circles as we do in the noisy circle example, and after taking random rotations for each circle, we attach the two circles to each other.
For each circle, we have coordinate information and it is used to make correlation scatter plots.

\begin{figure}[t]
\begin{subfigure}{.31\textwidth}
    \centering
    \includegraphics[width=\linewidth]{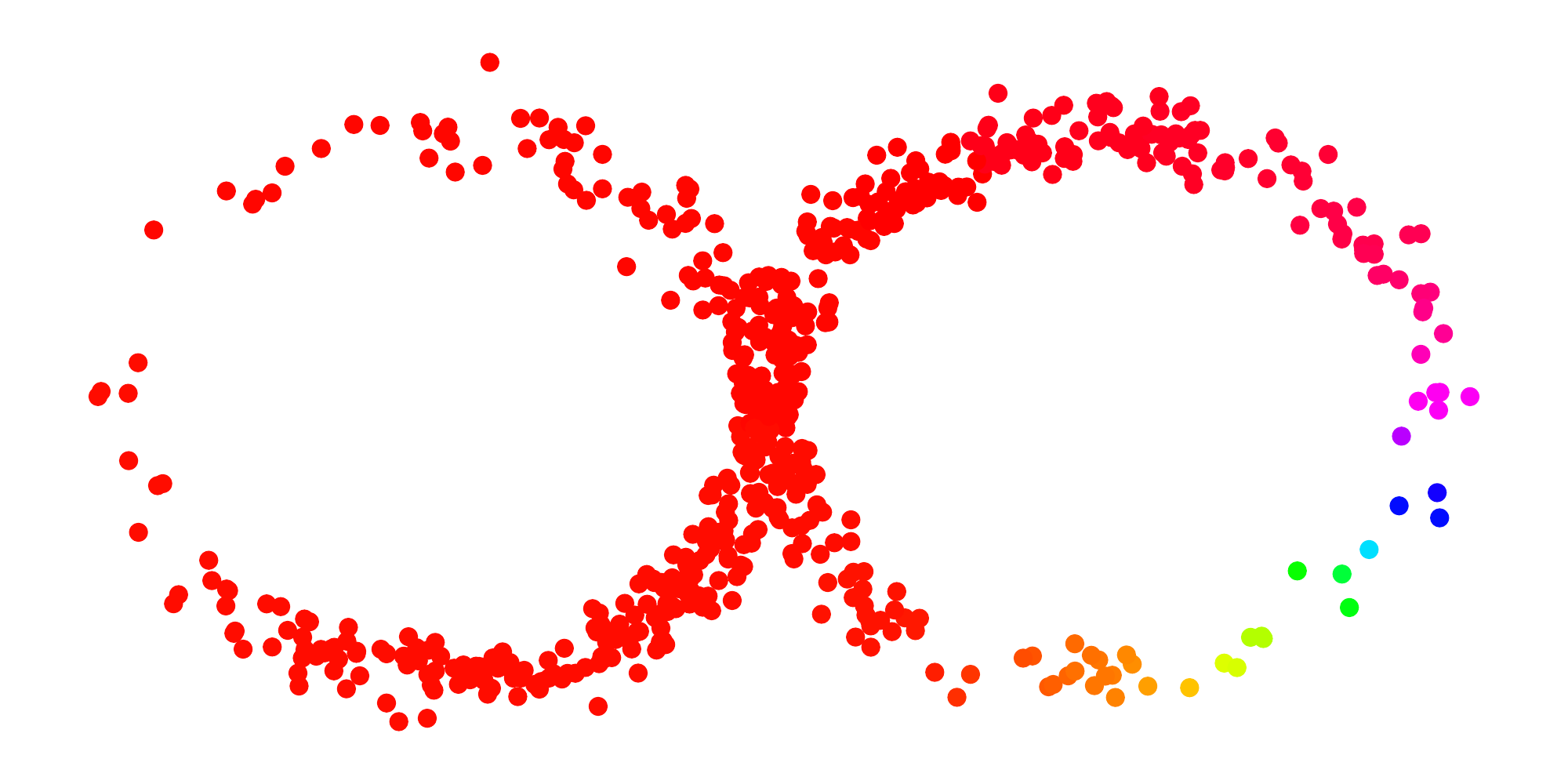}
    
\end{subfigure}
\begin{subfigure}{.21\textwidth}
    \centering
    \includegraphics[width=\linewidth]{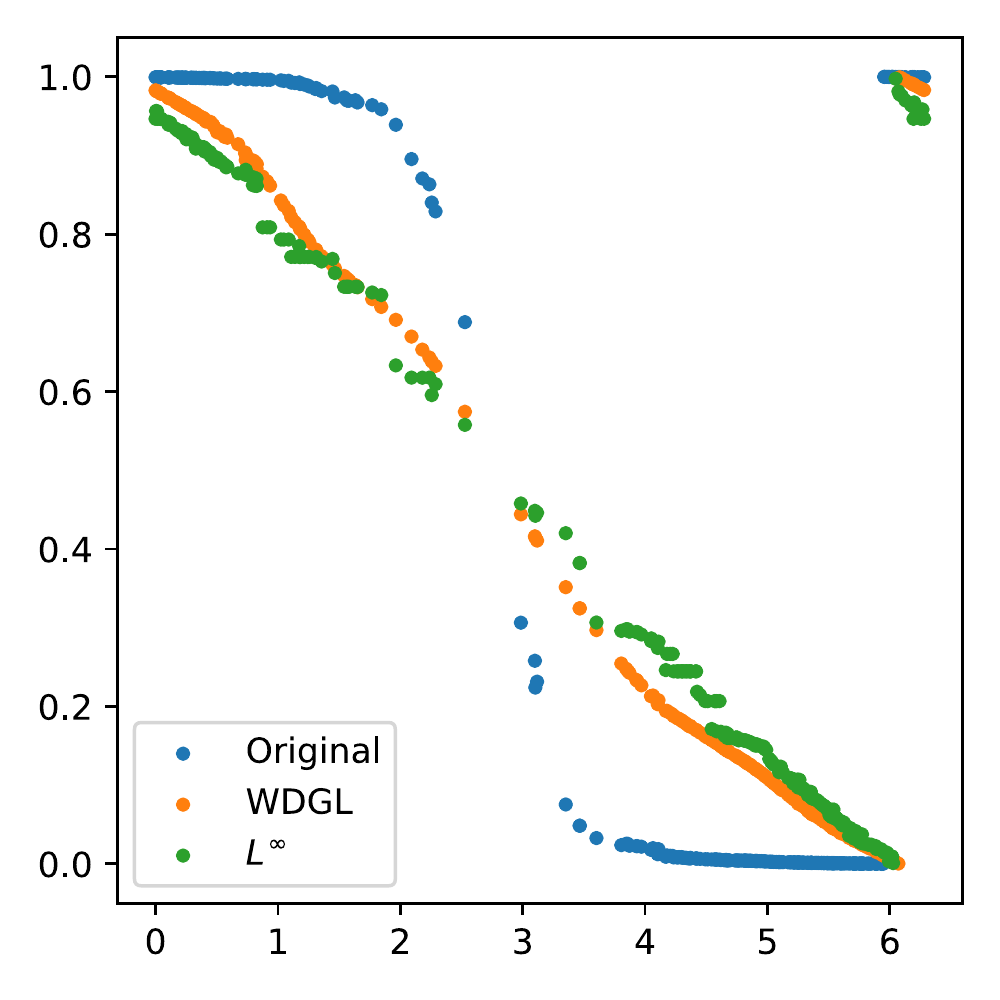}
\end{subfigure}
\begin{subfigure}{.31\textwidth}
    \centering
    \includegraphics[width=\linewidth]{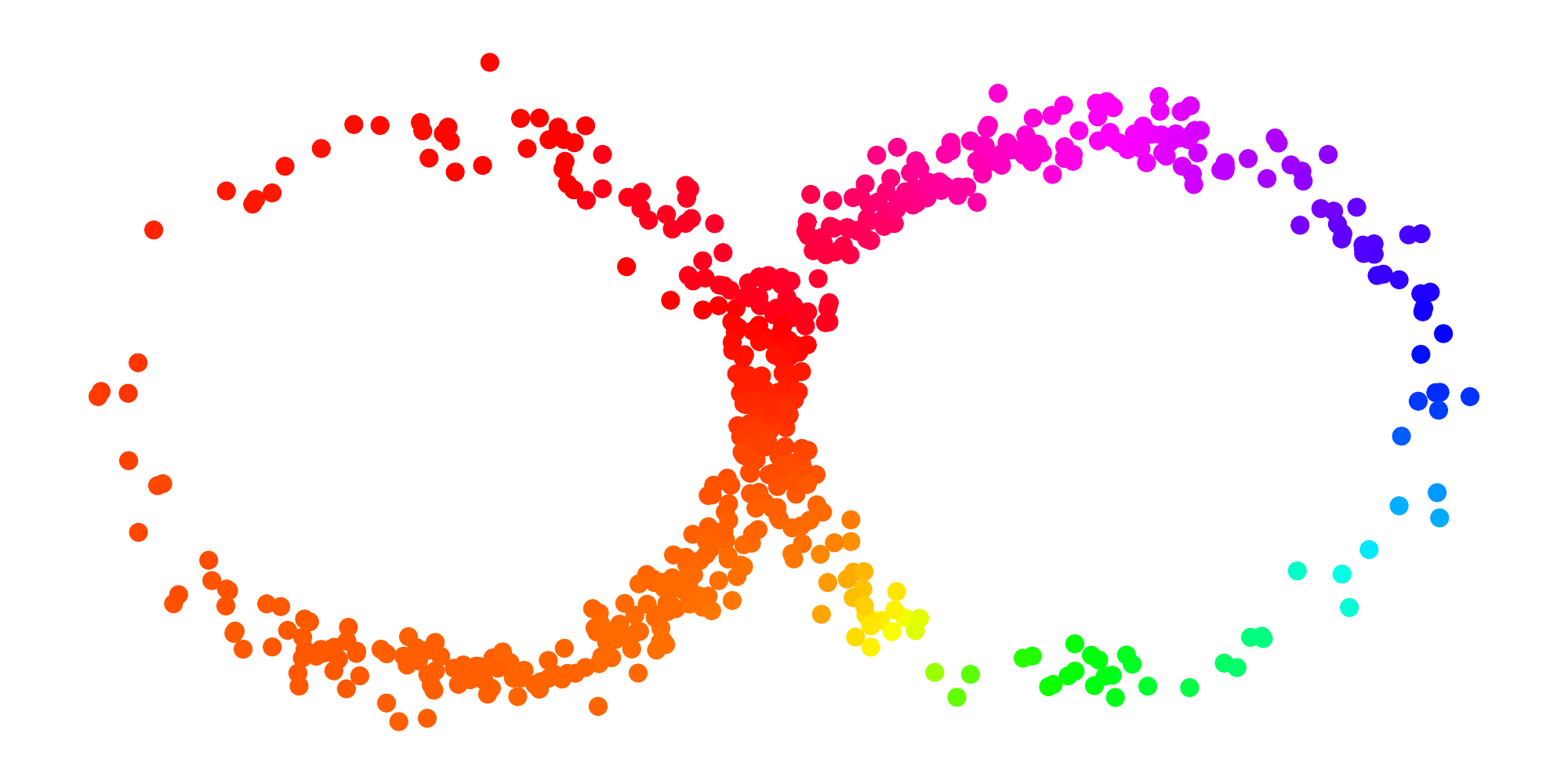}
\end{subfigure}
\begin{subfigure}{.31\textwidth}
    \centering
    \includegraphics[width=\linewidth]{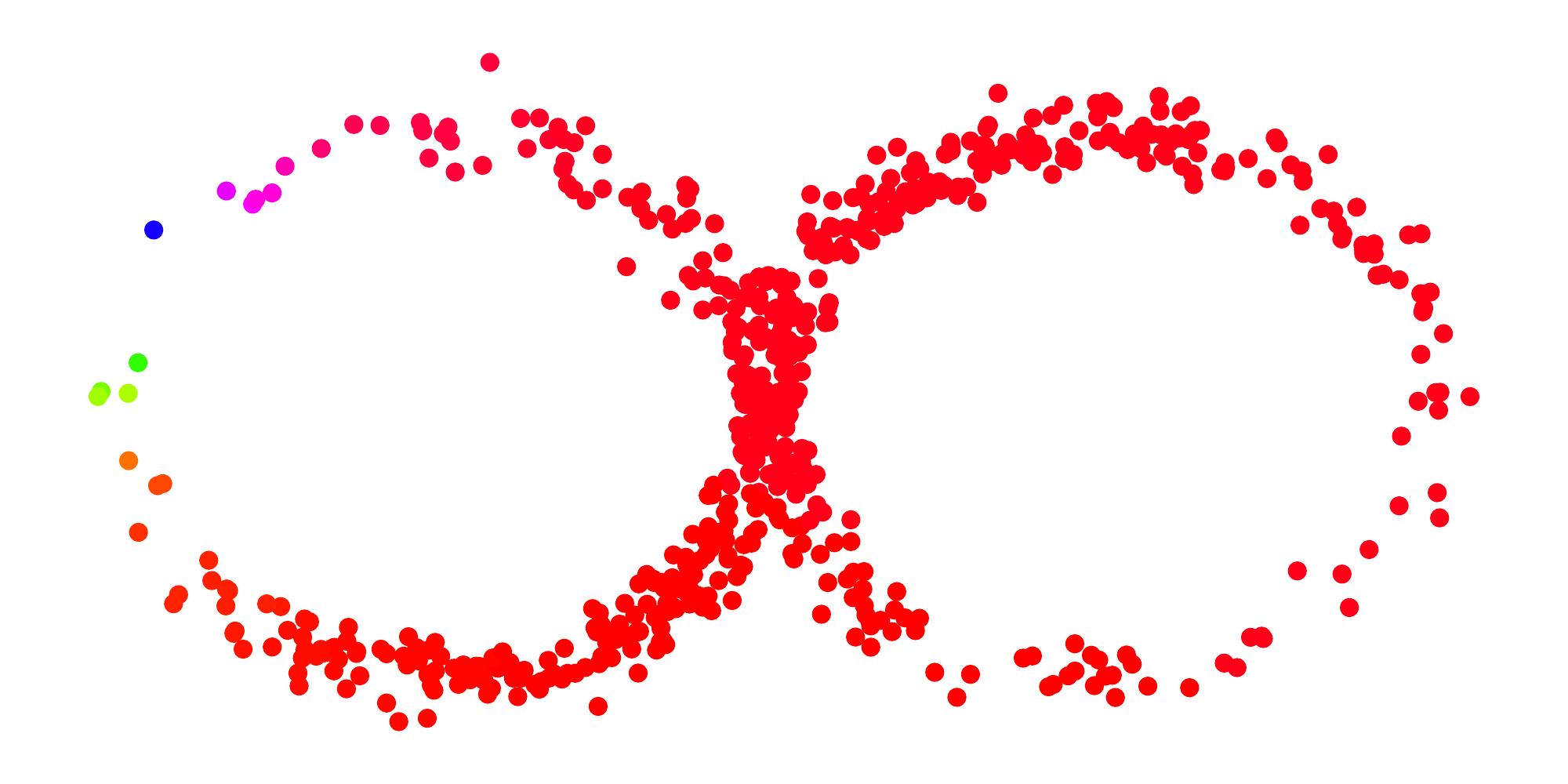}
\end{subfigure}
\begin{subfigure}{.21\textwidth}
    \centering
    \includegraphics[width=\linewidth]{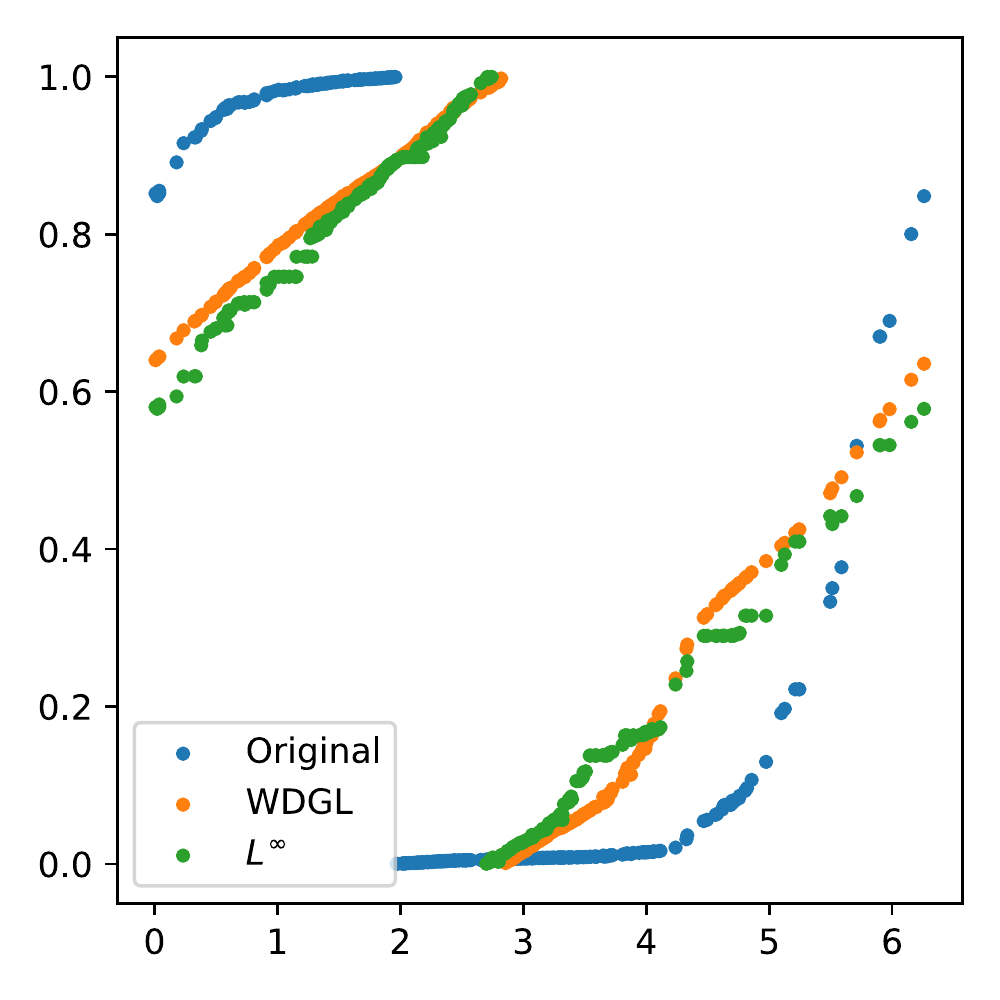}
\end{subfigure}
\begin{subfigure}{.31\textwidth}
    \centering
    \includegraphics[width=\linewidth]{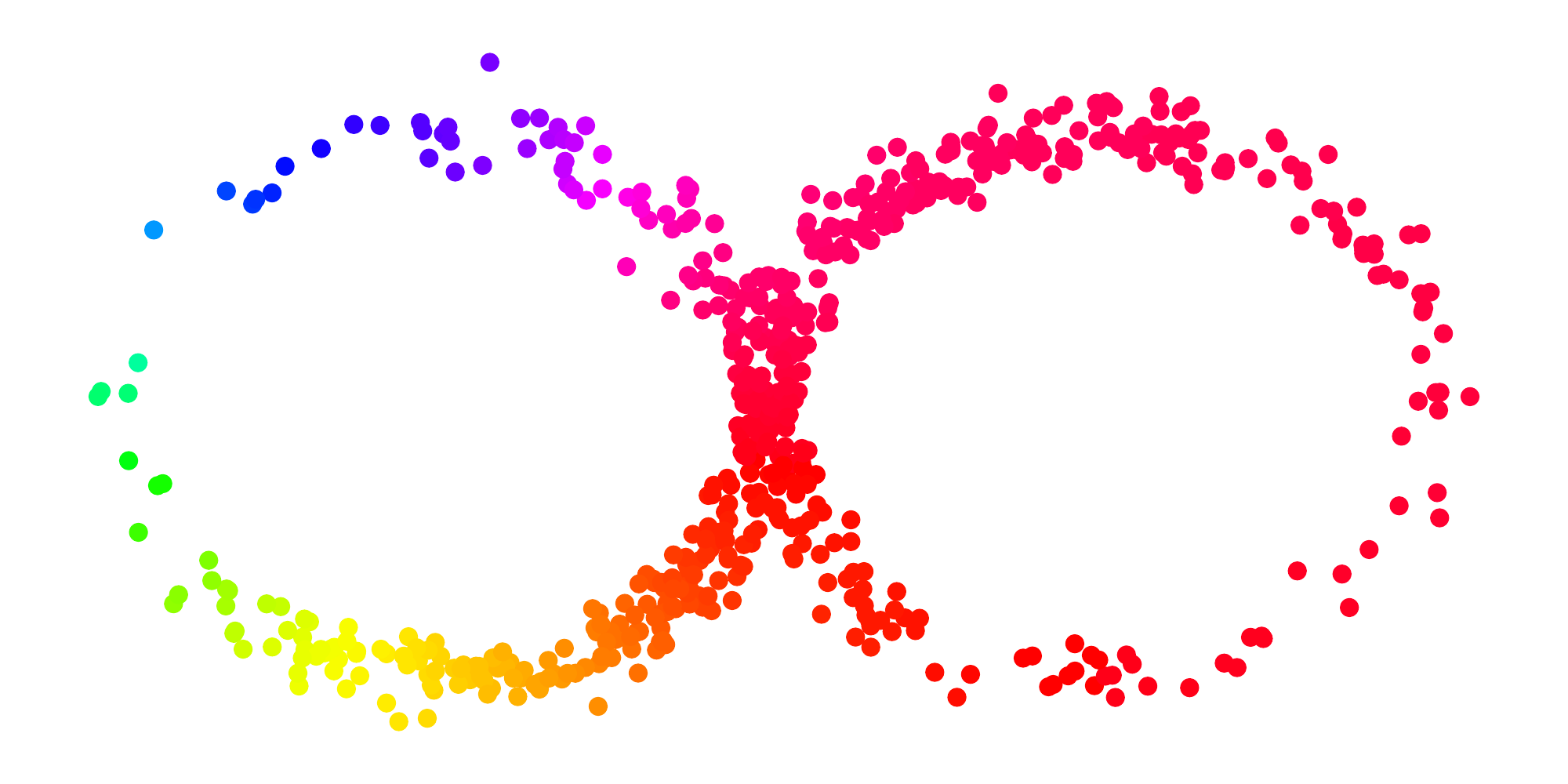}
\end{subfigure}
\caption{Results for two conjoined circles dataset; original circular coordinates (left column), corresponding circular coordinates (middle column), WDGL-circular coordinates (right column).}
\label{fig: conjoined_result_scatter_exp}
\end{figure}

We present the original circular coordinates, corresponding circular coordinates, and WDGL-circular coordinates for each cocycle in \Cref{fig: conjoined_result_scatter_exp}.

On this data, the $L^\infty$-circular coordinate does not look smooth, but it tends to be similar to the WDGL-circular coordinate overall. 

We conduct two supplementary experiments on the weighted circular coordinate and experiment with the $L^p$-norm variation approach, and we present the result in \Cref{fig: conjoined_result_scatter_app}.

\subsection{Torus}
\begin{figure}[t]
\begin{subfigure}{.328\textwidth}
    \centering
    \includegraphics[width=\linewidth]{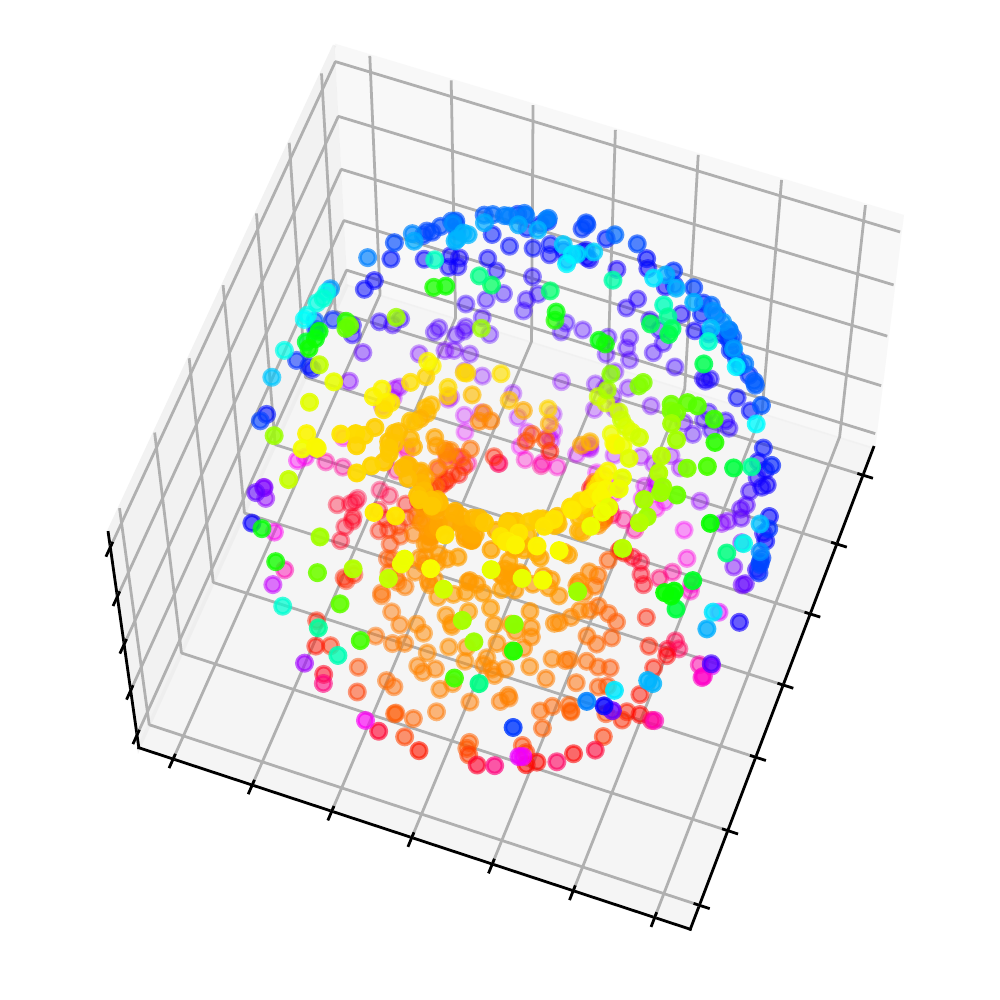}
\end{subfigure}
\begin{subfigure}{.328\textwidth}
    \centering
    \includegraphics[width=\linewidth]{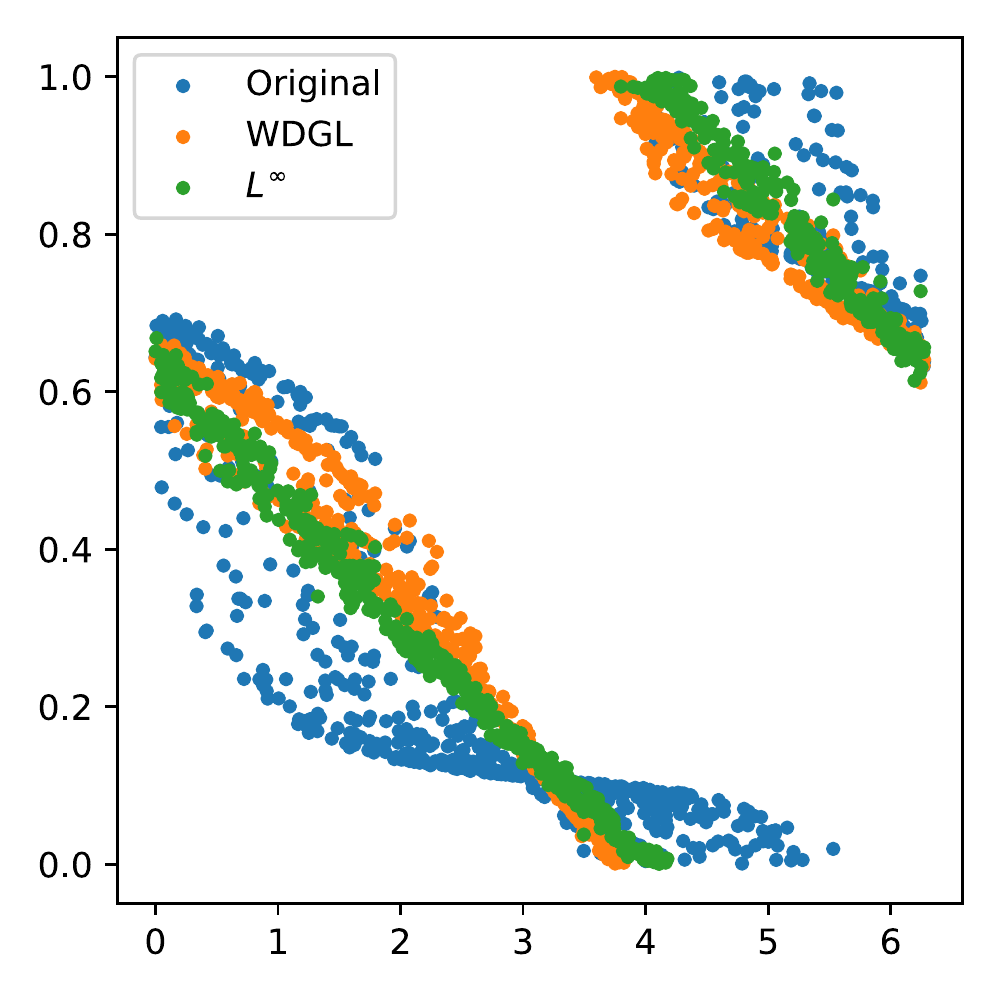}
\end{subfigure}
\begin{subfigure}{.328\textwidth}
    \centering
    \includegraphics[width=\linewidth]{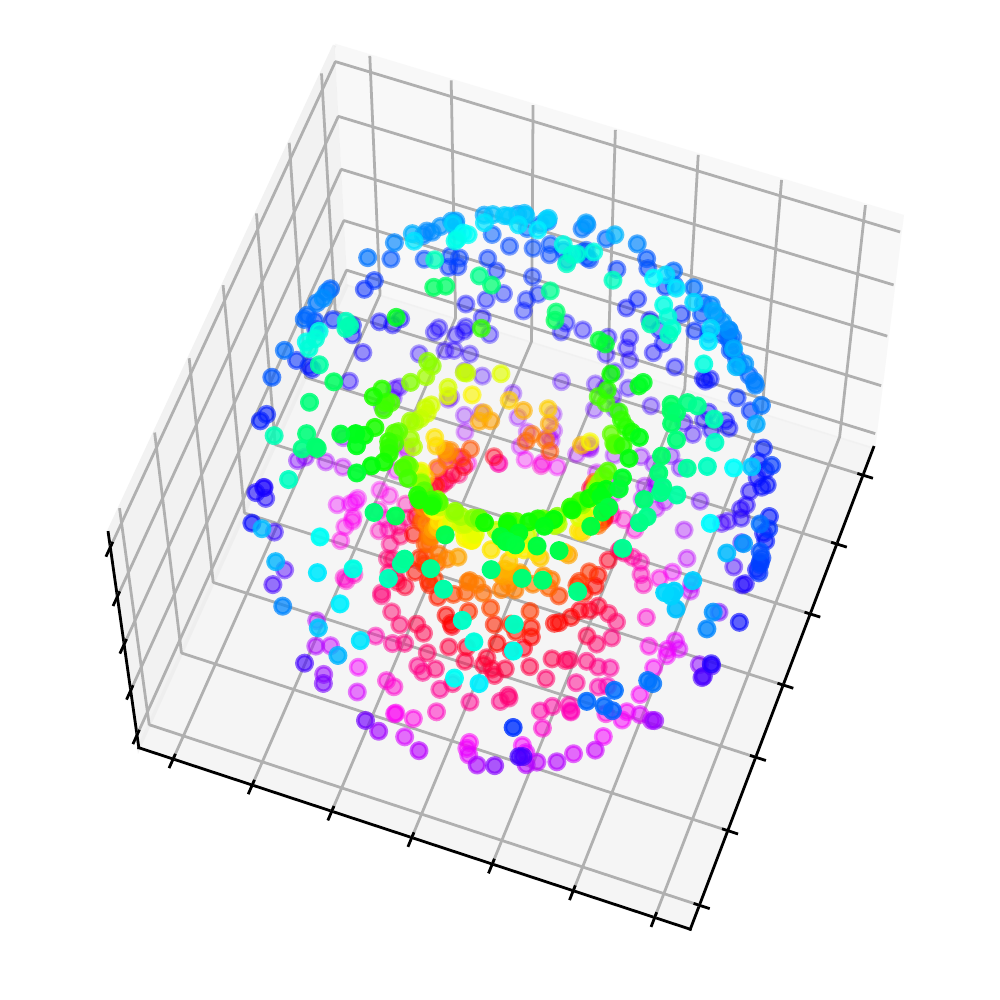}
\end{subfigure}
\begin{subfigure}{.328\textwidth}
    \centering
    \includegraphics[width=\linewidth]{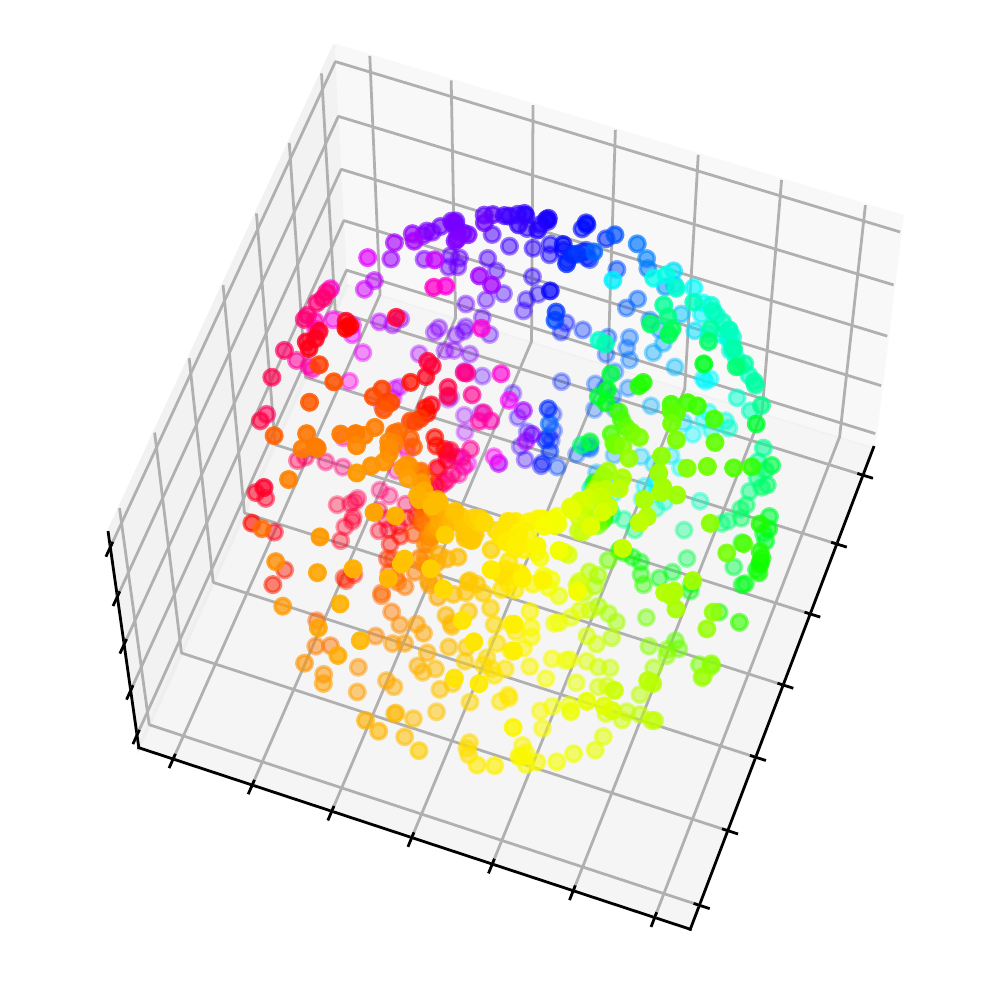}
\end{subfigure}
\begin{subfigure}{.328\textwidth}
    \centering
    \includegraphics[width=\linewidth]{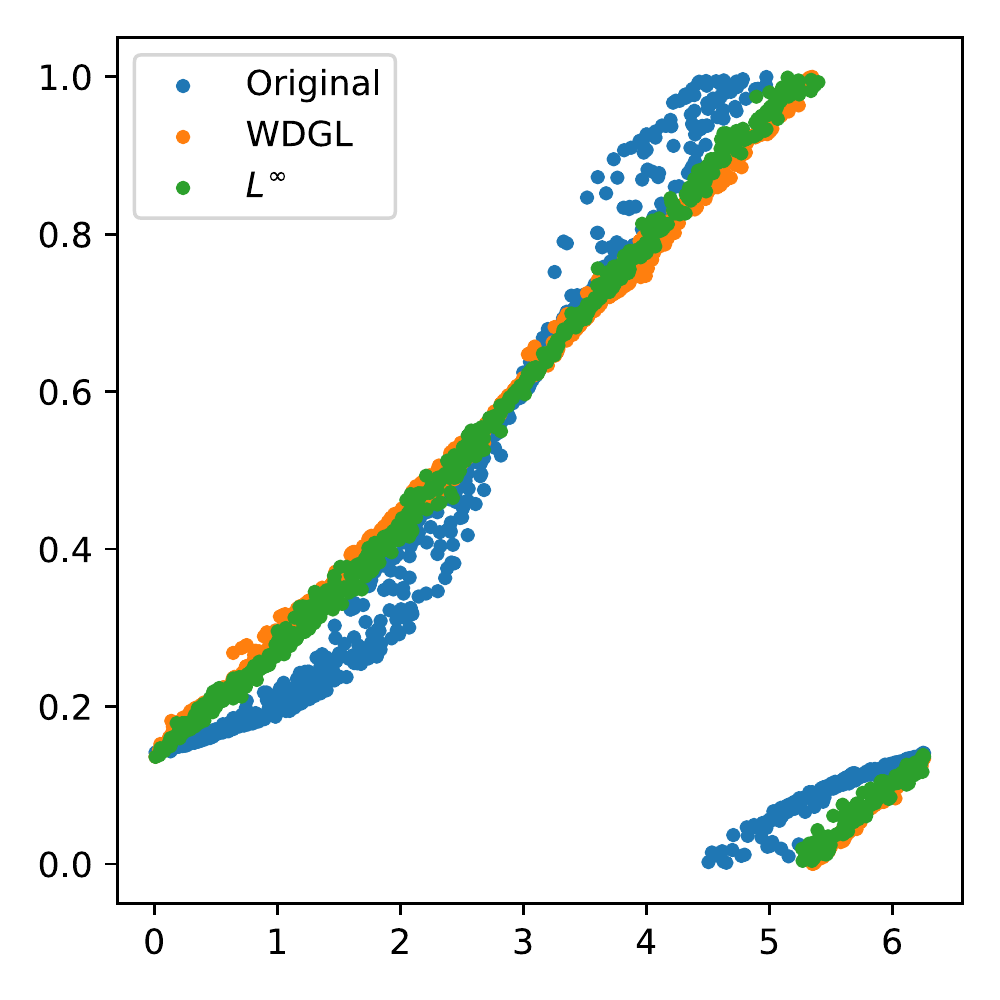}
\end{subfigure}
\begin{subfigure}{.328\textwidth}
    \centering
    \includegraphics[width=\linewidth]{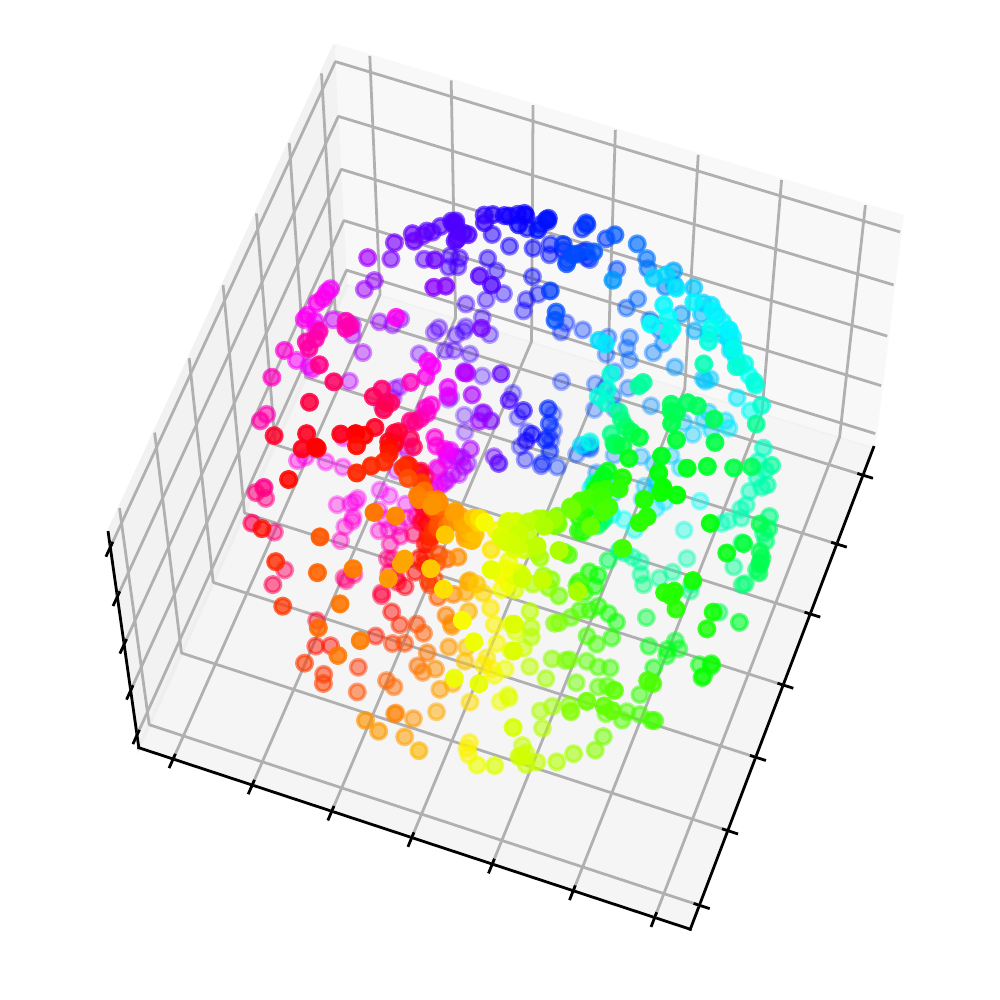}
\end{subfigure}
\caption{Results for torus dataset; original circular coordinates (left column), corresponding circular coordinates (middle column), WDGL-circular coordinates (right column).}
\label{fig: torus_result_scatter_exp}
\end{figure}
We use torus data as the final example of synthetic data to illustrate our methodology.
This dataset has two cocycles, as in the previous dataset; the Betti number of the torus is two, and this can be seen in the persistence diagram.
The torus is parametrized by
$$
\left\{\begin{array}{l}
x=(4 + 2\cos(s))\cos(t)\\
y=(4 + 2\cos(s))\sin(t)\\
z=2 \sin (s),
\end{array}\right.
$$
and we use Gaussian mixture model; we sample $(s, t)$ from a probability density function $p(x) = \frac{1}{2}\left(\mathcal{N}(x\mid (\pi,0), \Sigma) + \mathcal{N}(x\mid (0, \pi), \Sigma)\right)$ where $\Sigma$ is $\operatorname{diag}(0.4\pi, 0.4\pi)^2$.
In this experiment, we sample $800$ points.

We show the original circular coordinates, corresponding circular coordinates, and the WDGL-circular coordinate for each cocycle in \Cref{fig: torus_result_scatter_exp}.
There are $2$ cocycles corresponding to the meridian direction and the longitude direction.
Note that when using the WDGL method, the result of the meridian directional cocycle is not very linear in the correlation scatter plot.
The reason is that we are not dealing with a flat torus; the uniform probability distribution on the torus is not induced by the uniform probability distribution on $[0, 2\pi]\times [0, 2\pi]$.
However, using $L^\infty$-norm optimization, we get the linear result in the correlation scatter plot using the meridian directional cocycle.

As in the previous dataset, we perform two further experiments on the weighted circular coordinate and $L^p$-circular coordinates.
The results are shown in \Cref{fig: torus_result_scatter_app}.

\subsection{COIL-100}\label{subsec: exp-coil100}
In this experiment, we use the COIL-100 \cite{nene1996columbia} dataset as a real dataset to test our method.
Each object has 72 images, each taken every $5$ degrees as the object is rotated $360$ degrees.
Therefore, it is natural to consider that the discretized $S^1$ for each item is embedded in the space of the images.
The dataset was created by the Columbia Object Image Library (COIL) at Columbia University and is commonly used for object recognition and classification tasks in machine learning and computer vision research.

We obtain circular coordinates using Euclidean distance between the images.
To show the results we use t-SNE and PCA, which are popular techniques for dimensionality reduction to obtain a scatter plot of the images to visualize the high-dimensional data.

Using the circular coordinates obtained earlier, we express them in color on the t-SNE or PCA result for several cocycles.
This allows us to identify topological structures hidden in the result of the dimensionality reduction.

Furthermore, we select several examples of cocycles that show more hidden topological structures well in the $L^\infty$-circular coordinate than in the original circular coordinate method.

\begin{figure}[t]
\begin{subfigure}{.33\textwidth}
    \centering
    \includegraphics[width=\linewidth]{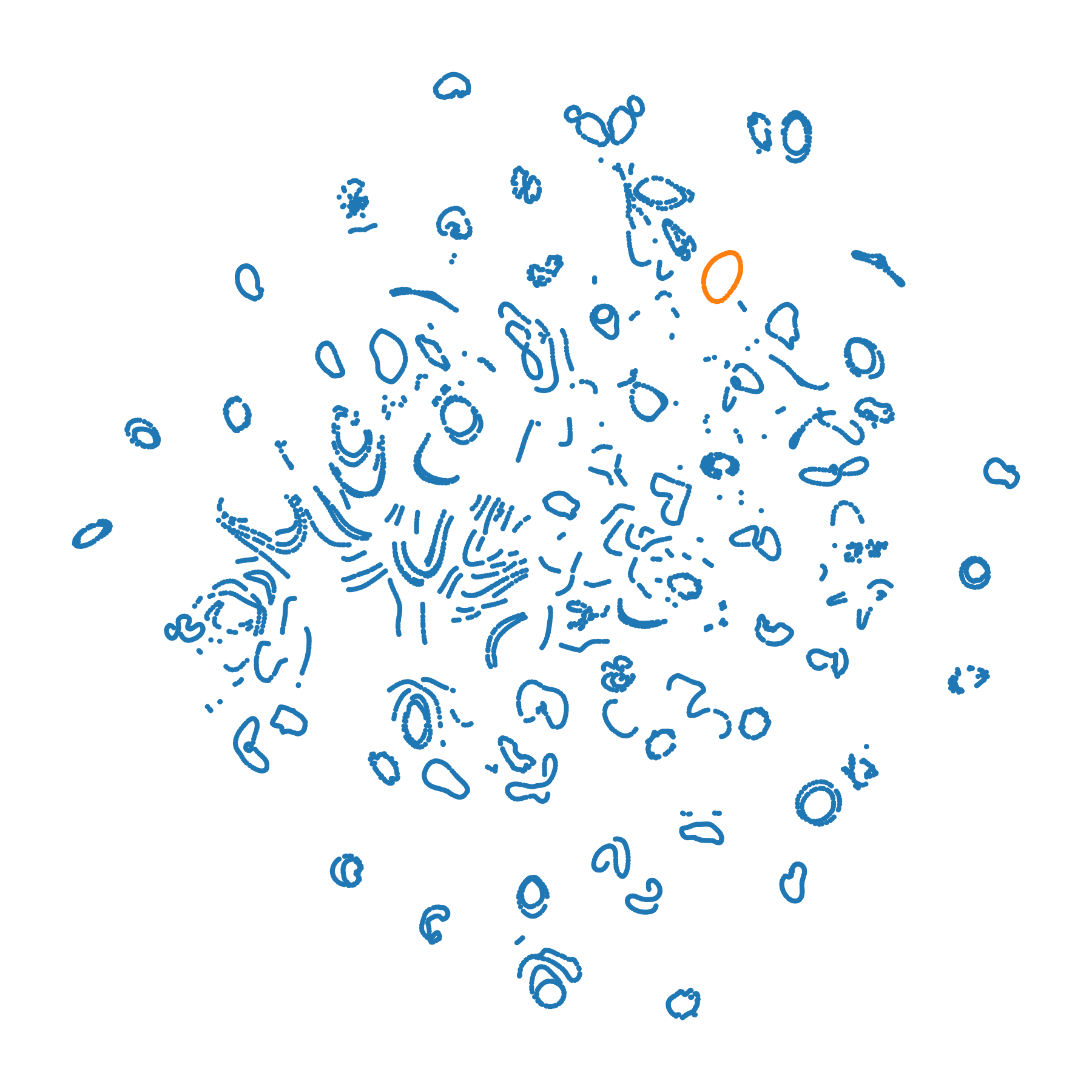}
    \caption{$S^1$ shape embedding.}
\end{subfigure}
\hspace{2cm}
\begin{subfigure}{.33\textwidth}
    \centering
    \includegraphics[width=\linewidth]{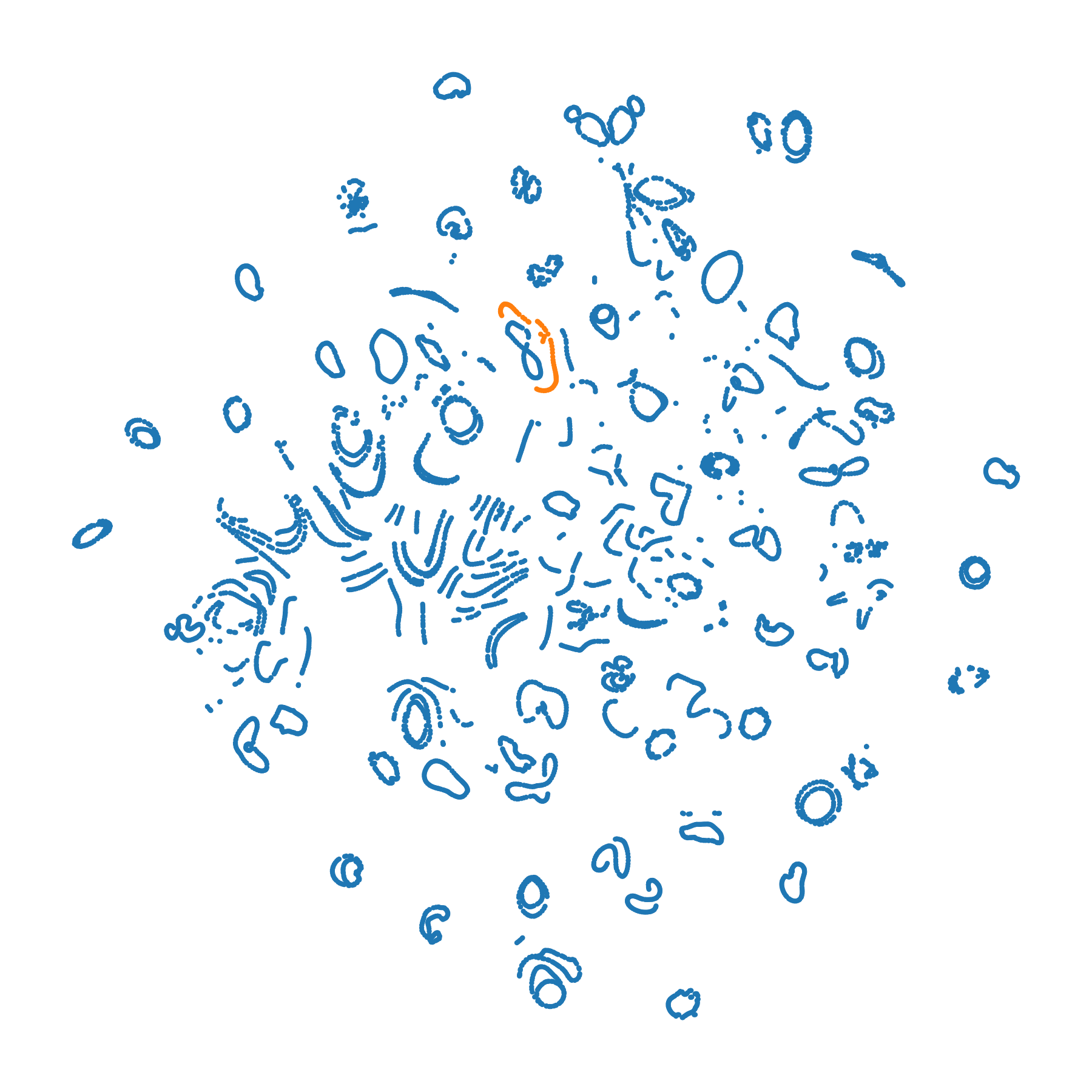}
    \caption{Broken $S^1$ shape embedding.}
\end{subfigure}
\caption{}
\label{fig: pca and tsne result}
\end{figure}

\begin{figure}[t]
    \begin{center}
        \includegraphics[width=0.35\linewidth]{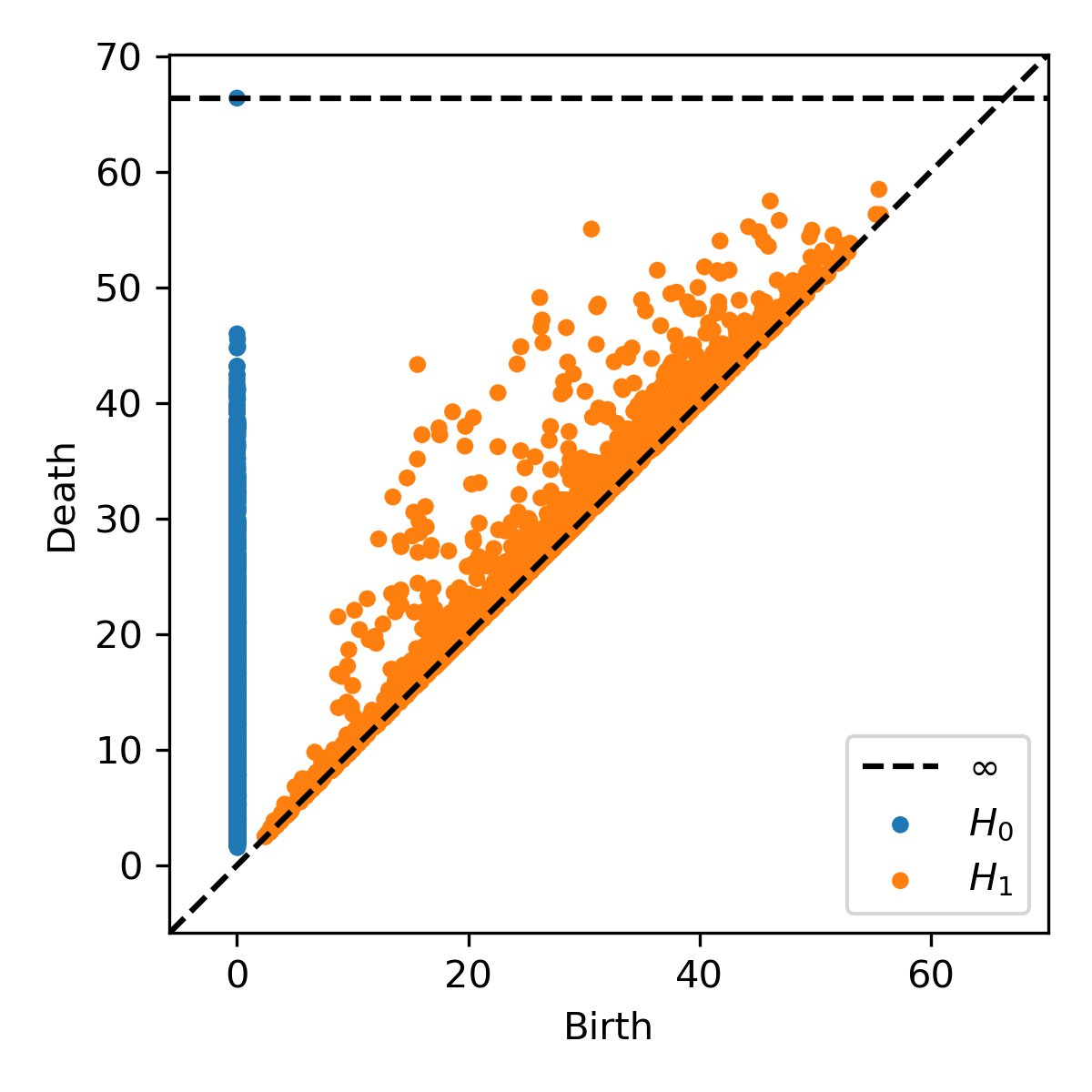}
    \end{center}
    \caption{Persistence diagram on COIL-100 data using Euclidean distance.}
    \label{fig: persistence diag}
\end{figure}
First, we present two circular coordinates on t-SNE in \Cref{fig: pca and tsne result}.
Many cycles are apparent in the t-SNE result as a ``$S^1$ shape''.
However, some cycles of some objects appear cut off even in t-SNE.
There are several potential causes for this; in the initialization of the embedding, far points are difficult to merge, and trying to embed many objects at once can cause conflicts with each other.
To identify the circular structures in the low-dimensional embeddings, we evaluate the circular coordinates for each cocycle in the persistence diagram.
We show the persistence diagram in \Cref{fig: persistence diag}.

We obtain circular coordinates by using the entries in the persistence diagram.
Then, we choose some of these circular coordinates to help us to see the topological structure in the embeddings.

\begin{figure}[b]
\begin{subfigure}{.24\textwidth}
    \centering
    \includegraphics[width=\linewidth]{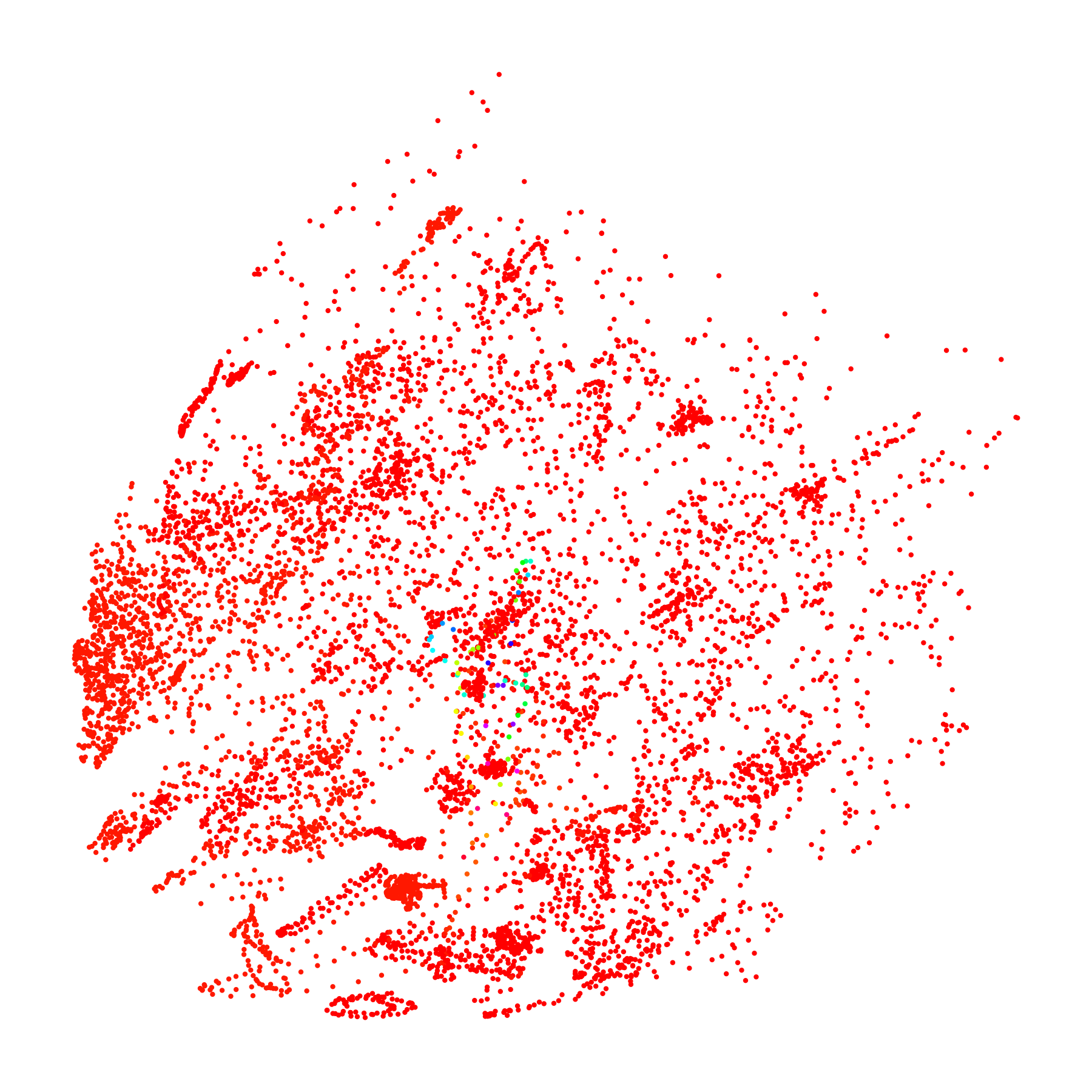}
    \caption{}
\end{subfigure}
\begin{subfigure}{.24\textwidth}
    \centering
    \includegraphics[width=\linewidth]{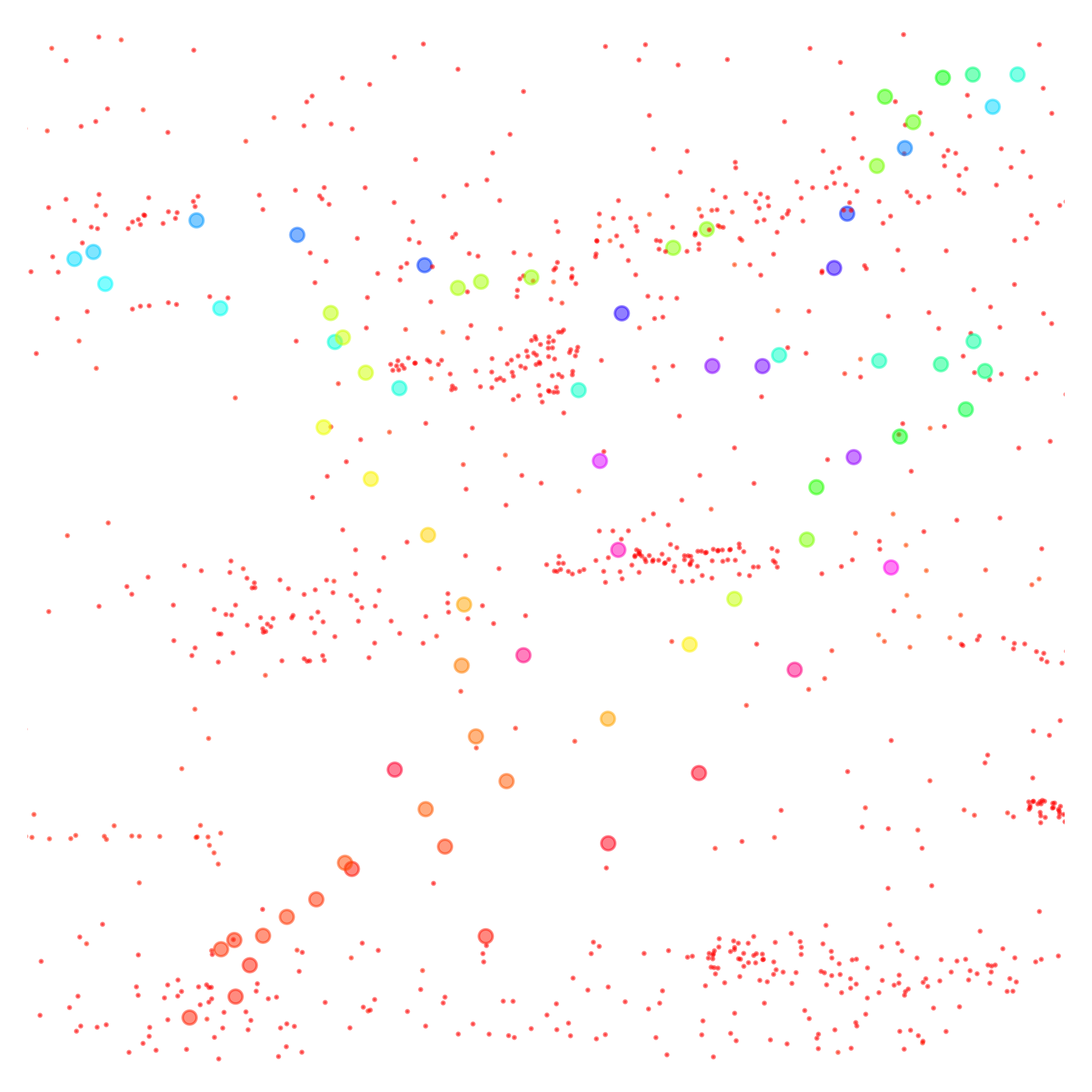}
    \caption{}
\end{subfigure}
\begin{subfigure}{.24\textwidth}
    \centering
    \includegraphics[width=\linewidth]{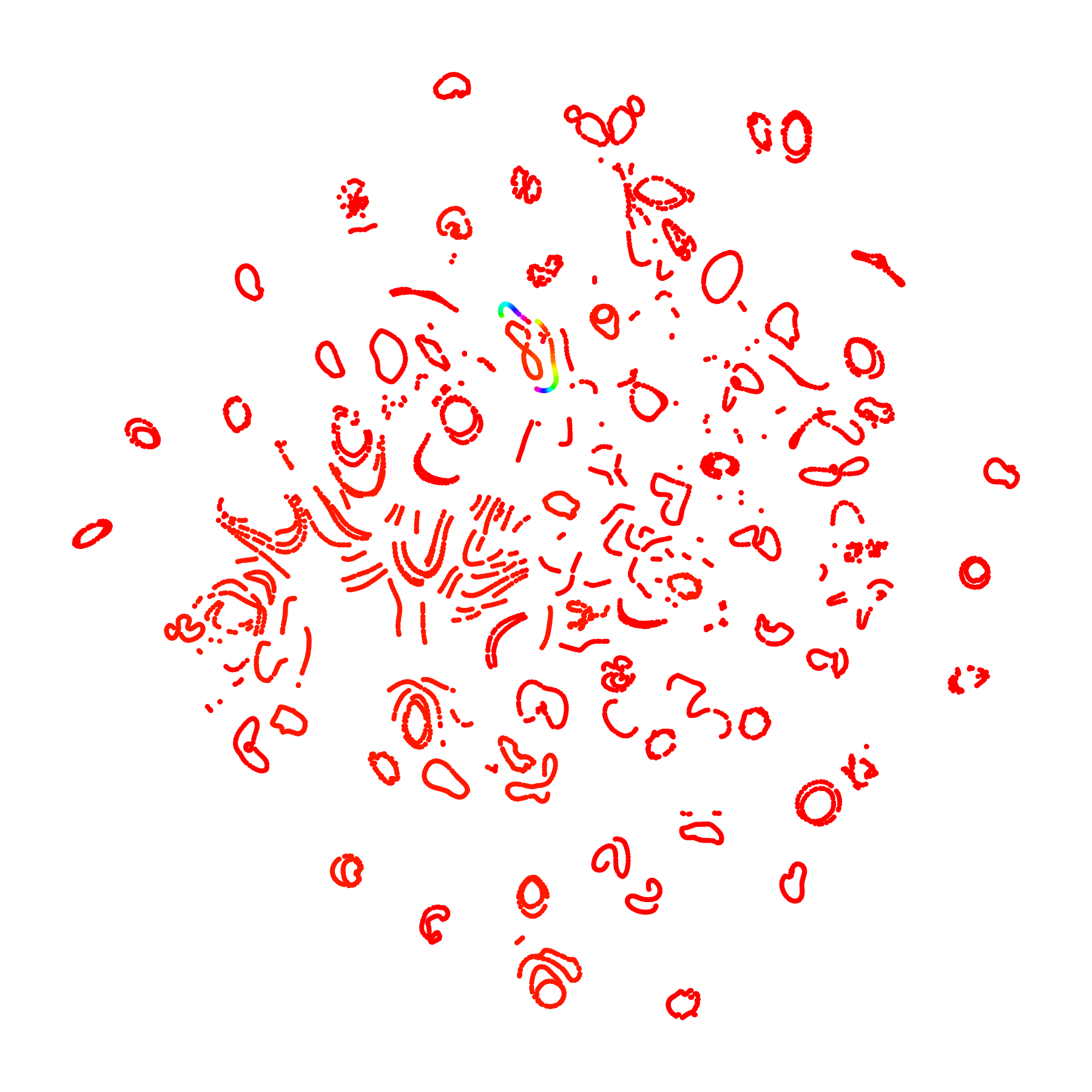}
    \caption{}
\end{subfigure}
\begin{subfigure}{.24\textwidth}
    \centering
    \includegraphics[width=\linewidth]{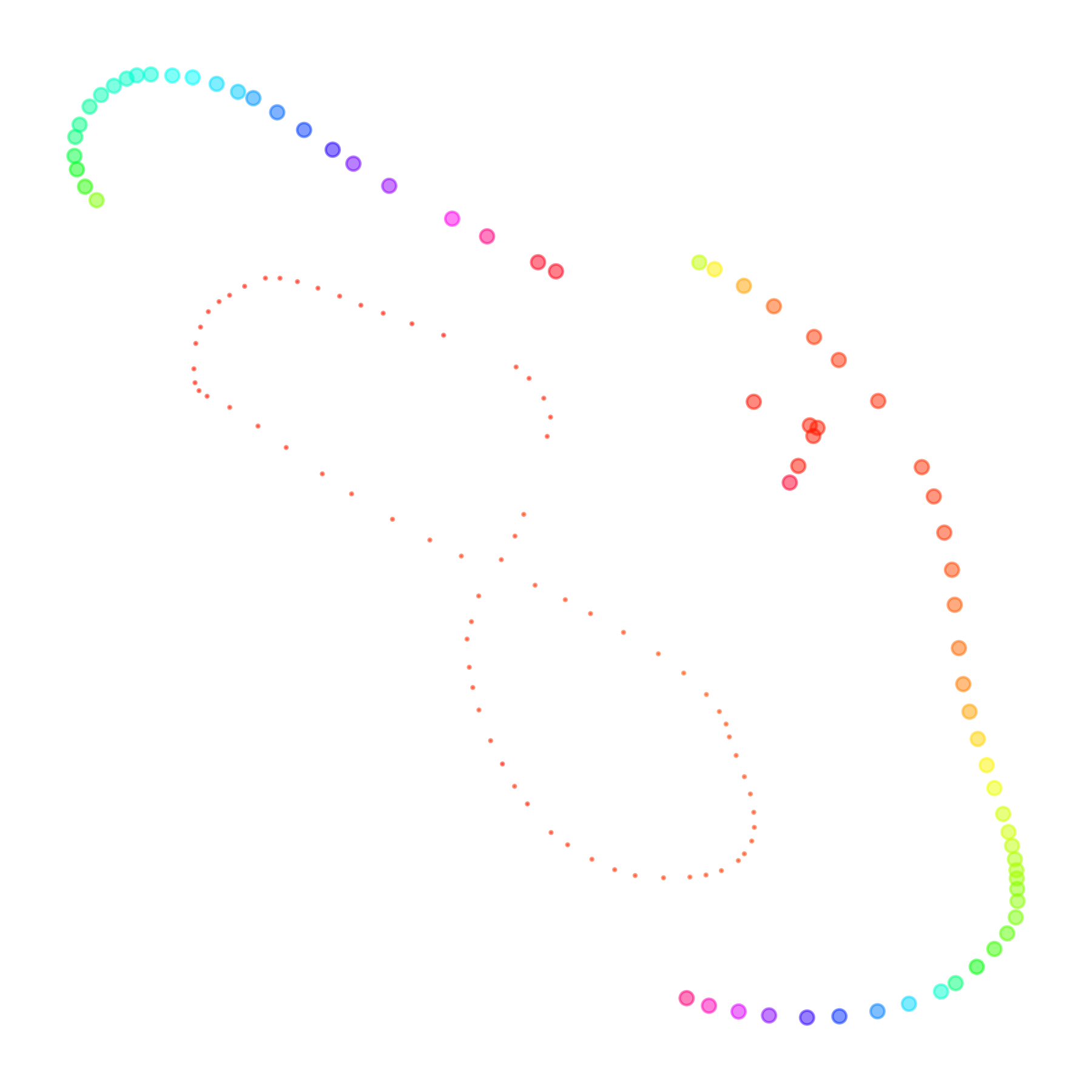}
    \caption{}
\end{subfigure}
\caption{The outcome of the circular coordinate is represented by the two figures on the left with PCA and t-SNE ((A), (C)).
We magnify where the points of the corresponding object are gathered and show them by adjusting the point size respectively ((B), (D)).
}
\label{fig: cc pro}
\end{figure}

After sampling one of the disconnected objects in the t-SNE, we conduct an experiment.
We present the result in \Cref{fig: cc pro}.
The cycle of this object looks entangled in PCA, but the circular coordinate identifies the cycle's location.
In the t-SNE result, Although the cycle looks to be broken in the embedding, the circular coordinate reveals that the two portions are actually from the same cycle.

\begin{figure}[ht]
\begin{subfigure}{.35\textwidth}
    \centering
    \includegraphics[width=\linewidth]{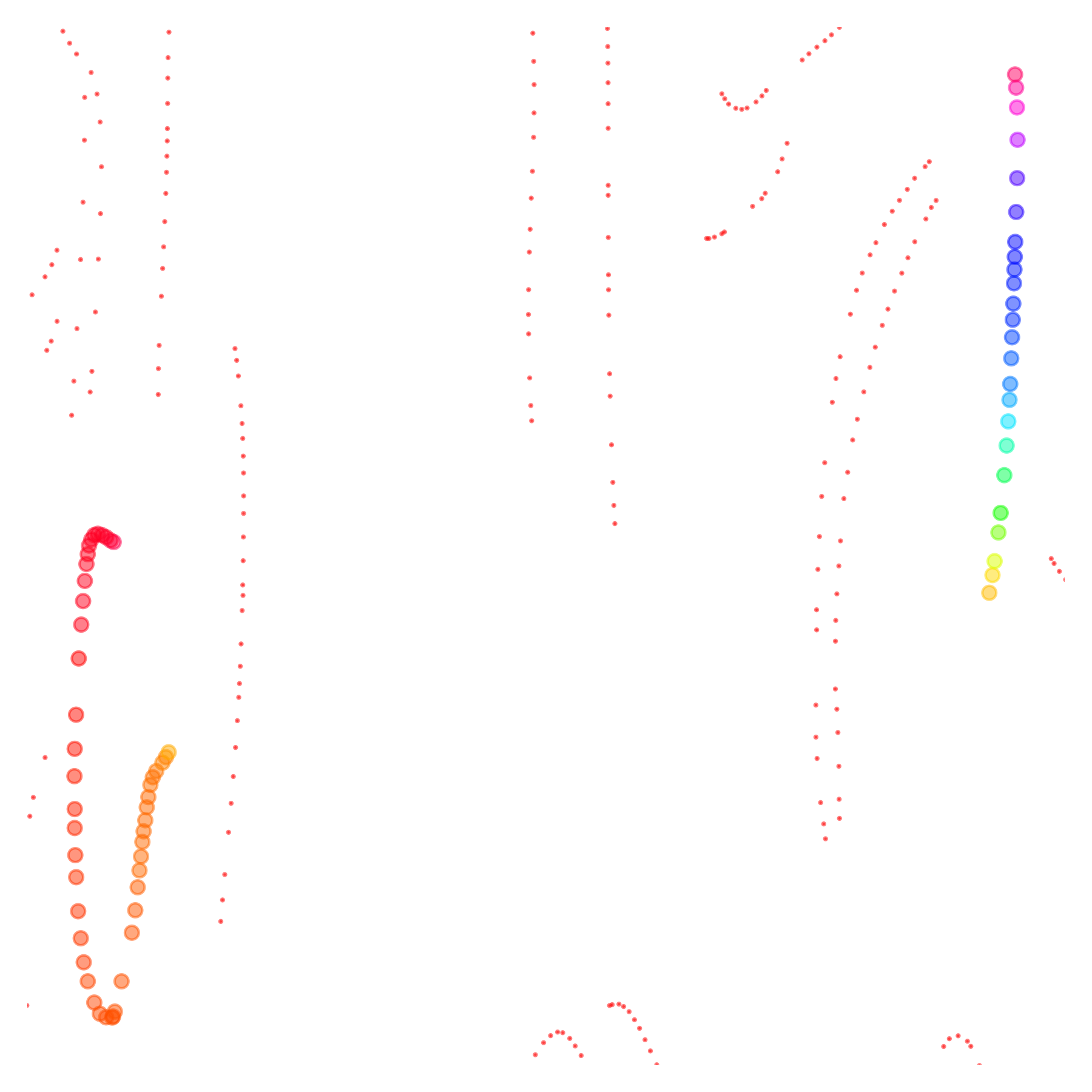}
    \caption{Original circular coordinate.}
\end{subfigure}
\hspace{3 cm}
\begin{subfigure}{.35\textwidth}
    \centering
    \includegraphics[width=\linewidth]{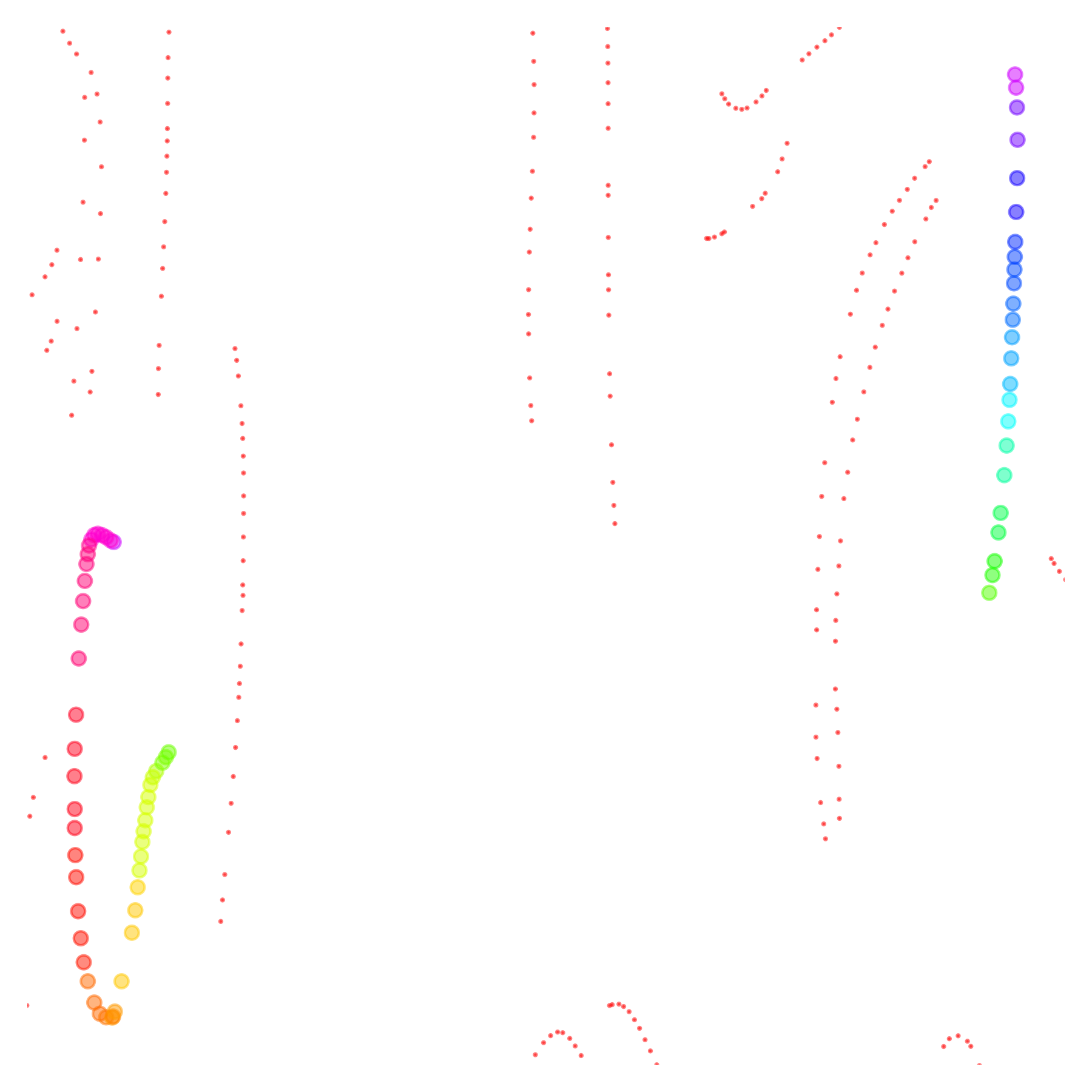}
    \caption{$L^\infty$-circular coordinate.}
\end{subfigure}
\caption{
t-SNE result and circular coordinates. 
We magnify where the points of the corresponding object are gathered and show them by adjusting the point size respectively.
}
\label{fig: Lp pro}
\end{figure}

The following example in \Cref{fig: Lp pro} illustrates a situation in which $L^\infty$-circular coordinate is useful.
In this instance, one cycle is split into two in the t-SNE result.
Unfortunately, in the original circular coordinate, the value of the circular coordinate in the left part is similar to the value of the circular coordinate in the non-cycle part, making it challenging to recognize that the left part is also a part of the cycle.

However, in the case of $L^\infty$-circular coordinate, the circular coordinate on the left and right portions is smoothly and evenly distributed, suggesting that the two parts constitute one cycle.

\section{Discussion}

In our study, we conducted a thorough investigation into the circular coordinate to analyze data.
As part of this investigation, we looked at how these approaches handle data that has uneven density, or where the distribution of data points is not consistent throughout the dataset.

Through our observations, we were able to identify areas where the results produced by these methods could be improved.
In order to address this issue, we developed two new approaches that are more robust to different densities of data.
These new approaches are designed to produce more accurate and reliable results even when the data has uneven density.

In the weighted circular coordinate method, the circular coordinate changes sensitively depending on how we choose weight.
Therefore, it is important to choose the weight for each edge.
In such a sense, we present a circular coordinate using WDGL, and it works as we expected in synthetic datasets.
But the variable $t_n$ in \Cref{thm: LB approxi} remains our hyper-parameter; when $n$ goes to infinity, we know that our method will work well if we use $t_n=n^{-\frac{1}{k+2+\alpha}}$ for any $\alpha>0$, but we need to determine $t_n$ manually for a relatively small limited amount of data.

For $L^p$-circular coordinate, optimization time is a major obstacle.
We know experimentally and theoretically that $L^p$-circular coordinate is robust to change of density functions as $p$ increases, but it takes quite a long time to optimize this.
For this reason, we studied how to speed up optimization and we achieved quite a lot of speed-up, but it still takes a lot of time to get $L^\infty$-circular coordinate.

\bibliographystyle{unsrt}
\bibliography{main}

\clearpage
\appendix

\section{Proof of theorems}
\label{sec: proofOfTheorems}
\begin{proof}[Proof of \Cref{thm: unique and exist of the solution}]
    Let $\alpha_H$ be the corresponding harmonic cocycle.
    Since $\alpha_H$ is cohomologous to $\alpha$, we have $f\in\mathcal{C}^0$ such that $\alpha - \alpha_H = d_0f$.
    Since $\alpha_H$ is in $\ker d_0^*$, we get
    $$
    d_0^*\alpha = d_0^*(\alpha - \alpha_H) = d_0^*d_0f = \Delta_0 f,
    $$
    and $f$ is one of the solutions to the Dirichlet problem.
    Let us assume that the graph is connected and $g$ is another solution to the problem.
    Then, we get $\Delta_0(f-g) = 0$.
    Note that $\ker\Delta_0$ is isomorphic to $H^0(X^\epsilon;\mathbb{R})$
    Since $H^0(X^\epsilon;\mathbb{R})$ is the space of constant functions for each connected component, $f-g$ is a constant.
\end{proof}

\begin{proof}[Proof of \Cref{thm: LB approxi}]
Note that $\frac{1}{n}\mathbf{L}_n^t$ is an approximation of 
$$
\mathbf{L}^t f(p) = \int_M w(p, y)(f(p) - f(y))P(y) \; dV(y)
$$
where $dV$ is a volume form of the manifold $M$ satisfying $\int_M P(y)\;dV(y) = 1$.
From Hoeffding's inequality \cite{hoeffding1994probability}, since we have 
$$
\left| \frac{1}{t}w(x, x_i)(f(x) - f(x_i)) \right| < \frac{1}{t(4\pi t)^{k/2}P(x_i)}\left|f(x) - f(x_i)\right| \leq \frac{a}{t(4\pi t)^{k/2}}
$$
for a constant $a>0$, we have
$$
\mathbb{P}\left[\left|\frac{1}{nt}\mathbf{L}_n^t f(p)-\frac{1}{t}\mathbf{L}^tf(p) \right|>\epsilon_0\right] \leq 2 \exp\left(-\frac{\epsilon_0^2 n (4\pi)^k}{2a^2}t^{k+2}\right).
$$
Therefore, by our choice of $t_n$ in \Cref{thm: LB approxi},
\begin{equation}\label{eq: probability_limit}
\lim_{n \rightarrow \infty} \mathbb{P}\left[\left|\frac{1}{nt_n}\mathbf{L}_n^{t_n} f(p)-\frac{1}{t_n}\mathbf{L}^{t_n} f(p)\right|>\epsilon_0\right]=0.
\end{equation}

Let us define a function $w_t^\prime (x, y) = \frac{1}{(4\pi t)^{k/2}P(y)}e^{-\frac{\|x-y\|^2}{4t}}$.
If we take an epsilon ball $B_p(\epsilon)\subset \mathbb{R}^n$ and denote $B_p(\epsilon) \cap M$ by $B$, then from \cite{belkin2008towards}, we have 
$$
(4\pi t)^{k/2}\left|\mathbf{L}^t f(p)-\int_B w^\prime(p, y)(f(p)-f(y))P(y)\;dV(y)\right|=o\left(t^b\right)
$$
for any $b>0$ since $w^\prime(p, y) = w(p, y)$ for every $y\in B$.
Since we know
$$
\frac{1}{t}\int_B w^\prime(p, y)(f(p) - f(y)) P(y)\; dV(y) =
\frac{1}{t(4\pi t)^{k/2}}\int_B e^{-\frac{\|p-y\|^2}{4t}} (f(p) - f(y))\;dV(y)
$$
converges to $\Delta_M f(p)$ from \cite{belkin2008towards}, with \eqref{eq: probability_limit}, the claim holds.
\end{proof}

\begin{proof}[Proof of \Cref{thm: PL_nwt is dtd}]
    By the definition of the adjoint operator, $d_0^*$ is the unique function such that the diagram
    \begin{center}
    \begin{tikzcd}
        (\mathcal{C}^1)^* \arrow[r, "\partial_1"]                    & (\mathcal{C}^0)^*                    \\
        \mathcal{C}^1 \arrow[r, "d_0^*"] \arrow[u, "\varphi_1"] & \mathcal{C}^0 \arrow[u, "\varphi_0"]
    \end{tikzcd}
    \end{center}
    commutes where $\varphi_i$ is defined for $i = 1, 2$ so that $\varphi_0(v)(w) = v^tP^{-1}w$ and $\varphi_1(v)(w) = v^tQ_1w$ and $\partial_1$ is defined as the dual map of $d_0$.
    Note that $\varphi_i$ is isomorphism for $i = 1, 2$.
    If we represent the function $\varphi_0$, $\varphi_1$, and $\partial_1$ as matrix forms, we get $P^{-1}$, $Q_1$, and $[d_0]^t$ respectively where $[d_0]$ is the matrix form of $d_0$.
    Therefore, the matrix form of $d_0^* d_0$ is $P[d_0]^tQ_1d_0$.
    Since $[d_0]^tQ_1[d_0] = L_{n, W}^t$, we have the claim.
\end{proof}

\begin{proof}[Proof of \Cref{thm: lp_converge}]
    If $\|f\|_\infty =0$, then the theorem is obvious.
    Assume that $\|f\|_\infty > 0$.
    Note that $\|f\|_r> 0$.
    On one hand, let $t$ be a real number with $t<\|f\|_\infty$.
    Then we have
    $$
    \mu(\{|f|>t\})^\frac{1}{p}t \leq \left(\int_{|f|>t} |f|^p\;d\mu\right)^\frac{1}{p}\leq \|f\|_p.
    $$
    If we take the limit for $p$ going to infinity, we get $t \leq \underset{p\rightarrow\infty}{\liminf}\|f\|_p$, and therefore $\|f\|_\infty\leq \underset{p\rightarrow\infty}{\liminf}\|f\|_p$ since $t$ is an any number less than $\|f\|_\infty$.

    On the other hand, let $p$ be a real number with $0<r<p$.
    Then we have
    $$\|f\|_p^p = \int |f|^p\;d\mu \leq \int|f|^r|f|^{p-r}\;d\mu \leq  \|f\|_r^r\|f\|_\infty^{p-r}.$$
    Therefore, we get
    $$
    \|f\|_p\leq \left(\|f\|_r^r\right)^\frac{1}{p}\|f\|_\infty^{1-\frac{r}{p}}.
    $$
    Since $0<\|f\|_r <\infty$, if we take the limit for $p$, then we get
    $$
    \underset{p\rightarrow\infty}{\limsup}\|f\|_p \leq \|f\|_\infty.
    $$
    From the inequalities $ \underset{p\rightarrow\infty}{\limsup}\|f\|_p \leq \|f\|_\infty\leq\underset{p\rightarrow\infty}{\liminf}\|f\|_p$, we have $\|f\|_\infty = \underset{p\rightarrow\infty}{\lim}\|f\|_p$.
\end{proof}

\begin{proof}[Proof of \Cref{thm: softmax}]
Without loss of generality, we assume that $v_i\geq 0$ for $i=1, \dots, n$, and we are going to prove
$$
s(tv) \cdot v\longrightarrow \underset{i=1,\dots, n}{\max\{v_i\}}.
$$
Let $M$ be the maximum of $\{v_i\}$ and $K$ be the number of $v_i$ which attains the maximum.
Then we have
\begin{align*}
    s(tv)\cdot v &= \sum_{i=1}^n \frac{v_i\exp{(tv_i)}}{\sum_{j=1}^n \exp{(tv_j)}} = 
    \sum_{i=1}^n \frac{v_i\exp{(t(v_i - M))}}{\sum_{j=1}^n \exp{(t(v_j- M))}}\\
    &= \sum_{v_i = M} \frac{v_i\exp{(t(v_i - M))}}{\sum_{j=1}^n \exp{(t(v_j- M))}} + \sum_{v_i < M} \frac{v_i\exp{(t(v_i - M))}}{\sum_{j=1}^n \exp{(t(v_j- M))}}\\
    &=\sum_{v_i= M}\frac{M}{\sum_{j=1}^n \exp{(t(v_j-M))}} + \sum_{v_i < M} \frac{v_i\exp{(t(v_i - M))}}{\sum_{j=1}^n \exp{(t(v_j- M))}}\\
    &=\frac{KM}{\sum_{j=1}^n \exp{(t(v_j-M))}} + \sum_{v_i < M} \frac{v_i\exp{(t(v_i - M))}}{\sum_{j=1}^n \exp{(t(v_j- M))}}.
\end{align*}
If $v_i < M$, then $\exp{(t(v_i - M))}\longrightarrow 0$ as $t\longrightarrow\infty$.
Therefore, $\sum_{j=1}^n\exp{(t(v_j-M))}$ converges to $K$, and
$$s(tv)\cdot v \longrightarrow M$$
as $t$ goes to infinity.
\end{proof}

\newpage
\section{Algorithms}

\RestyleAlgo{ruled}
\SetKwComment{Comment}{/* }{ */}

\begin{algorithm}[H]
\caption{$L^p$-circular coordinate}\label{alg: p_circular}
\KwIn{A simplicial complex $X$ with $n$ vertices and $e$ edges, a cochain $\alpha\in \mathrm{C}^1(X;\bR)\simeq \bR^e$, the coboundary matrix $d_0$, an integer $p\geq 1$, and a learning rate $\eta$}
\KwOut{$L^p$-circular coordinate, $f$}
Initialize $f = (f_1, \dots, f_n) \in \bR^n$\;
\Repeat{$f$ converges}{
$g\gets \|\alpha + d_0 f\|_p$\;
$f\gets f - \eta\left(\frac{\partial g}{\partial f}\right)$\;
}
\end{algorithm}

\begin{algorithm}[H]
\label{alg: infty_circular_p}
\caption{$L^\infty$-circular coordinate(1)}
\KwIn{A simplicial complex $X$ with $n$ vertices and $e$ edges, a cochain $\alpha\in \mathrm{C}^1(X;\bR)\simeq \bR^e$, the coboundary matrix $d_0$, approximation range $(p, p+1, \dots, q)$, the number of epochs $N$, and a learning rate $\eta$}
\KwOut{$L^\infty$-circular coordinate, $f$}
Initialize $f = (f_1, \dots, f_n) \in \bR^n$\;
$I\gets p$\;
$J\gets 0$\;
\Repeat{$J= N$}{
    \Repeat{$I = q$}{
        \Repeat{$f$ converges}{
            $g\gets \|\alpha + d_0 f\|_I$\;
            $f\gets f - \eta\left(\frac{\partial g}{\partial f}\right)$\;
            $J\gets J+1$\;
        }
        $I\gets I+1$\;
    }
    \Repeat{$J = N$}{
    $g\gets \|\alpha + d_0 f\|_\infty$\;
    $f\gets f - \eta\left(\frac{\partial g}{\partial f}\right)$\;
    $J\gets J+1$\;
    }
}
\end{algorithm}

\newpage

\begin{algorithm}[H]
\caption{$L^\infty$-circular coordinate(2)}\label{alg: infty_circular_softmax}
\KwIn{A simplicial complex $X$ with $n$ vertices and $e$ edges, a cochain $\alpha\in \mathrm{C}^1(X;\bR)\simeq \bR^e$, the coboundary matrix $d_0$, the functions $h$ and $s$ in Theorem \ref{thm: softmax}, a starting temperature $t$, the number of epochs $N$, and a learning rate $\eta$}
\KwOut{$L^\infty$-circular coordinate, $f$}
Initialize $f = (f_1, \dots, f_n) \in \bR^n$\;
$I\gets 0$\;
\Repeat{$I = N$}{
\Repeat{$f$ converges}{
$g\gets (s\circ h)(tf)\cdot h(f)$\;
$f\gets f - \eta\left(\frac{\partial g}{\partial f}\right)$\;
$I \gets I+1$\;
}
$t\gets t+1$\;
}
\end{algorithm}

\newpage

\section{Additional experimental results}\label{sec: app_exp_results}

\begin{figure}[ht]
\begin{subfigure}{.47\textwidth}
    \centering
    \includegraphics[width=\linewidth]{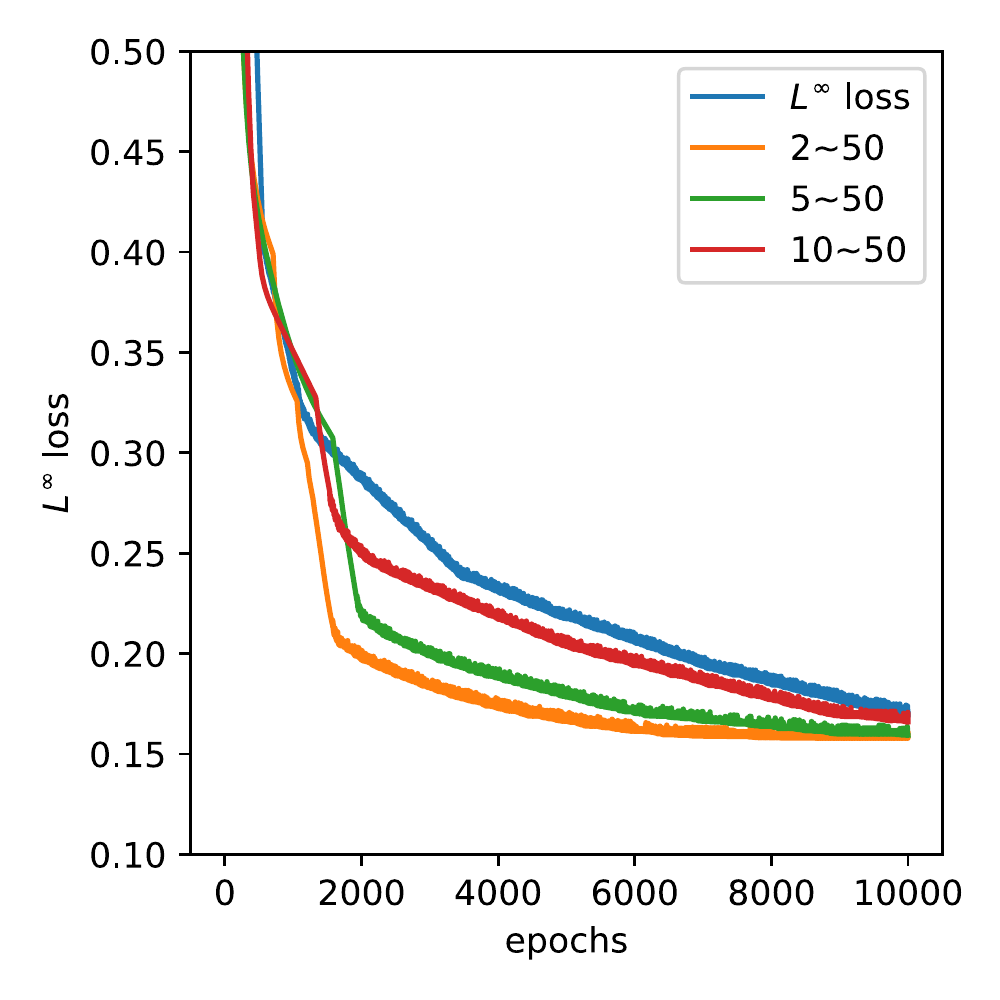}
    \caption{$\eta = 0.005, \tau = 0.0001$}
\end{subfigure}
\begin{subfigure}{.47\textwidth}
    \centering
    \includegraphics[width=\linewidth]{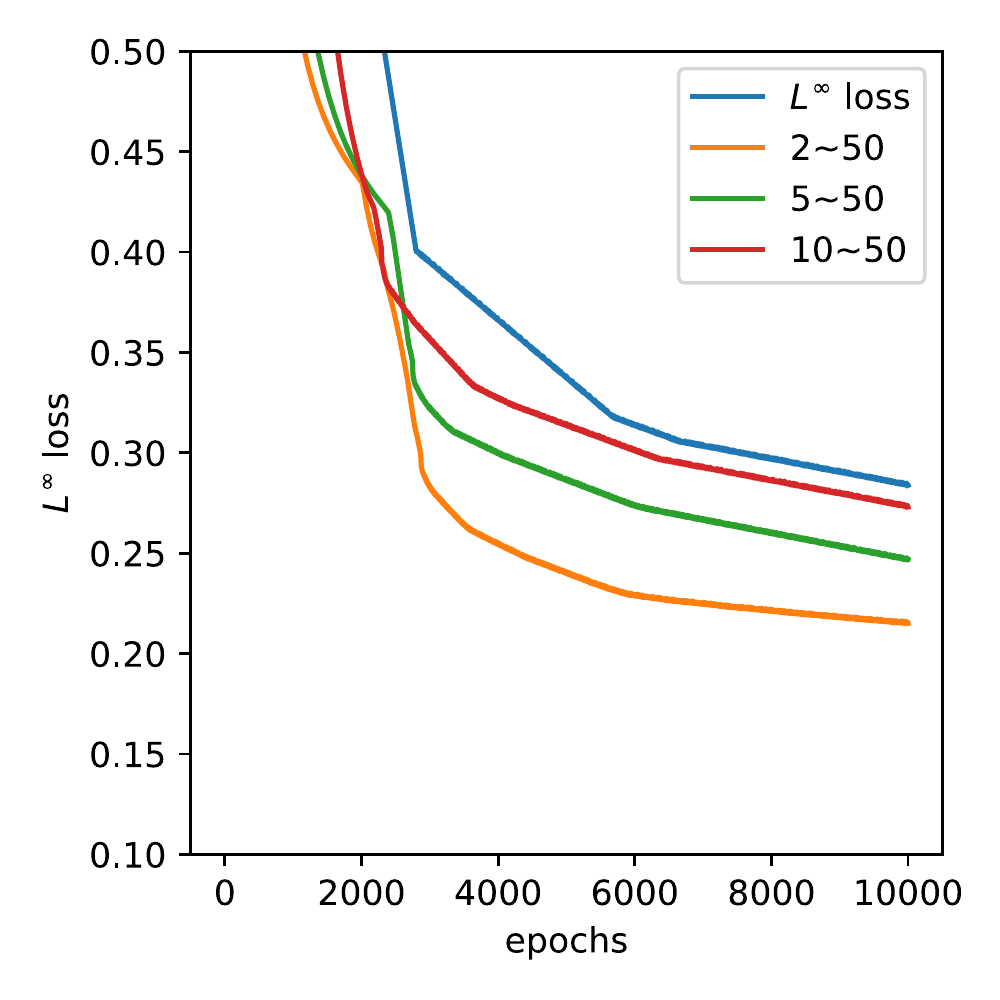}
    \caption{$\eta = 0.001, \tau = 0.0001$}
\end{subfigure}
\begin{subfigure}{.47\textwidth}
    \centering
    \includegraphics[width=\linewidth]{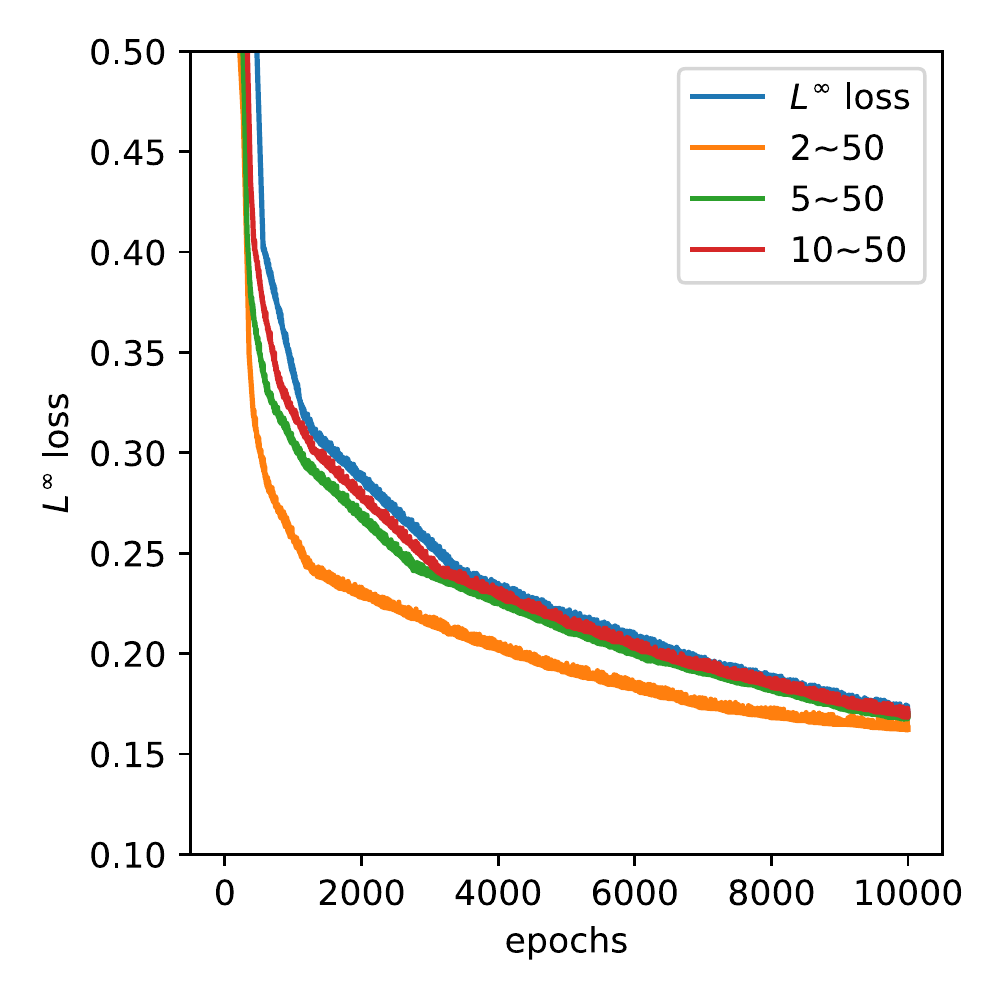}
    \caption{$\eta = 0.005, \tau = 0.001$}
\end{subfigure}
\begin{subfigure}{.47\textwidth}
    \centering
    \includegraphics[width=\linewidth]{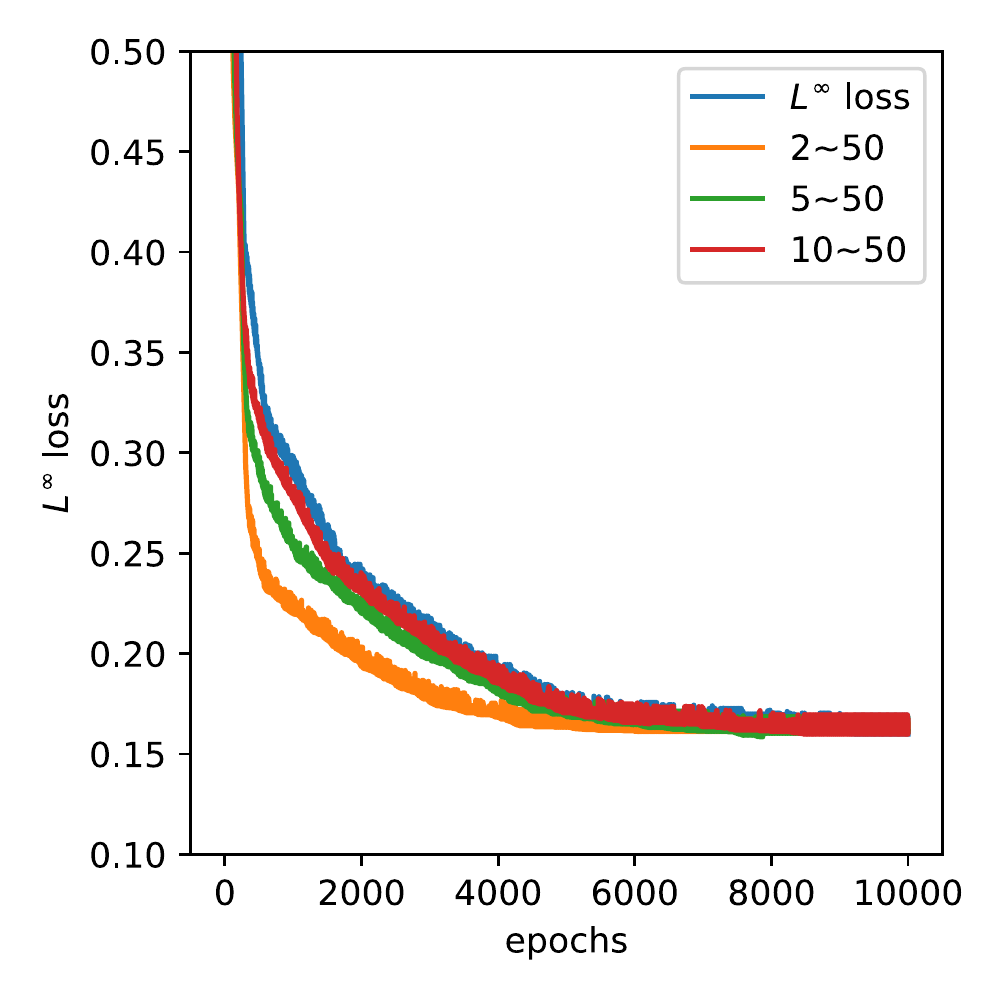}
    \caption{$\eta = 0.01, \tau = 0.001$}
\end{subfigure}
\caption{Comparisons of loss reduction in various hyperparameter settings. $\eta$ is the learning rate, and $\tau$ is the hyperparameter to determine whether a function is converged or not. The numbers on the label indicate the increasing range of $p$.}
\label{fig: lp_loss_comp1}
\end{figure}

\begin{figure}[ht]
\begin{subfigure}{.47\textwidth}
    \centering
    \includegraphics[width=\linewidth]{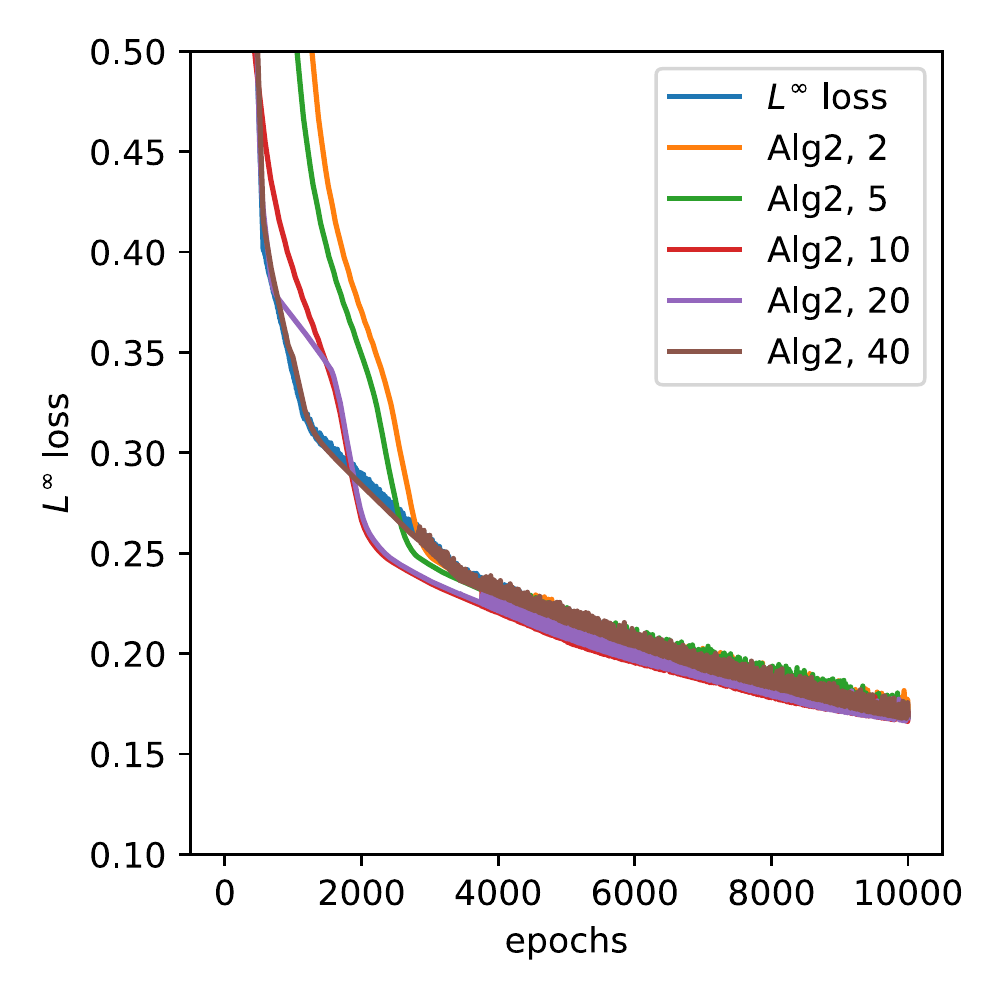}
    \caption{$\eta = 0.005, \tau = 0.0001$}
\end{subfigure}
\begin{subfigure}{.47\textwidth}
    \centering
    \includegraphics[width=\linewidth]{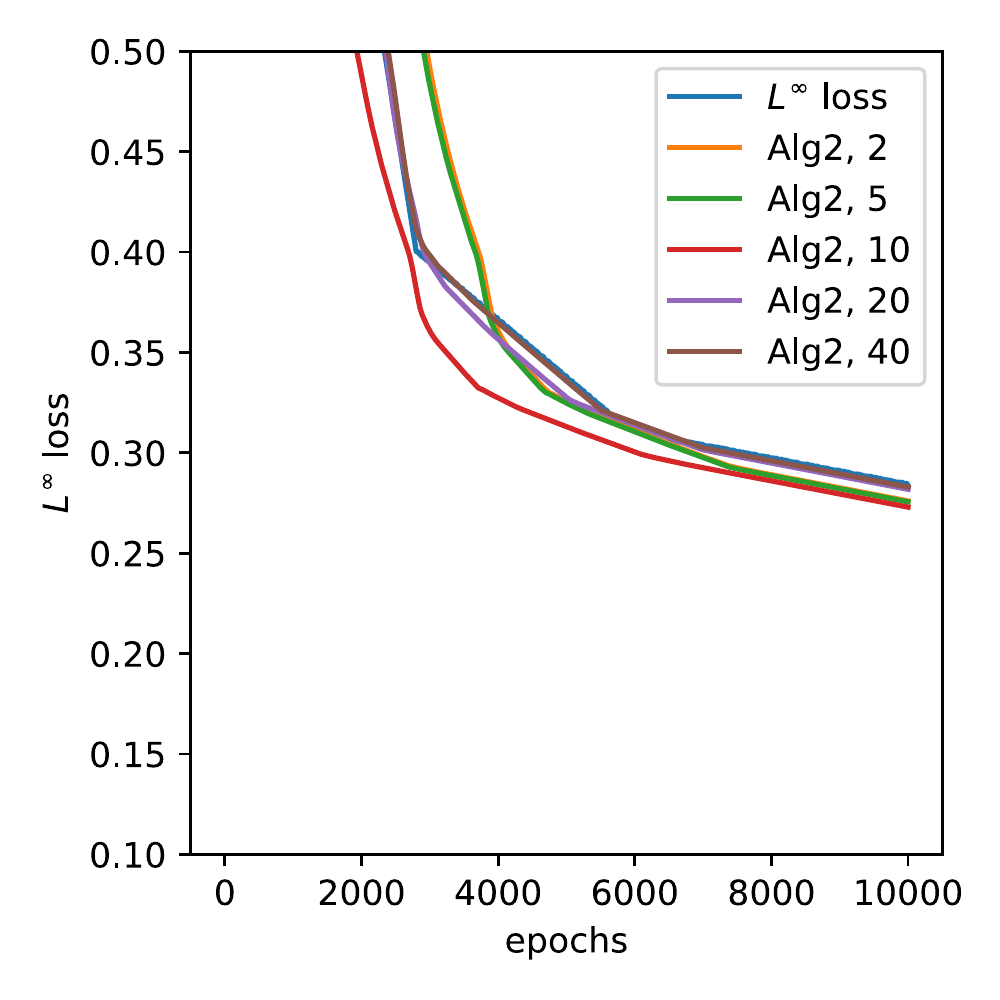}
    \caption{$\eta = 0.001, \tau = 0.0001$}
    \label{subfig: softmax_fast}
\end{subfigure}
\begin{subfigure}{.47\textwidth}
    \centering
    \includegraphics[width=\linewidth]{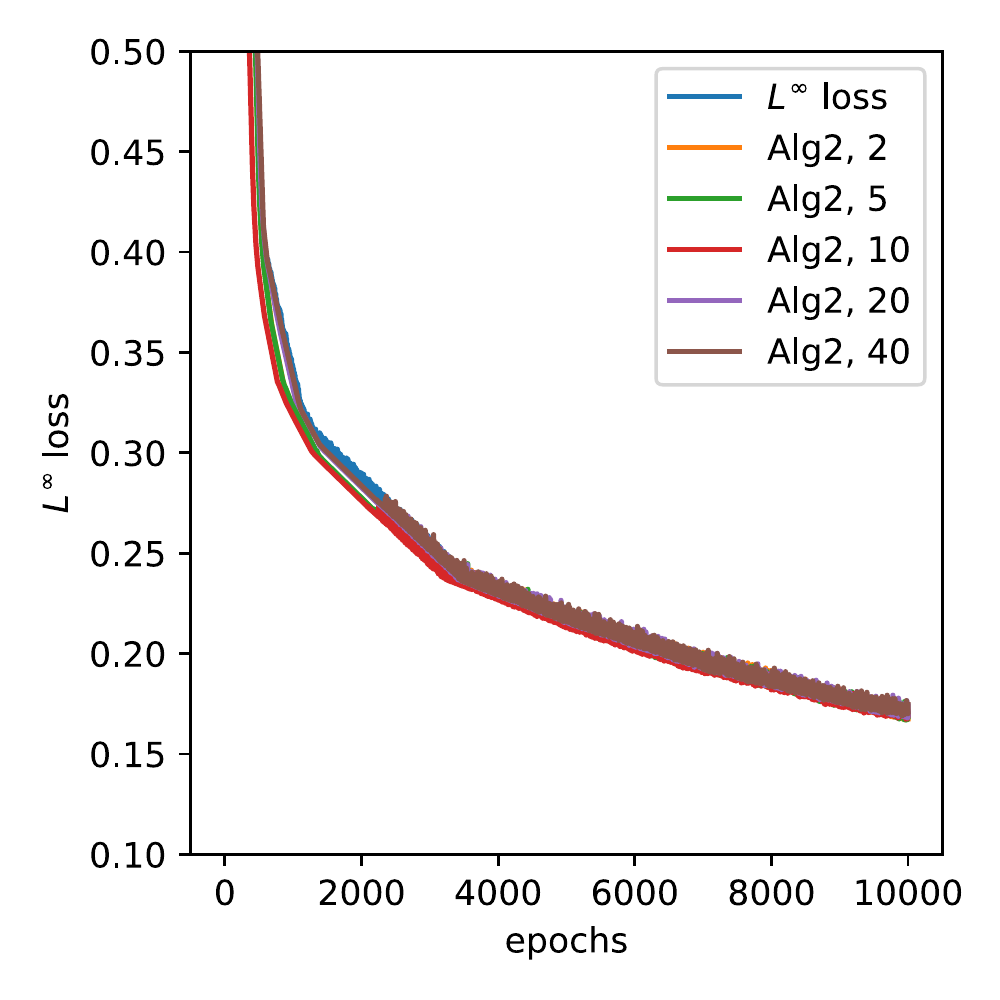}
    \caption{$\eta = 0.005, \tau = 0.001$}
\end{subfigure}
\begin{subfigure}{.47\textwidth}
    \centering
    \includegraphics[width=\linewidth]{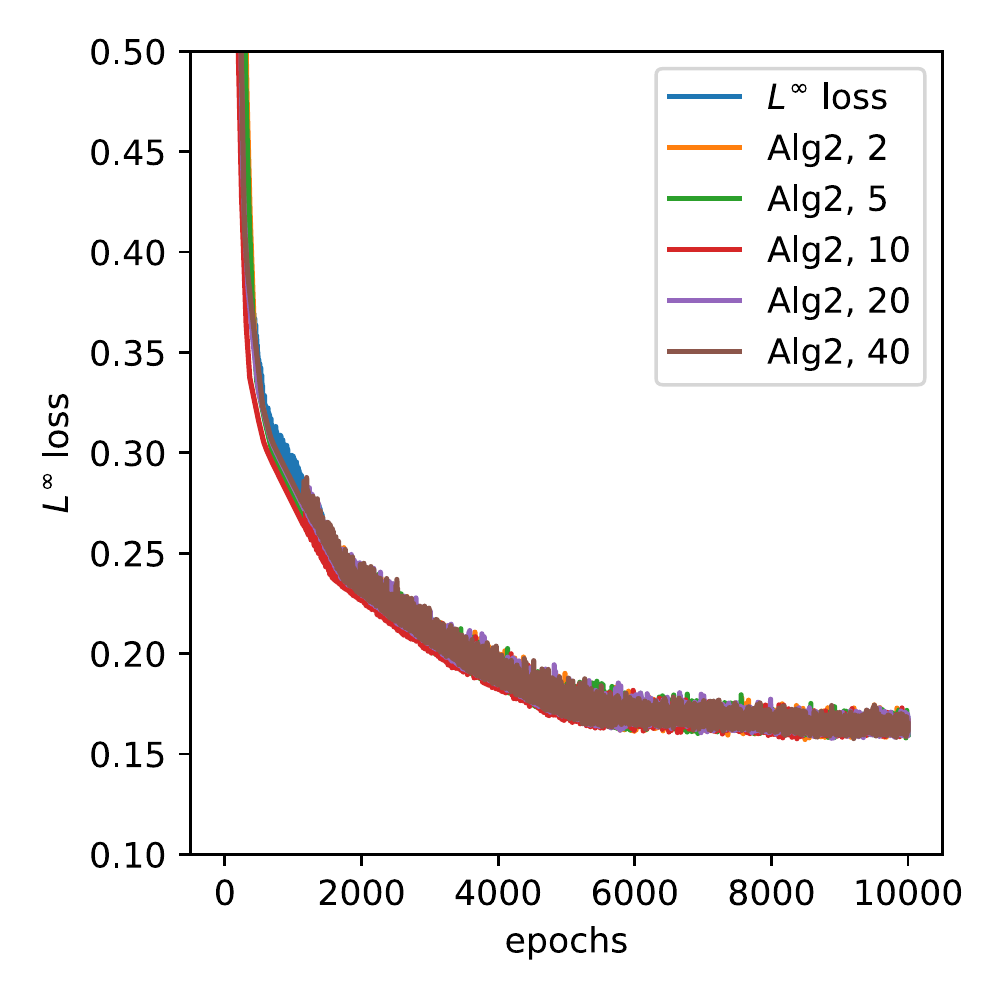}
    \caption{$\eta = 0.01, \tau = 0.001$}
\end{subfigure}
\caption{Comparisons of loss reduction in various hyperparameter settings. $\eta$ is the learning rate, and $\tau$ is the hyperparameter to determine whether a function is converged or not. The numbers on the label indicate the starting number of the temperature of the algorithm.}
\label{fig: lp_loss_comp2}
\end{figure}

\begin{figure}[ht]
\begin{subfigure}{.328\textwidth}
    \centering
    \includegraphics[width=\linewidth]{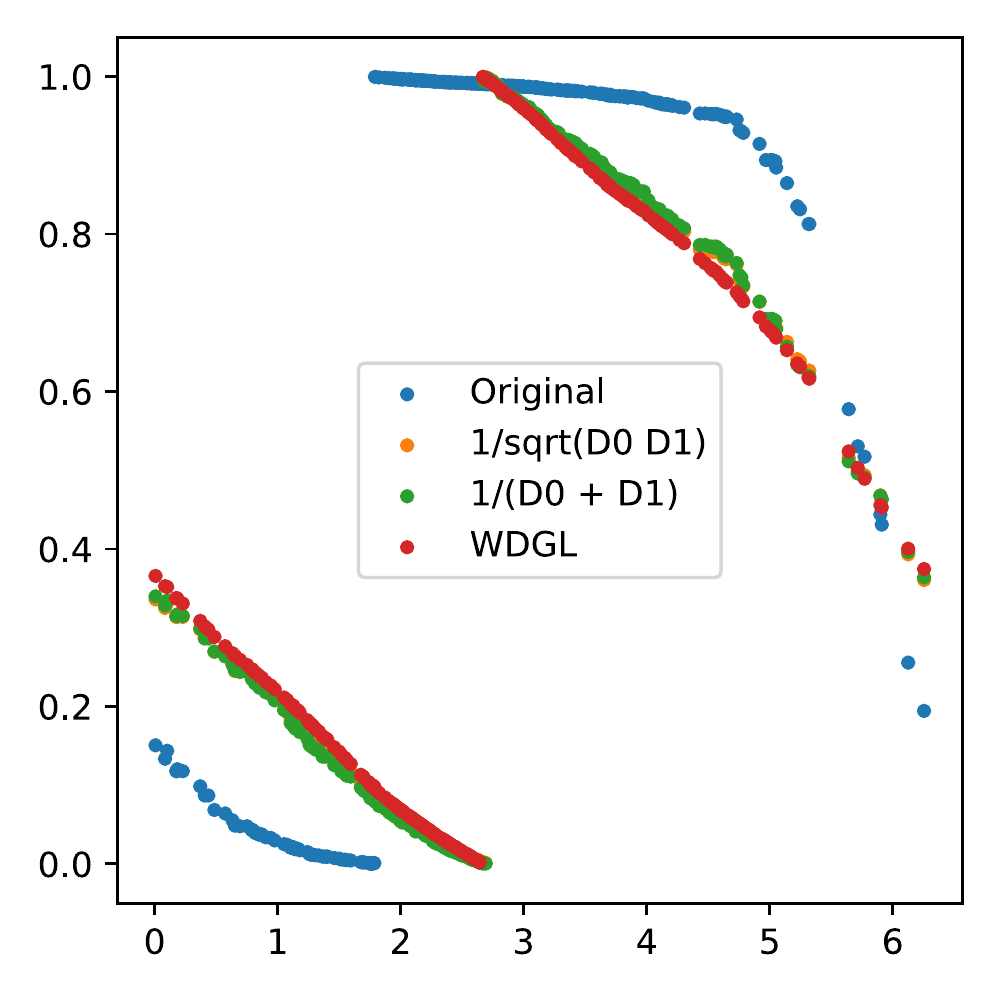}
\end{subfigure}
\begin{subfigure}{.328\textwidth}
    \centering
    \includegraphics[width=\linewidth]{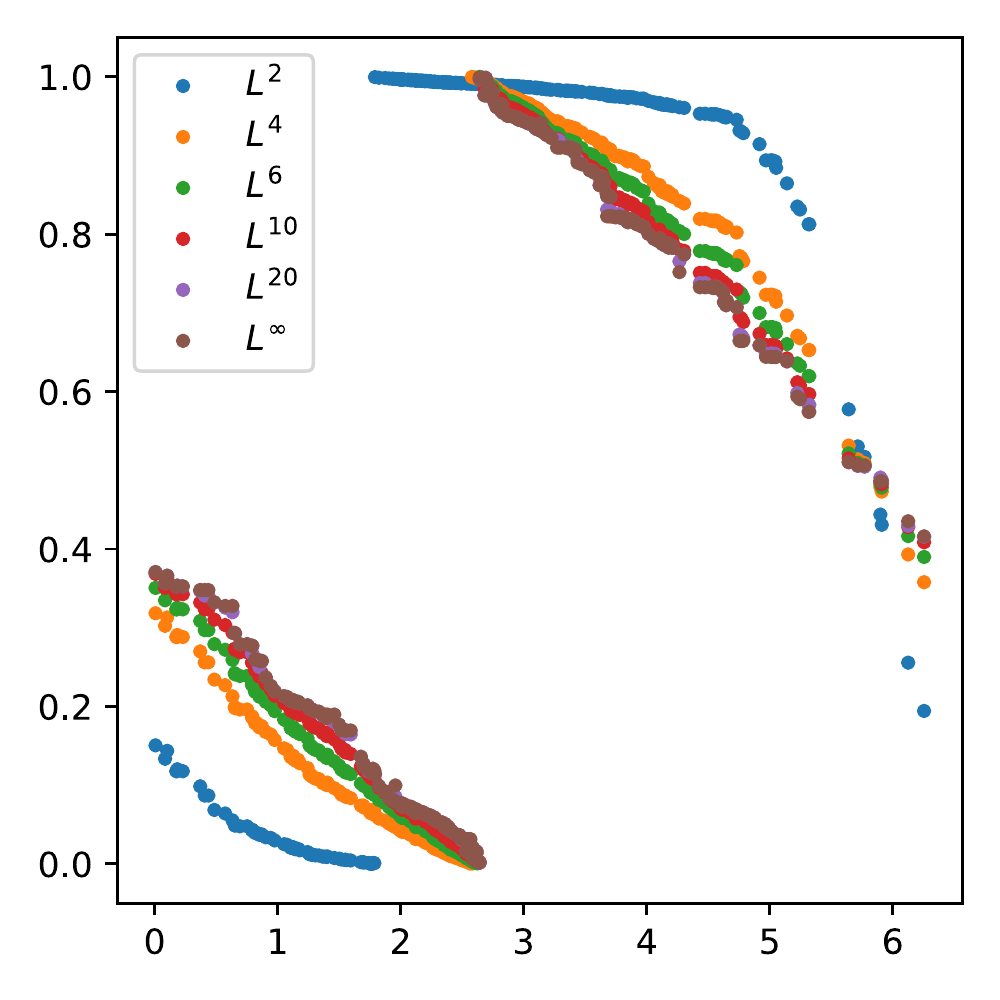}
\end{subfigure}
\caption{Results for noisy circle dataset; 
weighted circular coordinates (left), $L^p$-norm circular coordinates (right).
}
\label{fig: circle_result_scatter_app}
\end{figure}

\begin{figure}[ht]
\begin{subfigure}{.328\textwidth}
    \centering
    \includegraphics[width=\linewidth]{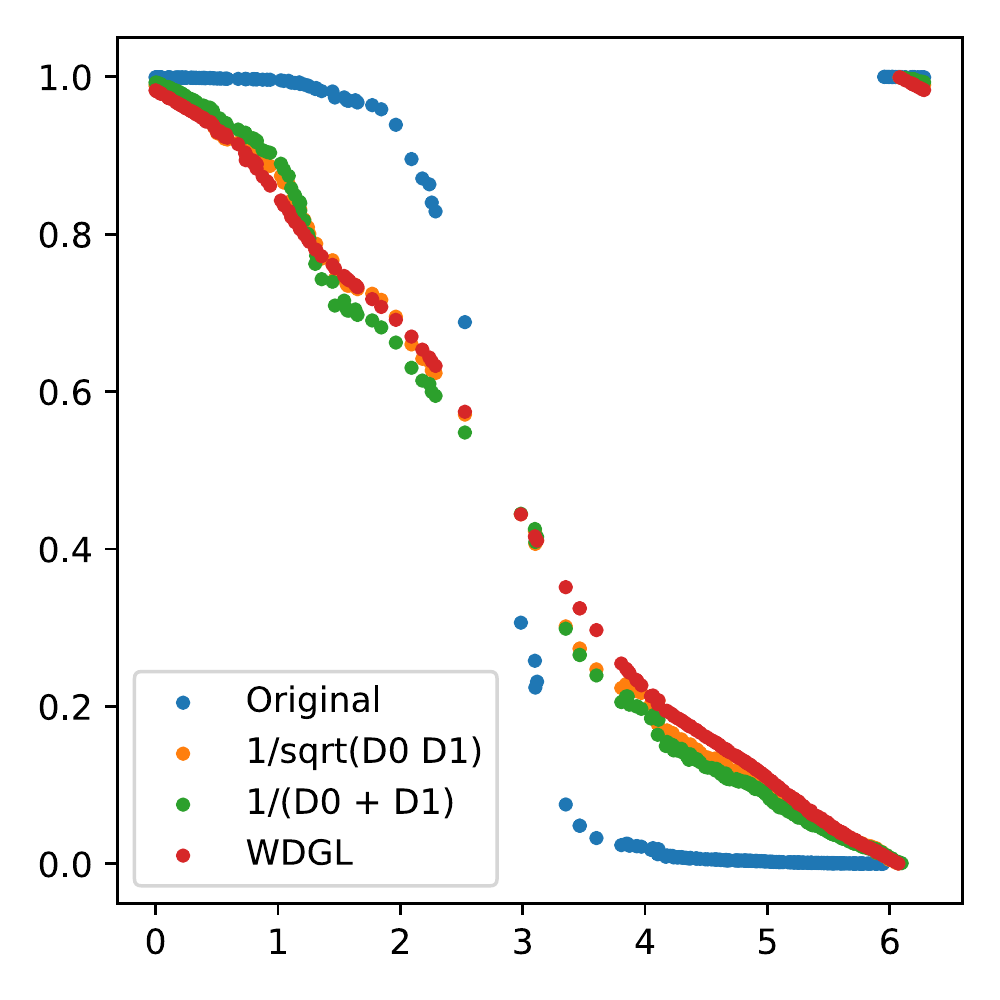}
\end{subfigure}
\begin{subfigure}{.328\textwidth}
    \centering
    \includegraphics[width=\linewidth]{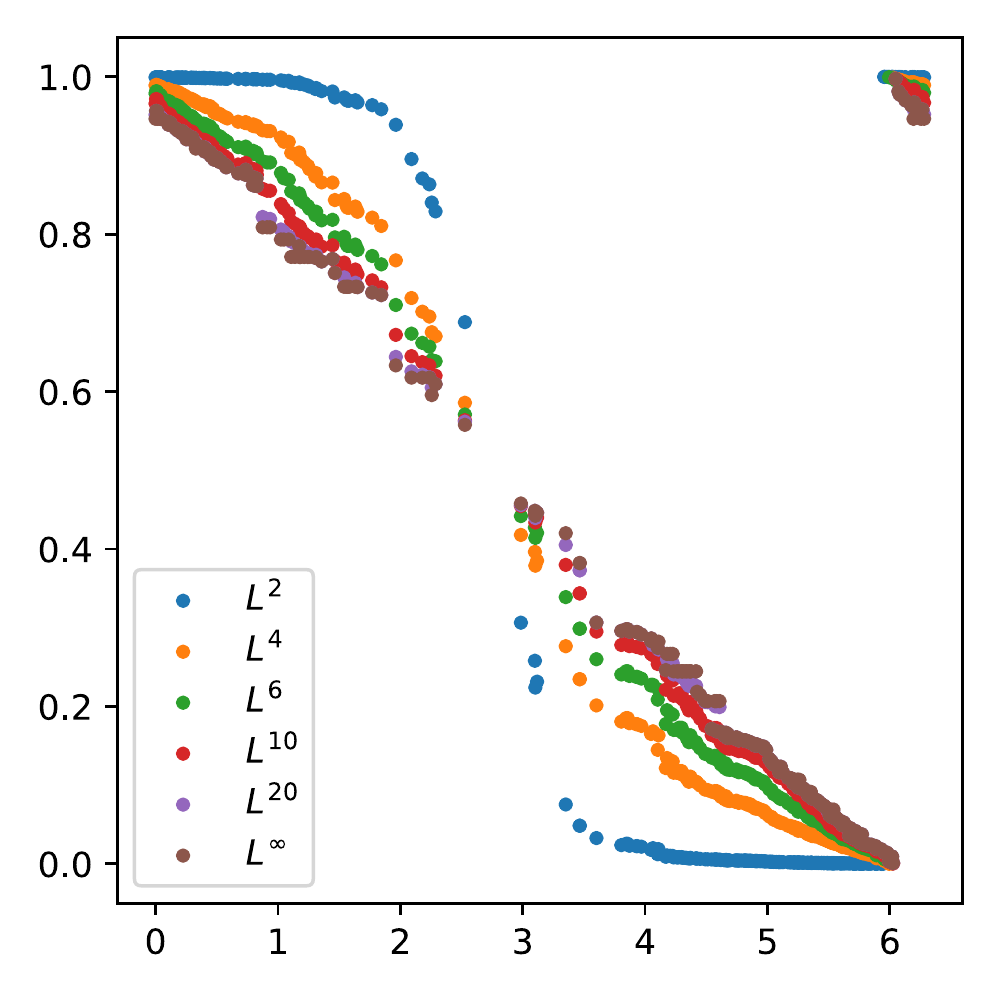}
\end{subfigure}
\\
\begin{subfigure}{.328\textwidth}
    \centering
    \includegraphics[width=\linewidth]{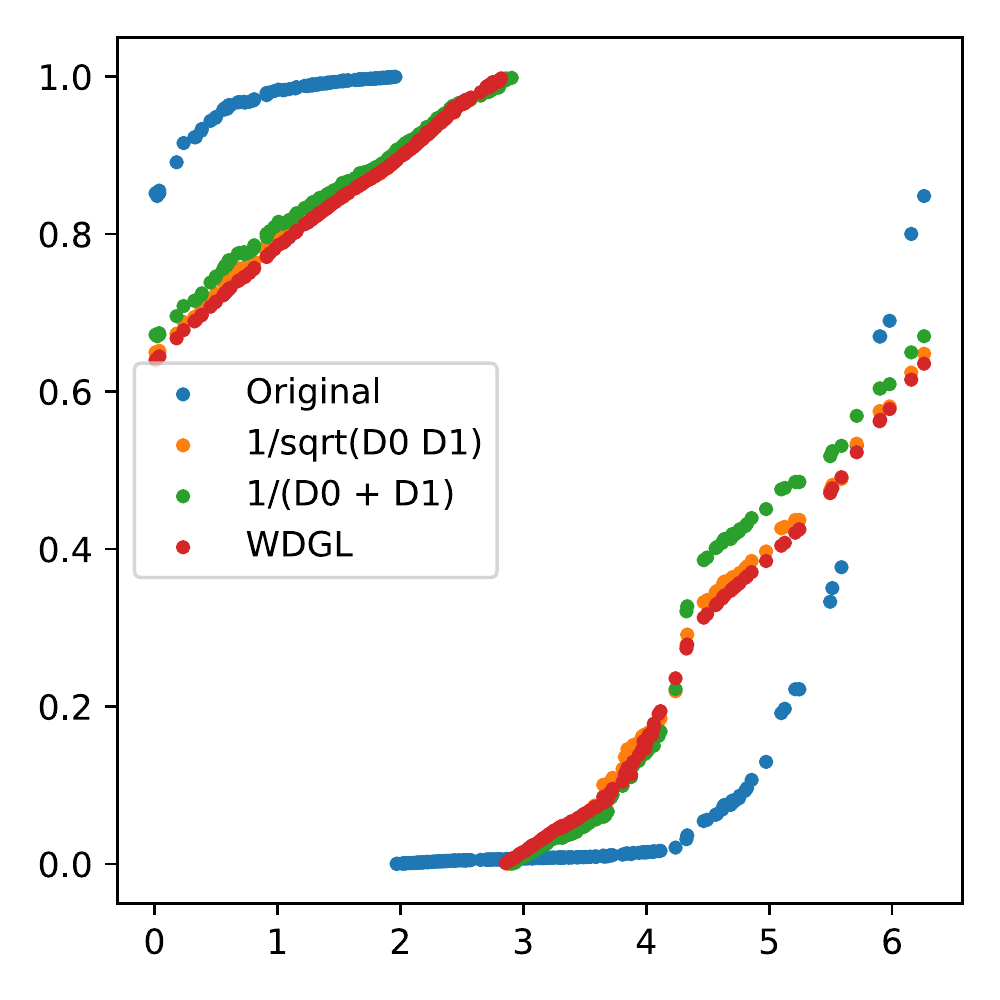}
\end{subfigure}
\begin{subfigure}{.328\textwidth}
    \centering
    \includegraphics[width=\linewidth]{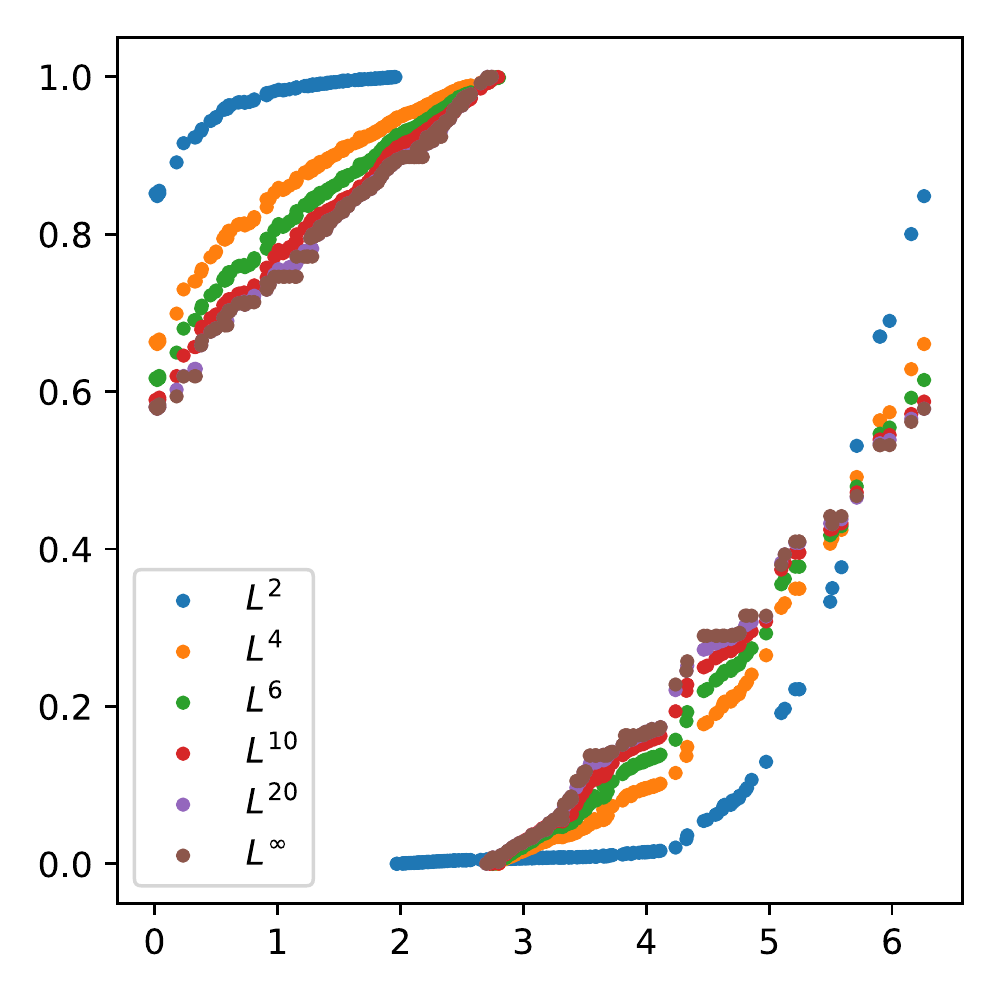}
\end{subfigure}
\caption{Results for two conjoined circles dataset; weighted circular coordinates (left column), $L^p$ circular coordinates (right column).}
\label{fig: conjoined_result_scatter_app}
\end{figure}

\begin{figure}[ht]
\begin{subfigure}{.328\textwidth}
    \centering
    \includegraphics[width=\linewidth]{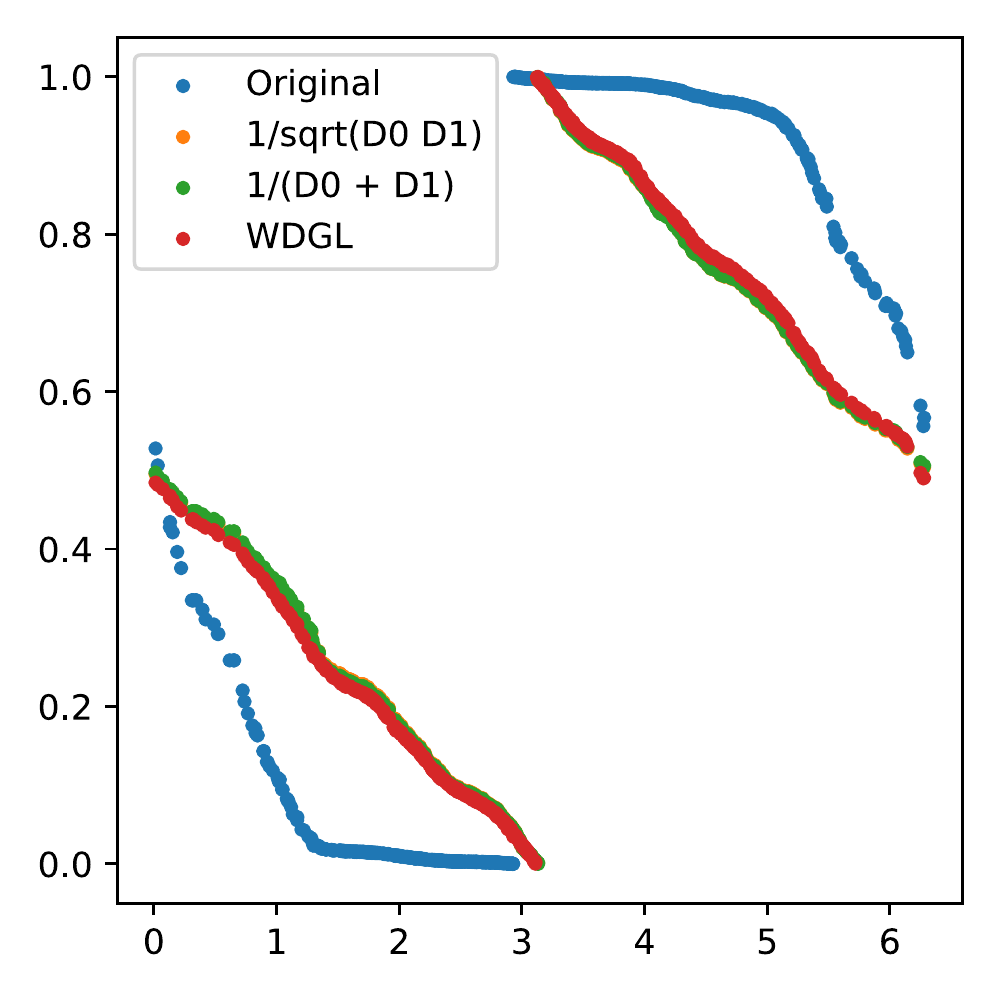}
\end{subfigure}
\begin{subfigure}{.328\textwidth}
    \centering
    \includegraphics[width=\linewidth]{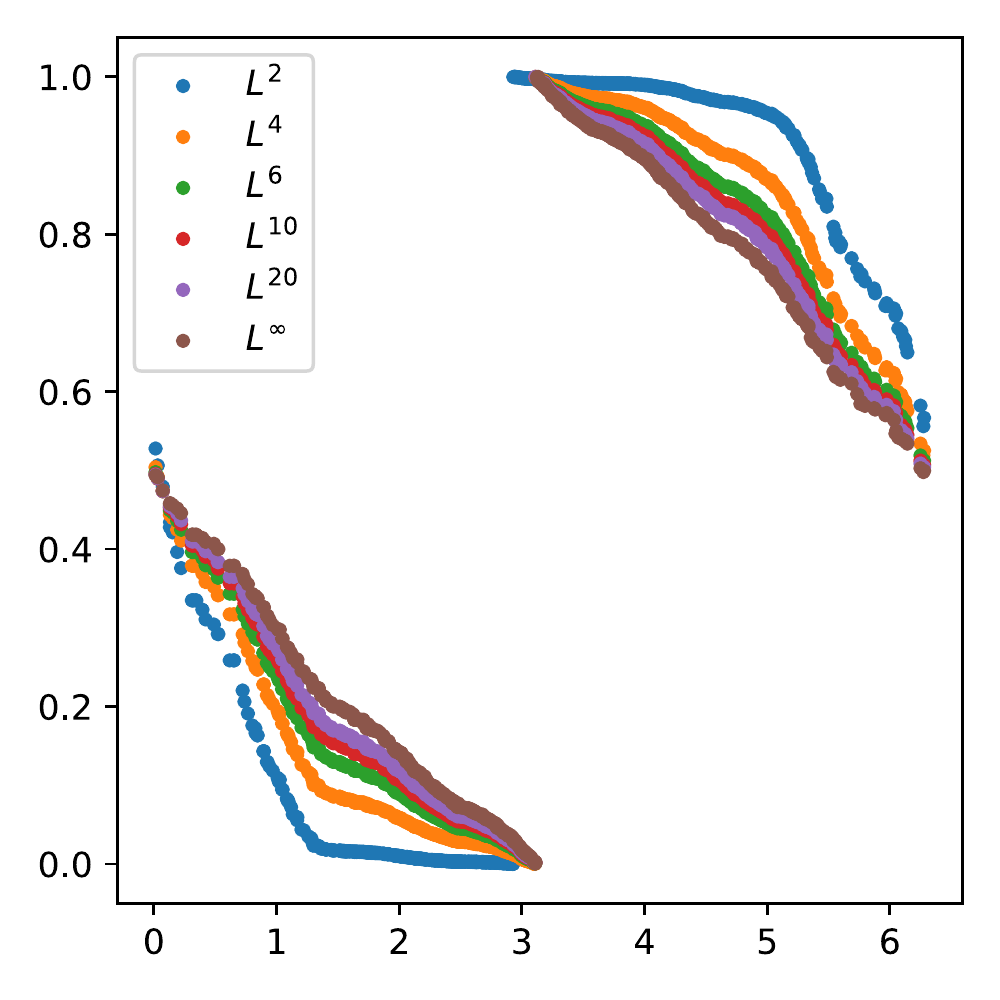}
\end{subfigure}
\caption{Results for noisy trefoil knot dataset; 
weighted circular coordinates (left), $L^p$-norm circular coordinates (right).
}
\label{fig: knot_result_scatter_app}
\end{figure}

\begin{figure}[ht]
\begin{subfigure}{.328\textwidth}
    \centering
    \includegraphics[width=\linewidth]{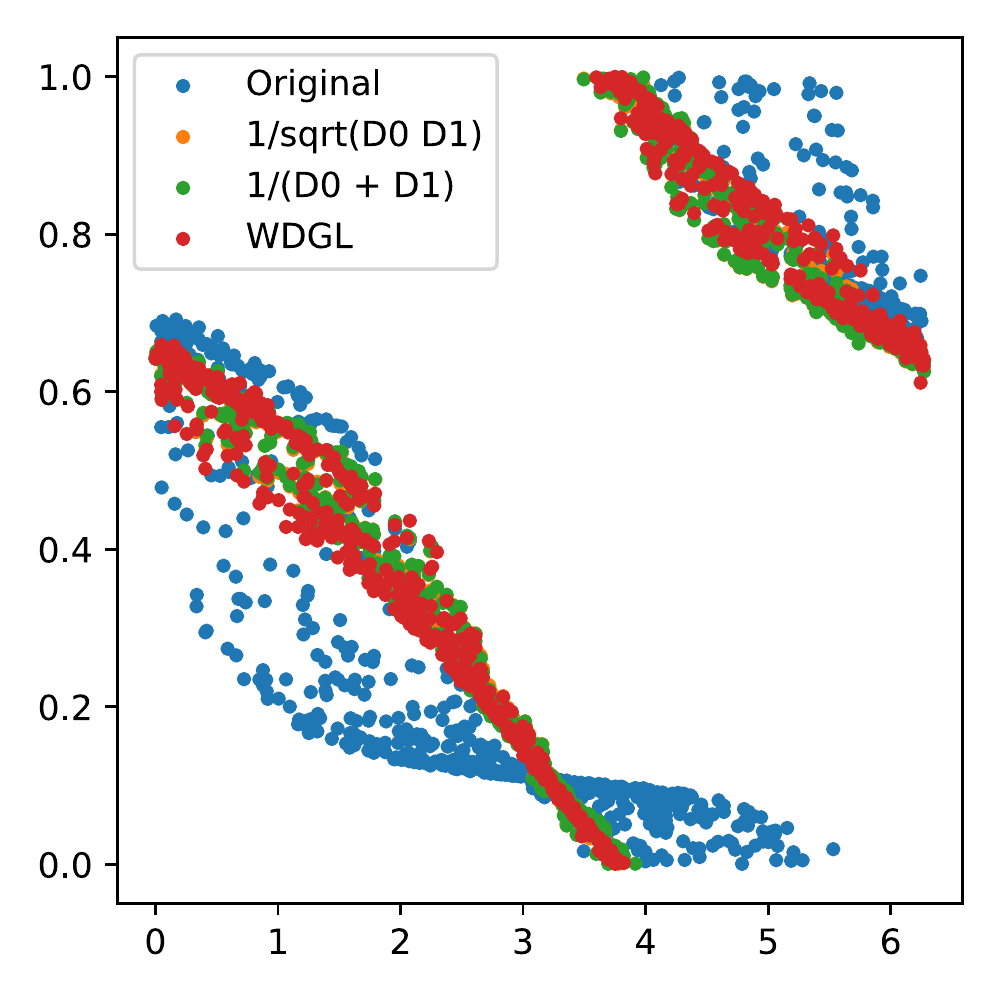}
\end{subfigure}
\begin{subfigure}{.328\textwidth}
    \centering
    \includegraphics[width=\linewidth]{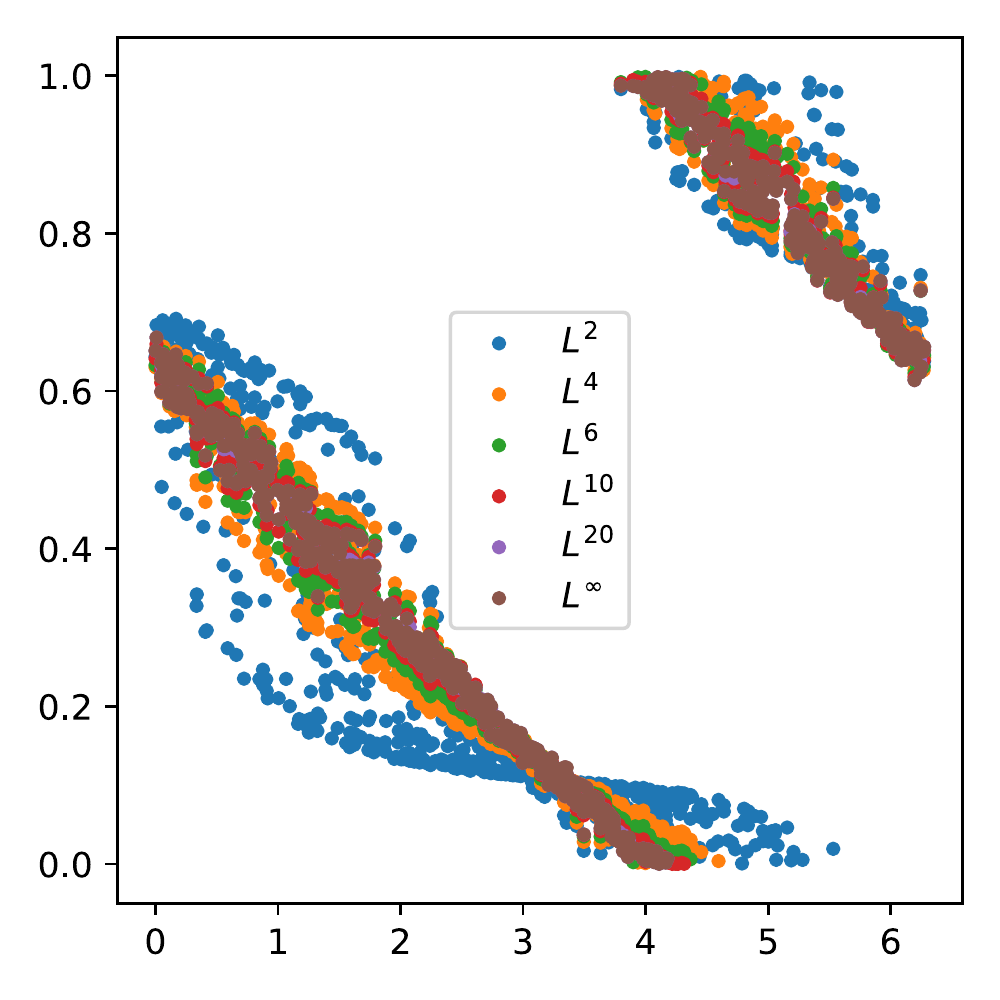}
\end{subfigure}
\\
\begin{subfigure}{.328\textwidth}
    \centering
    \includegraphics[width=\linewidth]{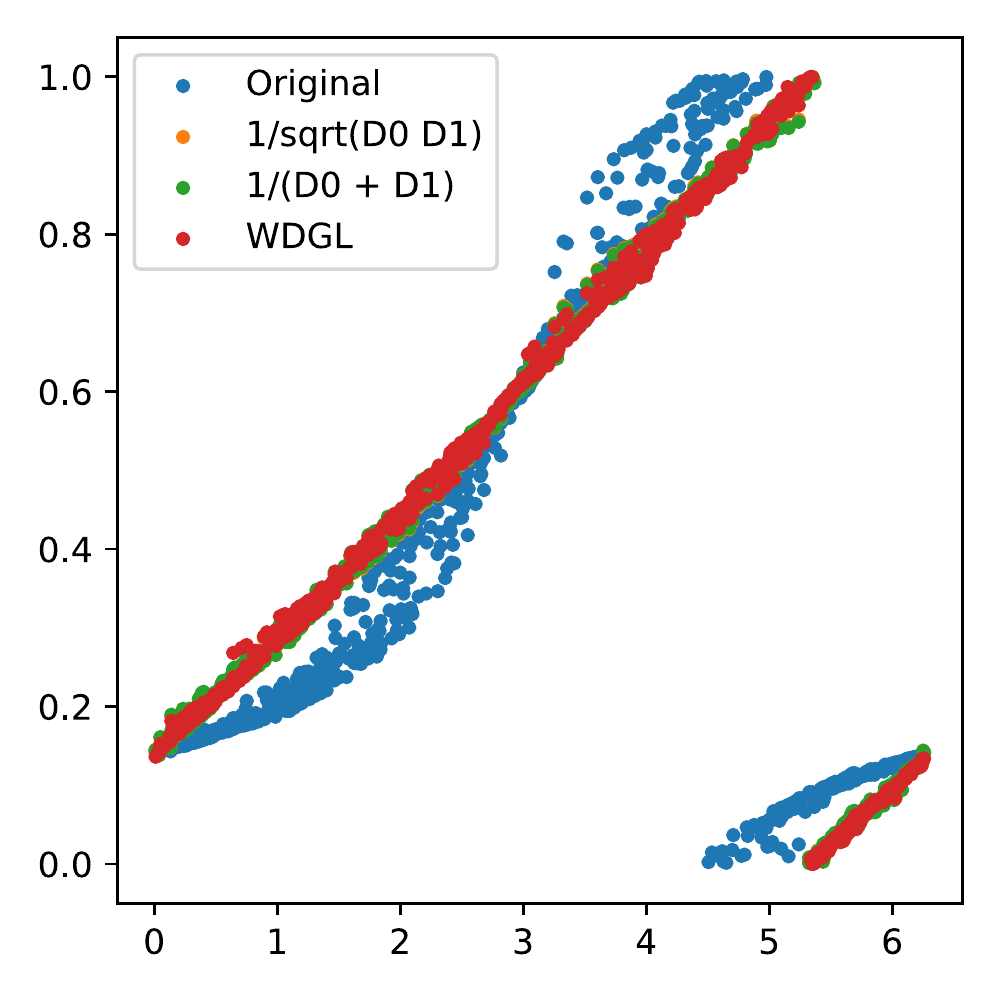}
\end{subfigure}
\begin{subfigure}{.328\textwidth}
    \centering
    \includegraphics[width=\linewidth]{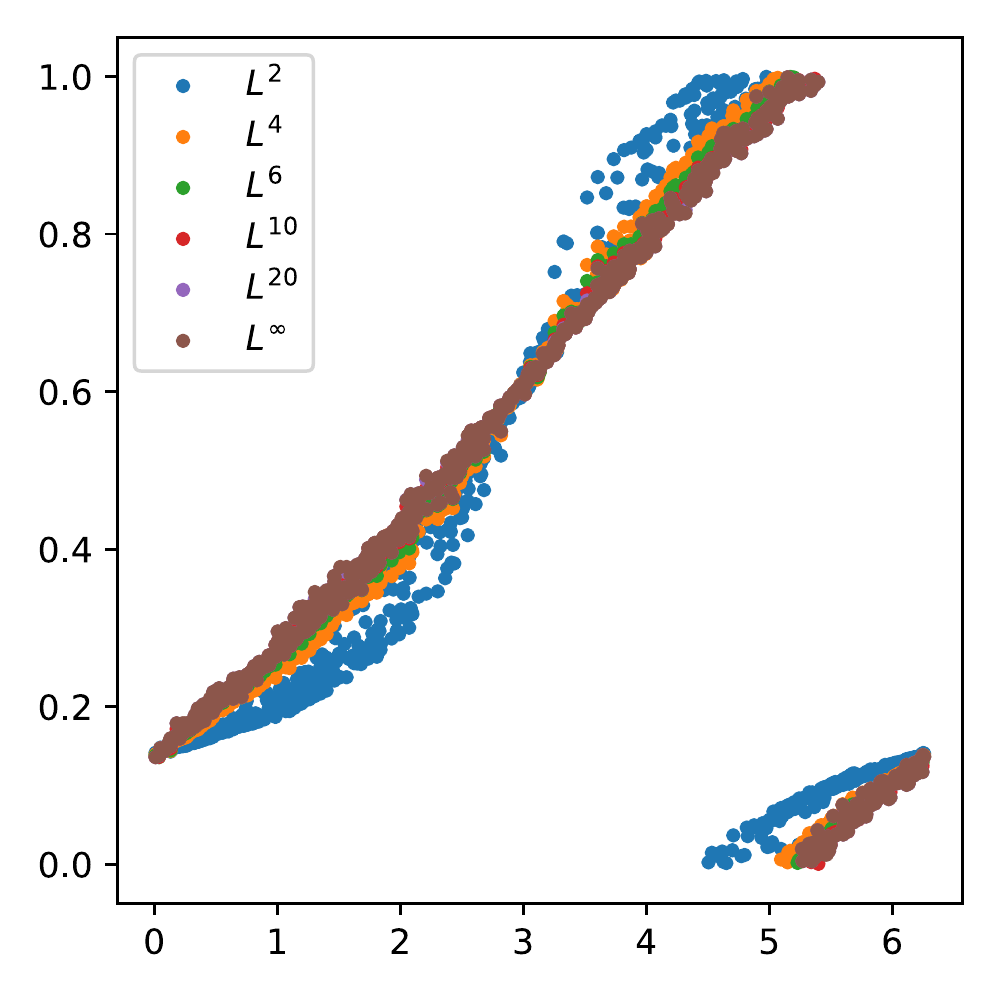}
\end{subfigure}
\caption{Results for torus dataset; weighted circular coordinates (left column), $L^p$ circular coordinates (right column).}
\label{fig: torus_result_scatter_app}
\end{figure}

\clearpage

\end{document}